%% file: main.tex
\documentclass[11pt]{article}

\input{packages}

\input{commands}

\def\softd{{\leavevmode\setbox1=\hbox{d}%
\hbox to 1.05\wd1{d\kern-0.4ex{\char039}\hss}}}

\begin{document}

%
%
	
\begin{center}	
    \LARGE{Convergent numerical schemes for the viscoelastic Giesekus model in two dimensions}
\end{center}
\bigskip
\begin{center}	
	\normalsize{Endre S{\"u}li}\\[1mm]
    \textit{Mathematical Institute,}\\[1mm]
	\textit{University of Oxford, United Kingdom}\\[1mm]
	\texttt{endre.suli@maths.ox.ac.uk}
\end{center}
\medskip
\begin{center}	
	\normalsize{Dennis Trautwein}\\[1mm]
    \textit{Department of Mathematics,}\\[1mm]
	\textit{University of Regensburg, Germany}\\[1mm]
	\texttt{dennis.trautwein@ur.de}
\end{center}
\bigskip

\begin{abstract}
\footnotesize
In this work, we develop a class of stable and convergent numerical methods for the approximate solution of the viscoelastic Giesekus model in two space dimensions. The model couples the incompressible Navier--Stokes equations with an evolution equation for an additional stress tensor accounting for elastic effects. This coupled evolution equation is stated here in terms of the elastic deformation gradient and models transport and nonlinear relaxation effects. In the existing literature, numerical schemes for such models often suffer from accuracy limitations and convergence problems, usually due to the lack of rigorous existence results or inherent limitations of the discretization. Therefore, our main goal is to prove the (subsequence) convergence of the proposed numerical method to a large-data global weak solution in two dimensions, without relying on cut-offs or additional regularization. This also provides an alternative proof of the recent existence result by Bul\'{\i}\v{c}ek et al.~(Nonlinearity, 2022). Finally, we verify the practicality of the proposed method through numerical experiments, including convergence studies and typical benchmark problems.
%
\medskip

\medskip
	
\end{abstract}

\input{results}

\vspace{-\baselineskip}

\section*{Acknowledgements}
This work was initiated during a visit of D.T.~to the University of Oxford in 2025. D.T.~is grateful for the hospitality and support received there. The work of D.T.~was partially supported by the Deutsche For\-schungs\-ge\-mein\-schaft (DFG, German Research Foundation) via the Gra\-duier\-ten\-kolleg 2339 IntComSin (Project-ID 321821685), and by the Swedish Research Council (grant no.~2021-06594) during residence at the Institut Mittag-Leffler, Djursholm, Sweden, in 2025.


\footnotesize\setlength{\parskip}{0cm}

\addcontentsline{toc}{section}{References}
\bibliography{references}

\end{document}

%% file: packages.tex
\usepackage[left=1.2in, right=1.2in, bottom=1in, top=1in, footskip=0.3in, 
]{geometry}

\usepackage{xcolor}
\definecolor{LinkColor}{rgb}{0,0,1}
\definecolor{LinkColor2}{rgb}{1,0,0}
\definecolor{lbcolor}{rgb}{0.85,0.85,0.85}
\definecolor{FrameColor}{rgb}{0.85,0.85,0.85}

\usepackage{todonotes}

\usepackage[ngerman, english]{babel}
\usepackage[utf8]{inputenc}
\usepackage{enumerate}
\usepackage[normalem]{ulem}
\usepackage{setspace}

\usepackage{multicol}
\usepackage[babel,french=guillemets,german=swiss]{csquotes}

\usepackage{array}
\usepackage{booktabs}
\usepackage{framed}
\usepackage{rotating}
\usepackage{longtable}
\usepackage{multirow}
\usepackage{tabularx}
\newcolumntype{L}[1]{>{\raggedright\arraybackslash}p{#1}} 
\newcolumntype{C}[1]{>{\centering\arraybackslash}p{#1}} 
\newcolumntype{R}[1]{>{\raggedleft\arraybackslash}p{#1}} 

\usepackage{graphicx,import}
\usepackage{caption}

\usepackage{amssymb}
\usepackage{MnSymbol}   
\usepackage{mathtools}  
\usepackage{amsmath}
\usepackage{dsfont}
\usepackage{nicefrac}
\usepackage{empheq}
\allowdisplaybreaks  

\usepackage[%
pdftitle={Titel},%
pdfauthor={Autor},%
pdfcreator={LaTeX, LaTeX with hyperref and KOMA-Script},
pdfsubject={Betreff}, 
pdfkeywords={Keywords}
]{hyperref} 

\hypersetup{%
	colorlinks=true,
	linkcolor=LinkColor,%
	anchorcolor=LinkColor,%
	citecolor=LinkColor2,%
	filecolor=LinkColor,%
	menucolor=LinkColor,%
	urlcolor=LinkColor,%
}

\usepackage[savemem]{listings}
\usepackage{paralist}
{\begin{list}{$\diamondsuit$}{}}%
	{\end{list}}

\usepackage{amsthm}
\usepackage{upgreek}

\newtheoremstyle{tstyle}
{15pt}	
{5pt}	
{\itshape}	
{}	
{\bfseries}	
{.}	
{0.5em}	
{}	

\theoremstyle{tstyle}

\newtheorem{theorem}{Theorem}[section]
\newtheorem{lemma}[theorem]{Lemma}

\newtheorem{definition}[theorem]{Definition}
\newtheorem{remark}[theorem]{Remark}

\newtheorem*{theorem*}{Theorem}

\newtheoremstyle{cstyle}
{15pt}	
{5pt}	
{}	
{}	
{\bfseries}	
{}	
{0.2222em}	
{}	

\theoremstyle{cstyle}

\makeatletter
\g@addto@macro{\thm@space@setup}{\thm@headpunct{}}
\renewenvironment{proof}[1][\proofname.]{\par
	\pushQED{\qed}%
	\normalfont \topsep6\p@\@plus6\p@\relax
	\trivlist
	\item[\hskip\labelsep
	\bfseries
	#1\@addpunct{\,}]\ignorespaces
}{%
	\popQED\endtrivlist\@endpefalse
}
\makeatother

\makeatletter
\g@addto@macro{\thm@space@setup}{\thm@headpunct{}}
\newenvironment{sketch-proof}[1][Sketch of the proof]{\par
	\pushQED{\qed}%
	\normalfont \topsep6\p@\@plus6\p@\relax
	\trivlist
	\item[\hskip\labelsep
	\bfseries
	#1\@addpunct{\,}]\ignorespaces
}{%
	\popQED\endtrivlist\@endpefalse
}
\makeatother

\newcounter{subeq}

\usepackage{fancyhdr}
\usepackage{tikz}
\usepackage{siunitx}

\usepackage{pdflscape}
\usepackage{subfig}
\usepackage{float}
\usepackage{hyperref}
\usepackage{pgfplots}
\pgfplotsset{compat=1.16} 
\usepackage{algorithmic}

\usepackage[section]{placeins}         

\numberwithin{equation}{section}
\usepackage[only,llbracket,rrbracket]{stmaryrd}



%% file: commands.tex
\makeatletter
\newcommand{\mylabel}[2]{%
  \protected@write \@auxout {}{\string \newlabel {#1}{{#2}{\thepage}{#2}{#1}{}}}%
  \hypertarget{#1}{#2}%
}
\newcommand{\mylabelHIDE}[2]{%
  \protected@write \@auxout {}{\string \newlabel {#1}{{#2}{\thepage}{#2}{#1}{}}}%
  \hypertarget{#1}{}%
}
\makeatother

\newcommand{\bc}{\ensuremath{\mathbf{c}}}

\newcommand{\bn}{\ensuremath{\mathbf{n}}}
\newcommand{\bu}{\ensuremath{\mathbf{u}}}
\newcommand{\bv}{\ensuremath{\mathbf{v}}}
\newcommand{\bw}{\ensuremath{\mathbf{w}}}
\newcommand{\bx}{\ensuremath{\mathbf{x}}}

\newcommand{\bz}{\ensuremath{\mathbf{z}}}

\newcommand{\bbA}{\ensuremath{\mathbb{A}}}
\newcommand{\bbB}{\ensuremath{\mathbb{B}}}

\newcommand{\bbF}{\ensuremath{\mathbb{F}}}
\newcommand{\bbG}{\ensuremath{\mathbb{G}}}
\newcommand{\bbI}{\ensuremath{\mathbb{I}}}
\newcommand{\bbN}{\ensuremath{\mathbb{N}}}

\newcommand{\bbR}{\ensuremath{\mathbb{R}}}
\newcommand{\bbS}{\ensuremath{\mathbb{S}}}

\newcommand{\mycal}[1]{\ensuremath{\mathcal{#1}}}
\newcommand{\calA}{\mycal{A}}
\newcommand{\calB}{\mycal{B}}

\newcommand{\calP}{\ensuremath{\mathcal{P}}}
\newcommand{\calS}{\ensuremath{\mathcal{S}}}
\newcommand{\calT}{\ensuremath{\mathcal{T}}}

\DeclareMathOperator{\Div}{div}
\DeclareMathOperator{\diam}{diam}

\DeclareMathOperator{\trace}{Tr}

\newcommand{\abs}[1]{\left\lvert{#1}\right\rvert}
\newcommand{\norm}[1]{\|{#1}\|}
\newcommand{\nnorm}[1]{\left\|{#1}\right\|}
\newcommand{\dualp}[2]{\left<{#1\, ,\, #2}\right>}
\newcommand{\skp}[2]{\left({#1\, ,\, #2}\right)}

\newcommand{\dv}[1]{\,{\mathrm d}#1}

\newcommand{\dt}{\dv{t}}


%% file: results.tex
\section{Introduction} \label{sec:introduction}

Viscoelastic fluids appear in many industrial and natural settings. Common examples range from biological fluids and geological flows to synthetic polymers. Unlike Newtonian fluids, these materials simultaneously store and dissipate energy. This dual nature leads to distinctive non-Newtonian phenomena, such as the Weissenberg effect (rod-climbing), die swell, and drag reduction \cite{bird_armstrong_hassager_1987}.
Mathematically, these elastic effects are modeled by introducing an extra stress tensor into the momentum equation. For polymeric fluids, this is typically formulated using a conformation tensor $\mathbb{B}$. This tensor acts as a macroscopic variable that describes how the underlying polymer chains stretch and orient during flow \cite{masmoudi_2018_handbook}.
At the microscopic level, these chains are often modeled as a collection of massless beads connected by elastic springs. While kinetic theory can describe the associated dynamics using Fokker--Planck equations, the high-dimensionality of the Fokker--Planck equation makes its numerical approximation computationally expensive \cite{beddrich_suli_wohlmuth_2024, bird_curtiss_armstrong_hassager_1987}. Consequently, closed constitutive laws at the macroscopic level are essential.
The specific choice of the spring law determines the macroscopic model \cite{barrett_suli_2011_partI, barrett_suli_2012_partII}. Assuming linear spring forces (Hooke's law) leads to the standard Oldroyd-B model \cite{Oldroyd_1950}. However, this model predicts infinite elongational viscosity at finite strain rates. To overcome this physical deficiency, the Giesekus model accounts for anisotropic drag between polymer chains, resulting in a nonlinear constitutive relaxation \cite{giesekus_1982}. 
In this work, we consider the resulting macroscopic system. While the Giesekus model is often expressed using the conformation tensor $\mathbb{B}$, we formulate the governing equations using the elastic deformation gradient $\mathbb{F}$. This approach allows us to exploit additional geometric structure, noting that the conformation tensor can be recovered via $\mathbb{B} = \mathbb{F}\mathbb{F}^\top$.

\paragraph{Governing equations.}
Let $T>0$ and $\Omega\subset\bbR^d$, $d\in\{2,3\}$, be a bounded domain with Lipschitz boundary $\partial\Omega$. 
We seek functions 
\begin{align*}
    (\bv, p, \bbF)\, \colon\, \Omega\times(0,T) \to \bbR^d \times \bbR \times \bbR^{d\times d},
\end{align*}
where $\bv$ denotes the fluid velocity, $p$ the pressure, and $\bbF$ the elastic deformation gradient. These functions satisfy the following system in $\Omega \times (0,T)$:
\begin{subequations}
\label{eq:system}
\begin{align}
    \label{eq:v}
    \rho \partial_t \bv + \rho (\bv\cdot\nabla) \bv
    + \nabla p - \nu \Delta \bv &= \Div( \mu \bbF \bbF^\top),
    \\
    \label{eq:div}
    \Div(\bv) &= 0, 
    \\
    \label{eq:F}
    \partial_t \bbF + (\bv\cdot\nabla) \bbF + \frac{\mu}{2\lambda} (\bbF \bbF^\top \bbF - \bbF) 
    &= (\nabla\bv) \bbF ,
\end{align}
\end{subequations}
and we require
\begin{align}
    \label{eq:detF}
    \det(\bbF)&>0
    \qquad \text{in } \Omega\times(0,T).
\end{align}
The system is complemented by the homogeneous Dirichlet boundary condition
\begin{align}
    \label{eq:u_bc}
    \bv=\mathbf{0} 
    \qquad \text{on } \partial\Omega\times(0,T),
\end{align}
and the following initial conditions:
\begin{align}
    \label{eq:system_init}
    \bv(\cdot,0) = \bv_0, \qquad \bbF(\cdot,0) = \bbF_0 \qquad \text{in } \Omega.
\end{align}
Here, $\rho, \nu, \mu, \lambda>0$ are fixed physical constants related to mass density, solvent viscosity, elastic modulus, and elastic relaxation, respectively.

\paragraph{Conformation tensor formulation.}
In many applications and existing analyses, the focus lies on the conformation tensor (or left Cauchy--Green tensor) $\mathbb{B} = \mathbb{F} \mathbb{F}^\top$, which is symmetric. 
The evolution equation for $\mathbb{B}$ corresponding to \eqref{eq:F} is:
\begin{align}
    \label{eq:B}
    \partial_t \bbB + (\bv\cdot\nabla) \bbB + \frac{\mu}{\lambda} (\bbB^2 - \bbB) 
    &= (\nabla\bv) \bbB + \bbB (\nabla\bv)^\top,
\end{align}
where $\mathbb{B}$ is required to be positive definite in $\Omega\times(0,T)$. The quadratic term $\mathbb{B}^2$ in \eqref{eq:B} represents the nonlinear damping characteristic of the Giesekus model, distinguishing it from the linear relaxation found in the Oldroyd-B model \cite{Hu_Lelievre_2007}, which corresponds to \eqref{eq:B} with $\frac{\mu}{\lambda} (\bbB^2 - \bbB)$ replaced by $\frac{\mu}{\lambda} (\bbB - \bbI)$, where $\bbI \in \bbR^{d\times d}$ denotes the identity matrix.

\paragraph{Energy estimates.}
The system \eqref{eq:system} admits a fundamental energy dissipation law. Formally, testing \eqref{eq:v} with $\mathbf{v}$ and \eqref{eq:F} with $\mu(\mathbb{F} - \mathbb{F}^{-\top})$, and using the incompressibility condition \eqref{eq:div}, yields the identity
\begin{align}
    \label{eq:stability} \nonumber
    &\int_{\Omega} \left( \frac\rho2 \abs{\bv(t)}^2 + \frac\mu2 \abs{\bbF(t)}^2 - \frac\mu2 \ln\det(\bbF(t)\bbF(t)^\top) \right) 
    \\ \nonumber
    &\quad 
    + \int_0^t \int_{\Omega} \left( \nu\abs{\nabla\bv}^2 
    + \frac{\mu^2}{2\lambda} \abs{\bbF \bbF^\top - \bbI }^2 \right) 
    \\
    &= \int_{\Omega} \left( \frac\rho2 \abs{\bv_0}^2 + \frac\mu2 \abs{\bbF_0}^2 - \frac\mu2 \ln\det(\bbF_0\bbF_0^\top) \right) , 
\end{align}
for all $t\in(0,T)$.
This identity ensures that the energy, consisting of the kinetic energy and the elastic energy, is dissipated by viscous effects and elastic relaxation. Note that the logarithmic part $-\ln\det(\mathbb{F}\mathbb{F}^\top)$ of the elastic energy penalizes configurations with small determinants and thus enforces the constraint $\det(\mathbb{F}) > 0$, assuming that $\det(\mathbb{F}_0) > 0$ .

\paragraph{Non-dimensionalization.}
For numerical simulations, we consider the dimensionless form of \eqref{eq:system}. Letting $v_c$ and $x_c$ be characteristic velocity and length scales, we define the Reynolds number ($\textit{Re}$), Weissenberg number ($\textit{Wi}$), and viscosity ratio ($\alpha$) as follows:
\begin{align}\label{eq:nondim_parameters}
    \textit{Re} = \frac{\rho v_c x_c}{\nu + \lambda},
    \quad
    \textit{Wi} = \frac{\lambda v_c}{\mu x_c},
    \quad
    \alpha = \frac{\lambda}{\nu + \lambda} \in (0,1).
\end{align}
Then, in nondimensional form, \eqref{eq:v} and \eqref{eq:F} read
\begin{subequations}\label{eq:system_nondim}
\begin{align}
    \label{eq:v_nondim}
    \textit{Re} \left( \partial_t \bv + (\bv\cdot\nabla) \bv \right)
    + \nabla p - (1-\alpha) \Delta \bv &= \frac{\alpha}{\textit{Wi}} \Div( \bbF \bbF^\top),
    \\
    \label{eq:F_nondim}
    \partial_t \bbF + (\bv\cdot\nabla) \bbF + \frac{1}{2 \textit{Wi}} (\bbF \bbF^\top \bbF - \bbF) 
    &= (\nabla\bv) \bbF .
\end{align}
\end{subequations}

\paragraph{State of the art and main goals.}

The analysis of viscoelastic fluids often relies on thermodynamically consistent laws. A thermodynamic framework using entropy estimates was introduced in \cite{Hu_Lelievre_2007}, while a rigorous approach for deriving viscoelastic models such as Oldroyd-B and Giesekus was developed in \cite{malek_prusa_2018}. 

Among these models, the Oldroyd-B model has been studied extensively. 
However, existence proofs often fail without regularization due to the hyperbolic nature of the equations and a lack of suitable a priori estimates. A common approach is to introduce artificial stress diffusion, where a thermodynamically consistent derivation was presented in~\cite{malek_2018_Oldroyd_diffusive}. A rigorous existence theory and a finite element approximation for a diffusive Oldroyd-B model were developed in~\cite{barrett_boyaval_2009}, while global regularity was established in~\cite{constantin_kliegl_2012}. Similar techniques were applied to more complex systems, including compressible flows~\cite{barrett_lu_sueli_2017}, two-phase flows~\cite{sieber_2020, LT_2022, GNT_2024}, and tumour growth models~\cite{GKT_2022_viscoelastic}. Related strategies exist for the FENE-P model \cite{barrett_2018_fene-p} and the Peterlin model \cite{lukacova_2015, Brunk_Lukacova_2021}.
For modified versions of the Oldroyd-B model, results on global existence were obtained in \cite{lions_masmoudi_2000}, while well-posedness was addressed in \cite{chemin_masmoudi_2001}. Strong solutions were investigated in \cite{chupin_2018}. In two dimensions, global existence of smooth solutions was shown in \cite{lin_zhang_zhang_2015}, while the limiting case of infinite Weissenberg number was analyzed in \cite{LinLiuZhang_2005}.

The Giesekus model, which is the focus of this work, has been addressed more recently. Weak solutions with stress diffusion are well-understood in both two and three dimensions (see \cite[Appendix B]{bulicek_2022_giesekus_2d}). Large-data existence in 3D using a modified energy combined with stress diffusion was established in \cite{bathory_2021_viscoelastic_3D}. Corresponding numerical approximations and their subsequence convergence were studied in \cite{T_2024_Enumath} and extended to two-phase models in \cite{GT_2023_DCDS}. Convergence to equilibria for thermo-viscoelastic fluids has also been studied in \cite{prusa_2021_equilibria_viscoelastic}. Furthermore, the existence of energy-variational solutions for a modified system in 3D was shown in \cite{agosti_lasarzik_rocca_2024}.
However, rigorous results for the Giesekus model without stress diffusion have only appeared in recent years. A preliminary step was taken by Masmoudi~\cite{masmoudi_2011}, who proved the weak sequential stability of solutions using the entropy framework from \cite{Hu_Lelievre_2007}. Although that work did not construct an approximating sequence needed for the analysis, the gap was recently closed in \cite{bulicek_2022_giesekus_2d}, which proved global existence of weak solutions in 2D. This result was later extended to 3D in \cite{bulicek_2025_giesekus_3d}. 

These analytical results are essential for understanding the numerical difficulties associated with viscoelastic flows. Simulating viscoelastic flows is challenging, especially at high Weissenberg numbers where the relaxation parameter is large. Standard methods often face stability issues in this regime, a phenomenon known as the high Weissenberg number problem \cite{alves_2021_viscoelastic_review, fattal_kupferman_2005_log_formulation}.
A promising strategy to address these difficulties relies on two key principles \cite{barrett_boyaval_2009}:
\begin{enumerate}
\item designing a discrete scheme with an energy dissipation property, and
\item rigorously proving a (subsequence) convergence result for the numerical approximation.
\end{enumerate}
The present work aims to complement the recent analytical advances for the Giesekus model in \cite{bulicek_2022_giesekus_2d} in 2D with a numerical method satisfying these two principles.

Our main result can be informally formulated as follows. For precise formulations, we refer to Theorems~\ref{theorem:existence_FE}, \ref{theorem:existence_time_discrete}, and~\ref{theorem:convergence}.
\begin{theorem*}
Let $\Omega$ be a polygonal Lipschitz domain for dimensions $d\in\{2,3\}$, and let $\bv_0$, $\bbF_0$ be suitable initial data. Then, under a mild restriction on the time step size, there exists a stable numerical approximation of \eqref{eq:system}. Moreover, in two dimensions, by first letting the spatial discretization parameter $h\to 0$ and then the time step size $\Delta t\to 0$, we obtain a subsequence that converges to a weak solution of \eqref{eq:system} subject to \eqref{eq:u_bc}--\eqref{eq:system_init}, satisfying the positivity constraint~\eqref{eq:detF}.
\end{theorem*}

\paragraph{Outline.}
The structure of this work is as follows. We end Section~\ref{sec:introduction} by introducing our notation.
Section~\ref{sec:weak_solution} recalls the definition of a weak solution to \eqref{eq:system} from \cite{bulicek_2022_giesekus_2d} and outlines our main strategy for constructing numerical solutions that converge (up to a subsequence) to such a weak solution in 2D.
Section~\ref{sec:numerical_scheme} introduces a fully discrete numerical scheme for \eqref{eq:system} and analyses its mathematical properties.
Section~\ref{sec:proof_2d} addresses the passage to the limit in the discretization parameters. For technical reasons (explained later), we first let the spatial discretization parameter $h \to 0$ while keeping the time step fixed, and then send the time step size $\Delta t \to 0$.
The main convergence result for $\Delta t \to 0$ is stated in Theorem~\ref{theorem:convergence}, where we identify subsequences converging to a weak solution of \eqref{eq:system} in two dimensions. The proof, based on techniques from \cite{bulicek_2022_giesekus_2d}, is developed in Sections~\ref{sec:limit_vanishing_dt}--\ref{sec:positivity_F}.
Finally, Section~\ref{sec:numerical_results} presents numerical experiments, including convergence studies and benchmark problems, to demonstrate the practicality of the proposed approach.

\paragraph{Notation.}
We now fix the notation used throughout this work.
Vector-valued and matrix-valued quantities are written in bold and blackboard bold font, respectively.
Let $d\in\{2,3\}$. For two vectors $\bv,\bw\in\bbR^d$, the Euclidean inner product is defined as
$\bv\cdot\bw = \bv^\top \bw = \bw^\top \bv$.
For two matrices $\bbA,\bbB\in\bbR^{d\times d}$, the Frobenius inner product is
$\bbA:\bbB = \trace(\bbA^\top\bbB) = \trace(\bbB^\top\bbA)$,
where $\trace(\bbA)$ denotes the trace of $\bbA$.
The tensor product of two vectors $\bv,\bw\in\bbR^d$ is given by
$(\bv\otimes\bw)_{i,j} = \bv_i \bw_j$ for $i,j\in\{1,\dots,d\}$.
For a vector $\bv\in\bbR^d$ and a matrix $\bbB\in\bbR^{d\times d}$, we set
$(\bv\otimes\bbB)_{i,j,k} = \bv_i \bbB_{j,k}$ for $i,j,k\in\{1,\dots,d\}$.
The Euclidean norm for vectors and the Frobenius norm for matrices are both denoted by $\abs{\cdot}$.
The zero vector is written as $\mathbf{0}\in\bbR^d$, the zero matrix as $\mathbb{O}\in\bbR^{d\times d}$, and the identity matrix as $\mathbb{I}\in\bbR^{d\times d}$.
For a vector $\bv\in\bbR^d$ and a matrix $\bbB\in\bbR^{d\times d}$, their gradients are defined by
$(\nabla\bv)_{i,j} = \partial_{x_j}\bv_i$ and $(\nabla\bbB)_{i,j,k} = \partial_{x_k}\bbB_{i,j}$.
The divergence of a vector $\bv$ is $\Div(\bv) = \sum_{j=1}^d \partial_{x_j}\bv_j$,
and the divergence of a matrix $\bbB$ is $(\Div\bbB)_i = \sum_{j=1}^d \partial_{x_j}\bbB_{i,j}$.
For a third-order tensor $\bbB\otimes\bv$ we set
$\Div(\bbB\otimes\bv) = \sum_{j=1}^d \partial_{x_j}(\bv_j\bbB)$.
The convective derivative of $\bbB\in\bbR^{d\times d}$ with respect to $\bv\in\bbR^d$ is defined by
$(\bv\cdot\nabla)\bbB = \sum_{k=1}^d \bv_k \partial_{x_k}\bbB$.
For $\bbA,\bbB\in\bbR^{d\times d}$ we also write
$\nabla\bbA:\nabla\bbB = \sum_{k=1}^d (\partial_{x_k}\bbA):(\partial_{x_k}\bbB)$.

If $X$ is a real Banach space, we denote its norm by $\norm{\cdot}_X$, its dual by $X'$, and the duality pairing by $\dualp{\cdot}{\cdot}_{X',X}$.
If $X=H$ is a Hilbert space, we denote the inner product by $\skp{\cdot}{\cdot}_H$.
For $p\in[1,\infty]$, $m\ge 0$ and an open set $U\subset\bbR^n$, $n\in\bbN$, we use the standard spaces $L^p(U;X)$ and $W^{m,p}(U;X)$.
If $p=2$ and $X$ is a Hilbert space, we also write $H^m(U;X)=W^{m,2}(U;X)$.
If $U=(0,T)$ with $T>0$, then $W^{m,p}(0,T;X)$ denotes $W^{m,p}((0,T);X)$.
If $X=\bbR$ and $U=\Omega$, we write $W^{m,p}(\Omega)$ instead of $W^{m,p}(\Omega;\bbR)$.
When $X\in\{\bbR,\bbR^d,\bbR^{d\times d}\}$ is clear from context, we simply write $W^{m,p}=W^{m,p}(\Omega;X)$.
The seminorm on $W^{m,p}$ is denoted by $\abs{\cdot}_{W^{m,p}}$.
The spaces $C([0,T];X)$ and $C_w([0,T];X)$ denote strongly and weakly continuous functions on $[0,T]$ with values in $X$.
We write $C^\infty(\overline\Omega;X)$ for smooth functions on $\overline\Omega$ and $C_c^\infty(\Omega;X)$ for smooth, compactly supported functions.
The space $H^1_0(\Omega;\bbR^d)$ is the closure of $C_c^\infty(\Omega;\bbR^d)$ in the $H^1$ norm.
We set $L^2_{\Div}(\Omega;\bbR^d)$ as the closure of $\{\bv\in C_c^\infty(\Omega;\bbR^d)\mid \Div(\bv)=0\}$ in $L^2$, and define $H^1_{0,\Div}(\Omega;\bbR^d) = H^1_0(\Omega;\bbR^d)\cap L^2_{\Div}(\Omega;\bbR^d)$.
Finally, $\mathcal{M}(\overline{\Omega_T})$ denotes the space of Radon measures on the closure of $\Omega_T=\Omega\times(0,T)$.

\section{Weak solution} \label{sec:weak_solution}
In this work, we let $T>0$ and $\Omega\subset\bbR^d$, $d\in\{2,3\}$, be a bounded Lipschitz domain with a polygonal (or polyhedral, respectively) boundary $\partial\Omega$. The restriction to polygonal boundaries is required for technical reasons in order to avoid boundary approximation techniques at the finite element level.

%
%
We now recall the definition of a weak solution to \eqref{eq:system} in two dimensions from \cite{bulicek_2022_giesekus_2d}. 

\begin{definition} 
\label{def:weak_solution}
Let $d=2$. Let the initial data $\bv_0 \in L^2_{\Div}(\Omega;\bbR^2)$, $\bbF_0 \in L^2(\Omega;\bbR^{2\times 2})$, with $\det(\bbF_0)>0$ a.e.~in $\Omega$ and $\ln\det(\bbF_0) \in L^1(\Omega)$, be given. We call $(\bv, \bbF)\colon\, \Omega\times(0,T) \to \bbR^2 \times \bbR^{2\times 2}$ a weak solution to \eqref{eq:system} if
\begin{align*}
    &\bv \in C([0,T]; L^2_{\Div}(\Omega;\bbR^2)) \cap L^2(0,T; H^1_{0,\Div}(\Omega;\bbR^2))\, ,
    \\
    &\partial_t \bv \in L^2(0,T; (H^1_{0,\Div}(\Omega;\bbR^2))')\, ,
    \\
    &\bbF \in C([0,T]; L^2(\Omega;\bbR^{2\times 2})) \cap L^4(\Omega_T;\bbR^{2\times 2})\, ,
    \\
    &\partial_t \bbF \in L^{4/3}(0,T; (H^1(\Omega;\bbR^{2\times 2}))')\, ,
    \\
    &\det(\bbF)>0 \text{ a.e.~in } \Omega_T\, ,
\end{align*}
such that, for all $\bw\in H^1_{0,\Div}(\Omega;\bbR^2)$ and $\bbG\in H^1(\Omega;\bbR^{2\times 2})$ and almost all $t\in(0,T)$, 
\begin{subequations}
\label{eq:system_weak}
\begin{align}
    \label{eq:v_weak}
    \rho \dualp{\partial_t \bv}{\bw}_{(H^1_{0,\Div})',H^1_{0,\Div}}
    - \rho \skp{\bv\otimes\bv}{\nabla\bw}_{L^2}
    + \nu \skp{\nabla\bv}{\nabla\bw}_{L^2}
    + \mu \skp{\bbF\bbF^{\top}}{\nabla\bw}_{L^2}
    &= 0 \, ,
    \\
    \label{eq:F_weak}
    \dualp{\partial_t \bbF}{\bbG}_{(H^1)',H^1}
    - \skp{\bbF\otimes\bv}{\nabla\bbG}_{L^2}
    - \skp{(\nabla\bv)\bbF}{\bbG}_{L^2}
    + \frac{\mu}{2\lambda} \skp{\bbF\bbF^\top \bbF - \bbF}{\bbG}_{L^2}
    &= 0 \, ,
\end{align}
\end{subequations}
and where the initial data are attained in the sense
\begin{align}
    \lim_{t\to0^+} \left( \norm{\bv(t) - \bv_0}_{L^2} 
    + \norm{\bbF(t) - \bbF_0}_{L^2} \right) 
    = 0 \, .
\end{align}
\end{definition}

As usual, the weak formulation \eqref{eq:system_weak} is derived from \eqref{eq:system} together with the boundary and initial conditions \eqref{eq:u_bc} and \eqref{eq:system_init} by taking inner products with suitable test functions and applying integration by parts over $\Omega$. In addition, Definition~\ref{def:weak_solution} already incorporates the correct regularity of the solution as well as the required positivity of $\det(\bbF)$.


Recently, the following existence theorem has been proved in \cite{bulicek_2022_giesekus_2d}.

\begin{theorem}
\label{theorem:weak_solution}
There exists a weak solution to \eqref{eq:system} in the sense of Definition~\ref{def:weak_solution}.
\end{theorem}

The main results of this work are formulated in Theorems~\ref{theorem:existence_FE}, \ref{theorem:existence_time_discrete}, and~\ref{theorem:convergence}. We present a stable numerical approximation which, by sending first the spatial discretization parameter $h\to 0$ and then the time step size $\Delta t\to 0$, converges (up to a subsequence) to a weak solution in two dimensions in the sense of Definition~\ref{def:weak_solution}. Thereby, we recover the result of Theorem~\ref{theorem:weak_solution}.

Our analysis relies on the fully discrete scheme presented in Section~\ref{sec:numerical_scheme}, where the main properties are summarized in Theorem~\ref{theorem:existence_FE}. The limit passage in the discretization parameters is addressed in Section~\ref{sec:proof_2d}. 

For technical reasons, we first consider the limit $h \to 0$. This step is necessary to adapt the arguments from \cite{bulicek_2022_giesekus_2d}. Specifically, identifying limits of products involving $\bbF$ requires the compactness of $\bbF$ (see Section~\ref{sec:compactness_F}). Unfortunately, this cannot be obtained directly from energy estimates without additional regularization of the form $\varepsilon\Delta \bbF$ with fixed $\varepsilon>0$. 
Hence, proving compactness requires testing the approximate equation for $\bbF$ with products of functions. At the space-discrete level, such functions are not admissible test functions, and projecting them to finite element spaces would make the subsequent computations infeasible. Therefore, passage to the limit $h \to 0$ is performed first in Section~\ref{sec:time_discrete} and summarized in Theorem~\ref{theorem:existence_time_discrete}.

The limit passage $\Delta t\to 0$ is established in Theorem~\ref{theorem:convergence}, with proofs provided in Sections~\ref{sec:limit_vanishing_dt}--\ref{sec:positivity_F}.
We first apply weak-($*$) compactness results based on a priori estimates to obtain limit functions (Section~\ref{sec:limit_vanishing_dt}). However, identifying certain nonlinear terms is challenging due to the lack of compactness for $\bbF$. 
To resolve this, we adapt the strategy of \cite{bulicek_2022_giesekus_2d} to the time-discrete setting, which consists of two steps:
First, we show the strong continuity of $\bbF$ in Section~\ref{sec:continuity_F}.
Then, we prove the compactness of the deformation gradient in two dimensions in Section~\ref{sec:compactness_F}. This involves deriving a differential inequality for the difference between a nonlinear quantity involving $\bbF$ and its unidentified limit. This step requires testing the time-discrete velocity equation with non-divergence-free functions, necessitating the reconstruction of local pressures. A renormalization argument finally allows us to identify the nonlinear term and deduce the strong convergence of the deformation gradient.
Lastly, we deduce the positivity of $\det(\bbF)$ in Section~\ref{sec:positivity_F} under suitable assumptions on the initial data.

\section{Fully discrete numerical schemes} \label{sec:numerical_scheme}
In this section, we present a numerical approximation of \eqref{eq:system} that is based on the weak formulation from Definition~\ref{def:weak_solution}. We will prove the existence and stability of discrete solutions in two and three space dimensions in Theorem~\ref{theorem:existence_FE}. By stability, we mean uniform boundedness of the sequence of discrete solutions generated by the numerical method, with respect to the spatial and temporal discretization parameters.
In Section~\ref{sec:numerical_scheme_linear}, we also discuss an alternative numerical scheme that is stable, fully linear, and thus may be useful for experiments, but for which we were not able to rigorously prove (subsequence) convergence.

\subsection{Discrete setting}
We divide the time interval $[0,T)$ into equidistant subintervals $[t^{n-1}, t^n)$ with $t^n = n \Delta t$ and $t^{N_T} = T$, $N\in\bbN$ and $n\in\{1,\ldots,N_T\}$.
We require $\{\calT^h\}_{h>0}$ to be a family of conforming partitions of $\Omega$ into disjoint open simplices such that $\overline\Omega = \bigcup_{K\in\calT^h} \overline K$, and we set $h \coloneqq \max_{K\in\calT^h} \diam(K)$.
We assume throughout that the family of triangulations is shape-regular 
in the sense of \cite{bartels_2016}, i.e., there exists a constant $C>0$ independent of $h$, such that
\begin{align*}
    \sup_{K\in\calT^h} \frac{h_K}{ r_K } \leq C,
\end{align*}
where $h_K = \diam(\overline{K})$, and $r_K$ denotes the diameter of the largest inscribed ball in $\overline{K}$.

For $k\in\bbN$, we introduce the following notation for the finite element spaces of continuous and piecewise polynomial functions
\begin{align*}
     \calS_k^h &\coloneqq \left\{ q \in C(\overline\Omega) \mid q_{|K} \in \calP_k(K) \quad \forall \,  K\in\calT^h \right\} 
     \ \subset \ W^{1,\infty}(\Omega) \, ,
\end{align*}
and the finite element space of discontinuous and piecewise constant functions
\begin{align*}
    \calS_0^h &\coloneqq \left\{ q \in L^1(\Omega) \mid q_{|K} \in \calP_0(K) \quad \forall \,  K\in\calT^h \right\} 
    \ \subset \ L^\infty(\Omega) \, ,
\end{align*}
where $\calP_k(K)$ denotes the set of polynomials of degree less than or equal to $k\geq0$, defined on $K\in\calT^h$.
Similarly, we introduce the notation 
\begin{align*}
    \mathbf{S}_k^h \coloneqq [\calS_k^h]^d,
    \qquad
    \mathbb{S}_k^h \coloneqq [\calS_k^h]^{d\times d},
\end{align*}
for vector and matrix valued finite element functions of degree less than or equal to $k\geq 0$, respectively. Moreover, for $k\geq 1$, we define
\begin{align*}
    \mathbf{U}_{k+1}^h &\coloneqq \mathbf{S}_{k+1}^h \cap H^1_0(\Omega;\bbR^d),
    \\
    \mathbf{U}_{k+1,\Div}^h &\coloneqq \{\bv_h \in \mathbf{U}_{k+1}^h \mid \skp{\Div \bv_h}{q_h}_{L^2} = 0 \quad \forall\, q_h\in \calS_{k}^h\} .
\end{align*}

Since the family of meshes $\{\calT^h \}_{h>0}$ is shape-regular, the pair of finite element spaces $(\mathbf{U}_{k+1}^h, \calS_k^h)$ with $k\geq 1$, also known as the Taylor--Hood element, satisfies the following inequality:
\begin{align}
    \label{eq:LBB}
    \sup_{\bw_h\in \, \mathbf{U}_{k+1}^h }
    \frac{\skp{\Div\bw_h}{q_h}_{L^2}}{ \norm{\bw_h}_{H^1} }
    \geq C_0 \norm{q_h}_{L^{2} } \quad \forall\, q_h\in \calS_k^h,
\end{align}
where $C_0>0$ denotes a constant that is independent of the mesh parameter $h>0$.
The inequality \eqref{eq:LBB} is often referred to as the inf-sup, or Babu{\v s}ka--Brezzi, condition. 
Sometimes, for technical reasons, an additional (but very mild) condition is imposed on the mesh, although the additional condition can be dropped if the triangulation $\calT^h$ is fine enough; see \cite[Chap.~8.8]{BBF2013_fem}.

For future reference, we recall the standard Lagrange interpolation operators \cite{brenner_scott_2008}, denoted by $\mathrm{I}^h_k \colon\, C(\overline\Omega) \to \calS_k^h$, $k\geq 1$, and we naturally extend their definitions to matrix-valued functions, i.e., $\bbI^h_1 \colon\, C(\overline\Omega;\bbR^{d\times d})\to \bbS_1^h$. 

The following local inverse estimate holds true for any $u_h \in \calS_k^h$, $k\geq 1$, $K\in\calT^h$, $s,m\in\{0,1\}$ with $s\leq m$ and $1\leq r \leq p \leq \infty$, see \cite{brenner_scott_2008},
\begin{align}
	\label{eq:inverse_estimate}
	\abs{u_h}_{W^{m,p}(K)} 
	&\leq 
	C h_K^{s - m + \frac{d}{p} - \frac{d}{r}} \abs{u_h}_{W^{s,r}(K)}.
\end{align}
%
%
%
%

\subsection{The numerical scheme}
We now introduce our numerical approximation of \eqref{eq:system}.
Afterwards, in Theorem~\ref{theorem:existence_FE}, we will show the existence of discrete solutions and derive a discrete energy inequality.

Let $k,m \in \bbN$. Given the discrete initial data 
\begin{align*}
    \bv_h^0 \in \mathbf{U}_{k+1,\Div}^h
    \quad \text{and} \quad 
    \bbF_h^0 \in \bbS_m^h,
\end{align*}
for $n\in\{1,\ldots,N_T\}$ we aim to find a solution 
\[(\bv_h^n, p_h^n, \bbF_h^n) \in \mathbf{U}_{k+1}^h \times \calS_k^h \times \bbS_m^h \]
such that, for all test functions $(\bw_h, q_h, \bbG_h) \in \mathbf{U}_{k+1}^h \times \calS_k^h \times \bbS_m^h$,
\begin{subequations}
\label{eq:system_FE}
\begin{align}
    \nonumber
    \label{eq:v_FE}
    0 &= 
    \frac\rho{\Delta t}\skp{\bv_h^n - \bv_h^{n-1}}{\bw_h}_{L^2}
    + \frac\rho{2} \skp{(\bv_h^{n-1} \cdot \nabla) \bv_h^n}{\bw_h}_{L^2}
    - \frac\rho{2} \skp{\bv_h^n}{(\bv_h^{n-1} \cdot \nabla) \bw_h}_{L^2}
    \\*
    &\quad
    + \nu \skp{\nabla \bv_h^n}{\nabla\bw_h}_{L^2}
    - \skp{p_h^n}{\Div\bw_h}_{L^2}
    + \mu \skp{\bbF_h^{n} (\bbF_h^{n})^\top}{\nabla\bw_h}_{L^2} ,
    \\
    \label{eq:div_FE}
    0 &= \skp{\Div\bv_h^n}{q_h}_{L^2} ,
    \\
    \label{eq:F_FE}
    \nonumber
    0 &= 
    \frac{1}{\Delta t} \skp{\bbF_h^n - \bbF_h^{n-1}}{\bbG_h}_{L^2}
    + \bc_h(\bv_h^{n-1}, \bbF_h^n, \bbG_h)
    + \frac{\mu}{2 \lambda} \skp{ \bbF_h^{n} (\bbF_h^{n})^\top \bbF_h^n- \bbF_h^{n}} {\bbG_h}_{L^2} 
    \\* 
    &\quad 
    - \skp{ (\nabla\bv_h^n) \bbF_h^n }{\bbG_h}_{L^2}
    + \Delta t \skp{\nabla\bbF_h^n}{\nabla\bbG_h}_{L^2}.
\end{align}
\end{subequations}

For given $\bv_h\in\mathbf{U}_{k+1,\Div}^h$ and $\bbF_h,\bbG_h\in\bbS^h_m$, the term $\bc_h(\bv_h,\bbF_h,\bbG_h)$ defines a numerical approximation of the convection term $\skp{\bbF\otimes\bv}{\nabla\bbG}_{L^2}$ in \eqref{eq:F_weak}, where $\bc_h(\cdot,\cdot,\cdot)$ is a trilinear form with the property
\begin{align}
    \label{eq:convective_chain_rule_def}
    \bc_h(\bv_h,\bbF_h,\bbF_h) = 0.
\end{align}
Moreover, let $\bv\in H^1_{0,\Div}(\Omega;\bbR^d)$, $\bbF\in H^1(\Omega;\bbR^{d\times d})$ and $\bbG\in C^\infty(\overline{\Omega};\bbR^{d\times d})$. Assume that $\bv_h \in \mathbf{U}_{k+1,\Div}^h$ and $\bbF_h,\bbG_h \in \bbS_m^h$ are given such that, as $h\to 0$, 
\begin{alignat*}{3}
    \bv_h &\to \bv \quad &&\text{strongly} \quad &&\text{in } L^2_{\Div}(\Omega;\bbR^d),
    \\
    \bbF_h &\rightharpoonup \bbF \quad &&\text{weakly} \quad &&\text{in } H^1(\Omega;\bbR^{d\times d}),
    \\
    \bbG_h &\to \bbG \quad &&\text{strongly} \quad &&\text{in } W^{1,\infty}(\Omega;\bbR^{d\times d}).
\end{alignat*}
Then we also require that
\begin{align}
    \label{eq:convective_convergence}
    \bc_h(\bv_h,\bbF_h,\bbG_h) \to - \skp{\bbF\otimes\bv}{\nabla\bbG}_{L^2}, \quad \text{as } h\to 0.
\end{align}
In Section~\ref{sec:convective}, we present three specific examples of $\bc_h(\cdot,\cdot,\cdot)$ that satisfy \eqref{eq:convective_convergence}.

\begin{remark}
The stabilization term $\Delta t \skp{\nabla\bbF_h^n}{\nabla\bbG_h}_{L^2}$ in \eqref{eq:F_FE} corresponds to a discrete Laplacian subject to a homogeneous Neumann boundary condition, commonly referred to as stress diffusion. This artificial diffusion is needed for the convergence analysis as the mesh size $h$ tends to zero. Specifically, it provides the necessary uniform control of the discrete deformation gradient $\bbF_h^n$ in the $H^1$ norm to handle the limit passage in the convective term $\bc_h(\bv_h^{n-1},\bbF_h^n,\bbG_h)$ (see Section~\ref{sec:convective}). While the stress diffusion term in \eqref{eq:F_FE} scales with $\Delta t$, our method also allows for more general scaling. One can use a coefficient $\phi(\Delta t)$, provided that $\phi \in C([0,T];\mathbb{R}_{\geq 0})$ is a monotonically increasing function, with $\phi(0)=0$. Once the spatial limit $h \to 0$ is established, we adapt the arguments from \cite{bulicek_2022_giesekus_2d} for time-discrete, spatially continuous solutions, and then pass to the limit $\Delta t\to 0$. The limit passages $h\to 0$ and $\Delta t\to 0$ are addressed in Section~\ref{sec:proof_2d}.
\end{remark}

\subsection{Approximation of the convective derivative}
\label{sec:convective}
In this work, we will often use the identities
\begin{align}
    \label{eq:chain_rule_continuous} 
    \skp{(\bv\cdot\nabla) \bbF}{\bbF}_{L^2}
    &= \frac12 \skp{\bv}{\nabla \abs{\bbF}^2}_{L^2}
    = -\frac12 \skp{\Div \bv}{\abs{\bbF}^2}_{L^2} 
    = 0,
\end{align}
for all $\bv\in H^1_{0,\Div}(\Omega;\bbR^d)$ and $\bbF\in H^1(\Omega;\bbR^{d\times d})$. However, replacing $\bv$ and $\bbF$ with finite element functions does not guarantee equalities analogous to those in \eqref{eq:chain_rule_continuous}.
Depending on the degrees $k\geq 1$ and $m\geq 1$ of the finite element spaces $(\mathbf{U}_{k+1}^h, \mathcal{S}_{k}^h)$ for the velocity and pressure and $\mathbb{S}_m^h$ for the deformation gradient, we propose different numerical approximations for the convective derivative featuring in \eqref{eq:F_weak}. We note that by integration by parts over $\Omega$, one has
\begin{align}
    \label{eq:convective_general}
    \nonumber
    \skp{(\bv\cdot\nabla) \bbF}{\bbG}_{L^2} 
    &= - \skp{\bbF}{(\bv\cdot\nabla) \bbG}_{L^2} 
    = - \skp{\bbF \otimes \bv }{\nabla \bbG}_{L^2} 
    \\
    &= \frac12 \skp{(\bv\cdot\nabla) \bbF}{\bbG}_{L^2} - \frac12 \skp{\bbF}{(\bv\cdot\nabla) \bbG}_{L^2},
\end{align}
for any $\bv\in H^1_{0,\Div}(\Omega;\bbR^d)$ and $\bbF,\bbG \in H^1(\Omega;\bbR^{d\times d})$.

\subsubsection{Skew-symmetric approximation}
Let $k\geq 1$ and $m\geq 1$. Then, the simplest approximation of the convective term \eqref{eq:convective_general} is given by the skew-symmetric trilinear form defined by
\begin{align}
    \label{eq:convective_skewsymm}
    \bc_h(\bv_h,\bbF_h,\bbG_h) = 
    \frac12 \skp{(\bv_h\cdot\nabla)\bbF_h}{\bbG_h}_{L^2}
    - \frac12 \skp{\bbF_h}{(\bv_h\cdot\nabla)\bbG_h}_{L^2},
\end{align}
for all $\bv_h \in \mathbf{U}^h_{k+1}$ and $\bbF_h,\bbG_h \in \bbS^h_m$.
By construction, we have
\begin{align*}
    \bc_h(\bv_h,\bbF_h,\bbF_h) = 0,
\end{align*}
for all $\bv_h\in\mathbf{U}^h_{k+1}$ and $\bbF_h\in\bbS_m^h$.

\subsubsection{Approximation with higher-order elements}
We now present a second approach for the approximation of $\skp{(\bv\cdot\nabla)\bbF}{\bbG}_{L^2}$, where we now use higher-order Taylor--Hood elements $(\mathbf{U}^h_{k+1}, \mathcal{S}^h_{k})$ with $k\geq 2$, while the deformation gradient is approximated using lower-order elements $\bbS^h_m$ with $2m\leq k$.
The associated numerical approximation is then defined as follows:
\begin{align}
    \label{eq:convective_higherorder}
    \bc_h(\bv_h,\bbF_h,\bbG_h) 
    = - \skp{\bbF_h}{(\bv_h\cdot\nabla)\bbG_h}_{L^2},
\end{align}
for all $\bv_h \in \mathbf{U}^h_{k+1}$ and $\bbF_h,\bbG_h \in \bbS^h_m$.

For any $\bbF_h \in \bbS_m^h$, the quantity $\frac12|\bbF_h|^2 \in \calS_{2m}^h$ is contained in the discrete pressure space $\calS_{k}^h$ whenever $2m\leq k$ is fulfilled. 
Hence, for $2m\leq k$ and for discretely divergence-free velocities $\bv_h \in \mathbf{U}^h_{k+1,\Div}$, we have with the chain rule and integration by parts over $\Omega$ that
\begin{align*}
    \bc_h(\bv_h,\bbF_h,\bbF_h) 
    = - \frac12 \skp{\bv_h}{\nabla \abs{\bbF_h}^2}_{L^2}
    = \frac12 \skp{\Div \bv_h}{|\bbF_h|^2}_{L^2}
    = 0.
\end{align*}

\subsubsection{Approximation based on discrete analogues of chain rules}
In contrast to \eqref{eq:convective_higherorder}, we now describe an approach that works for all Taylor--Hood elements $(\mathbf{U}_{k+1}^h,\mathcal{S}_k^h)$ with $k\geq 1$, including the lowest-order case $(\mathbf{U}_2^h,\mathcal{S}_1^h)$.
In this setting, the deformation gradient is approximated using linear finite elements $\bbS_m^h$ with $m=1$. Hence, the quantity $\frac12|\bbF_h|^2 \in \calS_{2}^h$, with $\bbF_h \in \bbS_1^h$, is generally not contained in the discrete pressure space $\calS_{k}^h$ for the lowest-order Taylor--Hood element with $k=1$.

Thus, we introduce the definition
\begin{align}
    \label{eq:convective_Lambda}
    \bc_h(\bv_h,\bbF_h,\bbG_h) 
    = - {\sum}_{i,j=1}^d \skp{(\bv_h)_i \, \mathbf{\Lambda}_{i,j}(\bbF_h)}{\partial_{x_j} \bbG_h}_{L^2},
\end{align}
for all $\bv_h \in \mathbf{U}_{k+1}^h$ and $\bbF_h,\bbG_h \in \bbS^h_1$. 
Here, $\mathbf{\Lambda}_{i,j}(\bbF_h) \in \bbS_0^h$, $i,j\in\{1,\ldots,d\}$, is a suitable approximation of $\delta_{i,j} \bbF_h\in \bbS_1^h$, where $\delta_{i,j}$ denotes the Kronecker delta.
For the precise definition and important properties, we refer to \eqref{eq:def_Lambda}, \eqref{eq:Lambda_bound}, and \eqref{eq:error_Lambda}.
Let us note here that for any $\bbF_h\in \bbS_1^h$, one has, on any $K\in\calT^h$ and for each $i\in\{1,\ldots,d\}$,
\begin{align}
    \label{eq:Lambda_chain_rule}
    \sum\limits_{j=1}^d 
    {\mathbf\Lambda}_{i,j} (\bbF_h) : \partial_{x_j} \bbF_h
    = \frac12 \partial_{x_i} \mathrm{I}_1^h\left[ \abs{\bbF_h}^2 \right] \, .
\end{align}
This allows us to compute
\begin{align*}
    \bc_h(\bv_h,\bbF_h,\bbF_h) 
    = - \frac12  \skp{\bv_h }{ \nabla \mathrm{I}_1^h\left[ \abs{\bbF_h}^2 \right]}_{L^2}
    = \frac12  \skp{\Div(\bv_h) }{ \mathrm{I}_1^h\left[ \abs{\bbF_h}^2 \right]}_{L^2}
    = 0,
\end{align*}
for any $\bv_h\in\mathbf{U}^h_{k+1,\Div}$ and $\bbF_h\in\bbS^h_1$.

The approximation \eqref{eq:convective_Lambda} of the convective term is based on a framework that provides discrete versions of chain rules for matrix-valued linear finite element functions. This framework was first developed in \cite{barrett_boyaval_2009} for the Oldroyd-B model and has later been applied to other viscoelastic models \cite{barrett_2018_fene-p, GT_2023_DCDS, T_2024_Enumath}.
Unlike \eqref{eq:convective_skewsymm} and \eqref{eq:convective_higherorder}, the approach \eqref{eq:convective_Lambda} can be extended to a wider class of nonlinear elastic energy densities. The only assumptions needed are that the elastic energy density is continuously differentiable and strictly convex on its domain of definition.

We now present the precise definition and important properties of $\mathbf{\Lambda}_{i,j}(\bbF_h) \in \bbS_0^h$, $i,j\in\{1,\ldots,d\}$. 
Let $\hat K$ denote the standard open reference simplex in $\bbR^d$ with vertices $\hat{\mathbf{p}}_0, \ldots, \hat{\mathbf{p}}_d$,
where $\hat{\mathbf{p}}_0$ is the origin and $\hat{\mathbf{p}}_i$, $i\in\{1,\ldots,d\}$, is the $i$-th standard basis vector in the $\bbR^d$. 
Given a simplex $K\in\calT^h$ with vertices $\mathbf{p}_0^K, \ldots, \mathbf{p}_d^K$, then there exists an invertible matrix $\calA_K \in \bbR^{d\times d}$ such that the affine transformation 
\begin{align}
	\label{eq:trafo_K}
	\calB_K \colon\, \hat K \to K, \quad
    \hat{\mathbf{x}}\mapsto \mathbf{p}_0^K + \calA_K \hat{\mathbf{x}}
\end{align}
maps vertex $\hat{\mathbf{p}}_i$ to vertex $\mathbf{p}_i^K$, $i\in\{0,\ldots,d\}$; see \cite[Sec.~3.1.3]{bartels_2016}.
Thanks to the assumed shape-regularity of $\{\calT^h\}_{h>0}$, we have
\begin{align}
	\label{eq:shape_regular}
	\abs{\calA_K^\top} \abs{(\calA_K^\top)^{-1}} \leq C \quad \forall\, K\in\calT^h.
\end{align}
Let $K\in\calT^h$ be an arbitrary non-degenerate simplex with vertices $\mathbf{p}_0^K, \ldots, \mathbf{p}_d^K$. 
Similar to \cite{barrett_boyaval_2009}, we introduce the notation
\begin{align*}
	\hat\bbF(\hat{\mathbf{x}}) \coloneqq \bbF(\calB_K(\hat{\mathbf{x}})),
	\quad\quad 
	(\hat{\mathbb{I}}_1^h \hat\bbF)(\hat{\mathbf{x}}) \coloneqq (\mathbb{I}_1^h \bbF)(\calB_K(\hat{\mathbf{x}})),
	\quad\quad \forall \,  \hat{\mathbf{x}}\in\hat K, \, \bbF\in C(\overline K; \bbR^{d\times d}),
\end{align*}
and we define $\bbF_j^K \coloneqq \bbF_h(\mathbf{p}_j^K)$, $j\in\{0,\ldots,d\}$, for any $\bbF_h\in\bbS_1^h$ and $K\in\calT^h$, where $\mathbf{p}_0^K, \ldots, \mathbf{p}_d^K$ denote the vertices of the simplex $K$. 
On the reference element $\hat K$ and for $i\in\{1,\ldots,d\}$, we define
\begin{align}
	\label{def:Lambda_hat}
    \hat{\mathbf\Lambda} _i(\hat\bbF_h) 
   & \coloneqq 
    \begin{cases}
    \bbF_i^K & \text{  if } \bbF_i^K = \bbF_0^K \, ,
    \\
    \bbF_i^K + \frac12 
    \big( \bbF_0^K - \bbF_i^K \big)
    & \text{ if } \bbF_i^K \not= \bbF_0^K \, .
    \end{cases}
\end{align}
It follows directly by the construction of $\hat{\mathbf\Lambda} _i(\hat\bbF_h)$, that
\begin{align}
    \hat{\mathbf\Lambda}_i(\hat\bbF_h)  : 
    \partial_{\hat x_i} \hat \bbF_h
    = \frac12 \partial_{\hat x_i} \hat{\mathrm{I}}_1^h \big[ 
    \abs{\hat \bbF_h} \big],
\end{align}
for all $i\in\{1,\ldots,d\}$. 
With the help of this definition on the reference element $\hat{K}$, we now define on an arbitrary element $K\in\calT^h$, for $i,j\in\{1,\ldots,d\}$, and for $\bbF_h\in\bbS_1^h$,
\begin{align}
    \label{eq:def_Lambda}
    {\mathbf\Lambda} _{i,j}(\bbF_h) _{|K}
	&\coloneqq
	\sum\limits_{m=1}^d [(\calA_K^\top)^{-1}]_{i,m} \,
	\hat{\mathbf\Lambda} _{m}(\hat\bbF_h) \, [\calA_K^\top]_{m,j} \, \in\bbR^{d\times d} \, .
\end{align}
We note that, for any $K\in\calT^h$ and for each $i,j\in\{1,\ldots,d\}$, ${\mathbf\Lambda} _{i,j}(\bbF_h)_{|K} \in \bbR^{d\times d}$ depends linearly on $\bbF_h \in \bbS_1^h$. 
Our definition is inspired by \cite{barrett_boyaval_2009, barrett_2018_fene-p, GT_2023_DCDS}. In those works, the study of more complex elastic energy densities requires replacing the constant factor $\frac12$ in \eqref{def:Lambda_hat} with nonlinear weight functions $\lambda_i(\hat{\bbF}_h) \in [0,1]$. This leads to a highly nonlinear dependency of ${\mathbf\Lambda}_{i,j}$ on $\bbF_h$. In contrast, our approach assumes $\lambda_i(\hat{\bbF}_h) = \frac12$, which ensures that ${\mathbf\Lambda}_{i,j}(\bbF_h)$ depends linearly on $\bbF_h \in \bbS_1^h$.

As the family of triangulations is shape-regular, we have
\begin{align}
    \label{eq:Lambda_bound}
    \max_{i,j\in\{1,\ldots,d\}} 
    \norm{ {\mathbf\Lambda} _{i,j}(\bbF_h) }_{L^\infty(K)} 
    \leq C \norm{ \bbF_h }_{L^\infty(K)} ,
\end{align}
for any $K\in\calT^h$ and $\bbF_h\in\bbS_1^h$.
Moreover, we have on any $K\in\calT^h$ and for any $i\in\{1,\ldots,d\}$, that
\begin{align*}
    \sum\limits_{j=1}^d 
    \mathbf\Lambda_{i,j} (\bbF_h) : \partial_{x_j} \bbF_h
    &=
    \sum\limits_{j,m=1}^d \Big( [(\calA_K^\top)^{-1}]_{i,m} \,
    \hat{\mathbf\Lambda}_{m}(\hat\bbF_h) \, [\calA_K^\top]_{m,j} \Big)
    : \Big( \sum\limits_{k=1}^d [(\calA_K^\top)^{-1}]_{j,k} \, \partial_{\hat x_k}  
    \hat \bbF_h \Big)
    \\
    &= \sum\limits_{k, m=1}^d 
    [(\calA_K^\top)^{-1}]_{i,m} \,
    \hat{\mathbf\Lambda}_{m}(\hat\bbF_h) \, 
    : \partial_{\hat x_k} \hat \bbF_h 
    \Big( 
    \sum\limits_{j=1}^d [\calA_K^\top]_{m,j} [(\calA_K^\top)^{-1}]_{j,k} \Big)
    \\
    &= \sum\limits_{m=1}^d 
    [(\calA_K^\top)^{-1}]_{i,m} \,
    \hat{\mathbf\Lambda}_{m}(\hat\bbF_h) \, 
    : \partial_{\hat x_m} \hat \bbF_h
    \\
    &= \frac12 \sum\limits_{m=1}^d 
    [(\calA_K^\top)^{-1}]_{i,m} \, 
    \partial_{\hat x_m} \hat{\mathrm{I}}_1^h \big[ 
    \abs{\hat \bbF_h}^2 \big]
    \\
    &= \frac12 \partial_{x_i} \mathrm{I}_1^h \left[ \abs{\bbF_h}^2 \right] \, ,
\end{align*}
where we have used that
\begin{align*}
    \partial_{x_i} q_h(\calB_K(\mathbf{\hat x})) = \sum_{m=1}^d [(\calA_K^\top)^{-1}]_{i,m} \, \partial_{\hat x_m} \hat q_h(\mathbf{\hat x}),
\end{align*}
for any $\mathbf{\hat x}\in \hat K$, $i\in\{1,\ldots,d\}$ and $q_h\in X \in \{\calS_1^h, \bbS_1^h\}$, in conjunction with the identity 
\begin{align*}
    \sum_{j=1}^d [\calA_K^\top]_{m,j} [(\calA_K^\top)^{-1}]_{j,k} = \delta_{m,k},
\end{align*}
where $\delta_{m,k}$ denotes the Kronecker delta, for any $m,k\in\{1,\ldots,d\}$. This shows the desired identity \eqref{eq:Lambda_chain_rule}.

Moreover, for any $i,j\in\{1,\ldots,d\}$, we have with Hölder's inequality, \eqref{eq:shape_regular}, and \eqref{eq:inverse_estimate}, that
\begin{align*}
    \norm{{\mathbf\Lambda} _{i,j}(\bbF_h) - \delta_{i,j} \, \bbF_h }_{L^2(K)}
    &= \nnorm{ \sum\limits_{m=1}^d [(\calA_K^\top)^{-1}]_{i,m} \,
    \big(\hat{\mathbf\Lambda}_{m}(\hat\bbF_h) - \bbF_h\big) \, [\calA_K^\top]_{m,j} }_{L^2(K)}
    \\
    &\leq C \sum\limits_{m=1}^d  
    \nnorm{ \hat{\mathbf\Lambda}_{m}(\hat\bbF_h) - \bbF_h }_{L^2(K)}
    \\
    &\leq C \abs{K}^{1/2} \max_{\ell,m\in\{0,\ldots,d\}} \abs{\bbF_h(\mathbf{p}_\ell^K) - \bbF_h(\mathbf{p}_m^K)}
    \\
    &\leq C h_K \norm{\nabla\bbF_h}_{L^2(K)},
\end{align*}
where $\mathbf{p}_\ell^K$, $\ell\in\{0,\ldots,d\}$, denote the vertices of the simplex $K\in\calT^h$. This gives us the error estimate
\begin{align}
\label{eq:error_Lambda}
    \max_{i,j\in\{1,\ldots,d\}} \norm{{\mathbf\Lambda} _{i,j}(\bbF_h) - \delta_{i,j} \, \bbF_h }_{L^2(K)}
    \leq C h_K \norm{\nabla\bbF_h}_{L^2(K)},
\end{align}
for any $K\in\calT^h$ and all $\bbF_h\in\bbS_1^h$.
The property \eqref{eq:error_Lambda} guarantees that \eqref{eq:convective_convergence} is satisfied in the limit $h\to 0$.

\subsection{Stability and existence of discrete solutions}

We first show a discrete stability result for the numerical scheme \eqref{eq:system_FE}.

\begin{lemma} \label{lemma:stability_FE}
Let $d\in\{2,3\}$ and $k,m\in\bbN$. Let the discrete initial data $\bv_h^0 \in \mathbf{U}_{k+1,\Div}^h$ and $\bbF_h^0 \in \bbS_m^h$ be given. Then, for any $n\in\{1,\ldots,N_T\}$, all solutions 
$(\bv_h^n, p_h^n, \bbF_h^n) \in \mathbf{U}_{k+1}^h \times \calS_k^h \times \bbS_m^h$
to \eqref{eq:system_FE}, if they exist, satisfy
\begin{align}
    \label{eq:stability_FE} \nonumber
    & \frac{\rho}{2 \Delta t} \norm{\bv_h^n}_{L^2}^2 
    + \frac{\rho}{2 \Delta t} 
    \norm{\bv_h^n - \bv_h^{n-1}}_{L^2}^2
    + \frac{\mu}{2\Delta t}  \norm{\bbF_h^n}_{L^2}^2
    + \frac{\mu}{2 \Delta t}  \norm{\bbF_h^n- \bbF_h^{n-1}}_{L^2}^2
    \\ \nonumber
    &\quad
    + \nu \norm{\nabla\bv_h^n}_{L^2}^2
    + \frac{\mu^2}{2\lambda} \norm{\bbF_h^n (\bbF_h^n)^\top}_{L^2}^2
    + \mu \Delta t \norm{\nabla\bbF_h^n}_{L^2}^2
    \\
    &=
    \frac{\rho}{2 \Delta t}  \norm{\bv_h^{n-1}}_{L^2}^2 
    + \frac{\mu}{2\Delta t} \norm{\bbF_h^{n-1}}_{L^2}^2
    + \frac{\mu^2}{2\lambda} \norm{\bbF_h^n}_{L^2}^2.
\end{align}
Moreover, under the assumption $\Delta t< \frac{\lambda}{\mu}$,
one has that
\begin{align}
    \label{eq:stability_FE_gronwall} \nonumber
    & \max_{n=1,\ldots,N_T} \left( 
    \norm{\bv_h^n}_{L^2}^2 
    + \norm{\bbF_h^n}_{L^2}^2 
    \right)
    + \sum_{j=1}^{N_T} \norm{\bv_h^j - \bv_h^{j-1}}_{L^2}^2
    + \sum_{j=1}^{N_T} \norm{\bbF_h^j- \bbF_h^{j-1}}_{L^2}^2
    \\* \nonumber
    &\quad
    + \Delta t \sum_{j=1}^{N_T} \norm{\nabla\bv_h^j}_{L^2}^2
    + \Delta t \sum_{j=1}^{N_T} \norm{\bbF_h^j (\bbF_h^j)^\top}_{L^2}^2
    + (\Delta t)^2 \sum_{j=1}^{N_T} \norm{\nabla\bbF_h^j}_{L^2}^2
    \\*
    &\leq 
    \exp(CT) \left( \norm{\bv_h^0}_{L^2}^2 
    + \norm{\bbF_h^0}_{L^2}^2 \right),
\end{align}
where the constant $C>0$ is independent of $h, \Delta t>0$.
\end{lemma}

\begin{proof}
The proof of \eqref{eq:stability_FE} directly follows from choosing the test functions $\bw_h = \bv_h^n$, $q_h = p_h^n$ and $\bbG_h = \mu \bbF_h^n$ in \eqref{eq:system_FE}, adding the resulting equations, applying \eqref{eq:Lambda_chain_rule} and noting the elementary identity 
\begin{align*}
    2a (a-b) = a^2 - b^2 + (a-b)^2 \quad \forall\, a,b\in\bbR\, .
\end{align*}
After multiplication of \eqref{eq:stability_FE} by $\Delta t$ and summation over $1\leq n \leq m$, where $1\leq m \leq N_T$, it follows that
\begin{align}
    \label{eq:stability_FE_sum} \nonumber
    & \frac{\rho}{2} \norm{\bv_h^m}_{L^2}^2 
    + \sum_{j=1}^m \frac{\rho}{2} \norm{\bv_h^j - \bv_h^{j-1}}_{L^2}^2
    + \frac{\mu}{2}  \norm{\bbF_h^m}_{L^2}^2
    + \sum_{j=1}^m \frac{\mu}{2}  \norm{\bbF_h^j- \bbF_h^{j-1}}_{L^2}^2
    \\ \nonumber
    &\quad
    + \nu \Delta t \sum_{j=1}^m \norm{\nabla\bv_h^j}_{L^2}^2
    + \frac{\mu^2}{2\lambda} \Delta t \sum_{j=1}^m \norm{\bbF_h^j (\bbF_h^j)^\top}_{L^2}^2
    + \mu (\Delta t)^2 \sum_{j=1}^m \norm{\nabla\bbF_h^j}_{L^2}^2
    \\
    &=
    \frac{\rho}{2}  \norm{\bv_h^0}_{L^2}^2 
    + \frac{\mu}{2} \norm{\bbF_h^0}_{L^2}^2
    + \frac{\mu^2}{2\lambda} \Delta t \sum_{j=1}^m \norm{\bbF_h^j}_{L^2}^2.
\end{align}
Under the assumption $\Delta t < \frac{\lambda}{\mu}$, the term $\frac{\mu^2}{2\lambda} \Delta t \norm{\bbF_h^m}_{L^2}^2$ on the right-hand side of \eqref{eq:stability_FE_sum} can be absorbed into the third term on the left-hand side. This allows us to apply a discrete Gronwall argument (cf.~\cite[pp.~401--402]{dahmen_reusken_numerik}), and then the inequality \eqref{eq:stability_FE_gronwall} directly follows.
\end{proof}

\begin{remark}
The ratio $\frac{\lambda}{\mu}$ appearing in Lemma~\ref{lemma:stability_FE} has the physical dimension of time and is commonly referred to as the elastic relaxation time. This parameter describes the characteristic time scale over which elastic stresses decay and is widely used as a reference in benchmark simulations for viscoelastic flows \cite{pimenta_alves_2017_viscoelastic}. Moreover, it plays a central role in defining the Weissenberg number \eqref{eq:nondim_parameters}, which quantifies the ratio of elastic forces to viscous effects.
Physically, accurate transient simulations generally require time steps smaller than the relaxation time. Numerically, the iterative solvers required for the nonlinear system \eqref{eq:system_FE} (e.g., Newton's method) typically demand small time steps to ensure convergence. Consequently, the condition $\Delta t < \frac{\lambda}{\mu}$ in Lemma~\ref{lemma:stability_FE} is not restrictive in practice.

An alternative stability estimate can be derived without invoking the discrete Gronwall lemma. We proceed by applying Hölder's and Young's inequalities to the right-hand side of \eqref{eq:stability_FE} and utilizing \eqref{eq:equivalence_L4}. This yields the bound
\begin{align*}
	\frac{\mu^2}{2\lambda} \norm{\bbF_h^n}_{L^2}^2
	\leq \frac{\mu^2}{2\lambda} \frac{1}{2d} \norm{\bbF_h^n}_{L^4}^4 
	+ \frac{\mu^2}{2\lambda} \frac{d}{2} \abs{\Omega}
	\leq \frac{\mu^2}{4\lambda} \norm{\bbF_h^n (\bbF_h^n)^\top}_{L^2}^2 
		+ \frac{\mu^2 d}{4\lambda} \abs{\Omega}.
\end{align*}
For any $n\in\{1,\ldots,N_T\}$, this implies
\begin{align}
    \label{eq:stability_FE_alternative} \nonumber
    & \frac{\rho}{2 \Delta t} \norm{\bv_h^n}_{L^2}^2 
    + \frac{\rho}{2 \Delta t} 
    \norm{\bv_h^n - \bv_h^{n-1}}_{L^2}^2
    + \frac{\mu}{2\Delta t}  \norm{\bbF_h^n}_{L^2}^2
    + \frac{\mu}{2 \Delta t}  \norm{\bbF_h^n- \bbF_h^{n-1}}_{L^2}^2
    \\ \nonumber
    &\quad
    + \nu \norm{\nabla\bv_h^n}_{L^2}^2
    + \frac{\mu^2}{4\lambda} \norm{\bbF_h^n (\bbF_h^n)^\top}_{L^2}^2
    + \mu \Delta t \norm{\nabla\bbF_h^n}_{L^2}^2
    \\
    &\leq
    \frac{\rho}{2 \Delta t}  \norm{\bv_h^{n-1}}_{L^2}^2 
    + \frac{\mu}{2\Delta t} \norm{\bbF_h^{n-1}}_{L^2}^2
    + \frac{\mu^2 d}{4\lambda} \abs{\Omega}.
\end{align}
Multiplying \eqref{eq:stability_FE_alternative} by $\Delta t$, summing over $1 \leq n \leq m$, and taking the maximum over $1 \leq m \leq N_T$ yields the unconditional estimate, for all $\Delta t > 0$:
\begin{align}
    \label{eq:stability_FE_withoutgronwall} \nonumber
    & \max_{n=1,\ldots,N_T} \left( \norm{\bv_h^n}_{L^2}^2 
    + \norm{\bbF_h^n}_{L^2}^2 \right)
    + \sum_{j=1}^{N_T} \norm{\bv_h^j - \bv_h^{j-1}}_{L^2}^2
    + \sum_{j=1}^{N_T} \norm{\bbF_h^j- \bbF_h^{j-1}}_{L^2}^2
    \\* \nonumber
    &\quad
    + \Delta t \sum_{j=1}^{N_T} \norm{\nabla\bv_h^j}_{L^2}^2
    + \Delta t \sum_{j=1}^{N_T} \norm{\bbF_h^j (\bbF_h^j)^\top}_{L^2}^2
    + (\Delta t)^2 \sum_{j=1}^{N_T} \norm{\nabla\bbF_h^j}_{L^2}^2
    \\*
    &\leq 
    C \left( T + \norm{\bv_h^0}_{L^2}^2 
    + \norm{\bbF_h^0}_{L^2}^2 \right).
\end{align}
Although \eqref{eq:stability_FE_withoutgronwall} holds for all $\Delta t > 0$, it fails to preserve the stationary equilibrium. Even if the system starts at rest ($\bv_h^0 = \mathbf{0}$, $\bbF_h^0 = \mathbb{I}$), the bound on the right-hand side of \eqref{eq:stability_FE_withoutgronwall} grows as an affine function with the final time $T$. Consequently, this estimate permits spurious energy accumulation and cannot guarantee that $\bv_h^n = \mathbf{0}$ and $\bbF_h^n = \mathbb{I}$ for all $n\in\{1,\ldots,N_T\}$, in contrast to the conditional estimate \eqref{eq:stability_FE_gronwall}.
\end{remark}

We recall the following lemma from \cite[Chap.~9.1]{evans_2010} which is a direct consequence of Brouwer's fixed point theorem. This will be useful for proving the existence of a solution to \eqref{eq:system_FE}.

\begin{lemma} \label{lemma:zeros_vector_field}
For $M\in\bbN$, assume that the continuous function $\mathbf{F}\colon\, \bbR^M \to \bbR^M$ satisfies 
\begin{align*}
    \mathbf{F}(\bx) \cdot \bx \geq 0 \quad \forall \, \bx\in\partial B_R(0) \subset\bbR^M,
\end{align*}
for some $R>0$. Then, there exists a point $\bx\in B_R(0)$ such that $\mathbf{F}(\bx) = \mathbf{0}$.
\end{lemma}

We now combine Lemma~\ref{lemma:zeros_vector_field} with the estimates from Lemma~\ref{lemma:stability_FE} to prove the existence of a solution to \eqref{eq:system_FE}.

\begin{theorem} \label{theorem:existence_FE}
Let $d\in\{2,3\}$ and $k,m\in\bbN$. Let the discrete initial data $\bv_h^0 \in \mathbf{U}_{k+1,\Div}^h$ and $\bbF_h^0 \in \bbS_m^h$ be given. In addition, suppose that $\Delta t < \frac\lambda\mu$ is satisfied. Then, for any $n\in\{1,\ldots,N_T\}$, there exists a solution
$(\bv_h^n, p_h^n, \bbF_h^n) \in \mathbf{U}_{k+1}^h \times \calS_k^h \times \bbS_m^h$
to \eqref{eq:system_FE} and it satisfies \eqref{eq:stability_FE} and \eqref{eq:stability_FE_gronwall}.
\end{theorem}

\begin{proof}
We prove the existence of solutions to \eqref{eq:system_FE} by combining the inequality \eqref{eq:stability_FE} and Lemma~\ref{lemma:zeros_vector_field}. 
However, we can not directly show the existence of solutions $(\bv_h^n, p_h^n, \bbF_h^n) \in \mathbf{U}_{k+1}^h \times \calS_k^h \times \bbS_m^h$, since we have no \textit{a priori} estimates for the pressure $p_h^n$. 
Therefore, we first prove the existence of functions $(\bv_h^n, \bbF_h^n) \in \mathbf{U}_{k+1,\Div}^h \times \bbS_m^h$ solving \eqref{eq:system_FE} with \eqref{eq:v_FE} replaced by
\begin{align}
    \label{eq:v_FE2} \nonumber
    0 &= 
    \frac{\rho}{\Delta t}\skp{\bv_h^n - \bv_h^{n-1}}{\bw_h}_{L^2}
    + \frac\rho{2} \skp{(\bv_h^{n-1} \cdot \nabla) \bv_h^n}{\bw_h}_{L^2}
    - \frac\rho{2} \skp{\bv_h^n}{(\bv_h^{n-1} \cdot \nabla) \bw_h}_{L^2}
    \\*
    &\quad
    + \nu \skp{\nabla \bv_h^n}{\nabla\bw_h}_{L^2}
    + \mu \skp{\bbF_{h}^{n} (\bbF_h^n)^\top}{\nabla\bw_h}_{L^2} ,
\end{align}
for all $\bw\in\mathbf{U}_{k+1,\Div}^h$.
Afterwards, we reconstruct the discrete pressure $p_h^n\in\calS_k^h$ using the inf-sup inequality \eqref{eq:LBB} and hence show that $\bv_h^n, p_h^n, \bbF_h^n$ solve \eqref{eq:system_FE}.

We define a continuous function $\mathbf{F}\colon\, \bbR^M \to \bbR^M$ with $M\coloneqq \dim(\mathbf{U}_{k+1,\Div}^h) +  \dim(\bbS_m^h)$, mapping the unknown coefficients $\bx\in\bbR^M$ of $(\bv_h^n, \mu \bbF_h^n) \in \mathbf{U}_{k+1,\Div}^h \times \bbS_m^h$ to the right-hand sides of \eqref{eq:v_FE2}, \eqref{eq:F_FE}. Then, a zero $\bx\in\bbR^M$ of $\mathbf{F}$ is equivalent to a solution $(\bv_h^n, \bbF_h^n) \in \mathbf{U}_{k+1,\Div}^h \times \bbS_m^h$ of \eqref{eq:v_FE2} and \eqref{eq:F_FE}.
Our aim is to show $\mathbf{F}(\bx) \cdot \bx \geq c_1 \abs{\bx}^2 - c_2$ with $\abs{\bx} = R$, for some $R>0$ and some constants $c_1,c_2\in\bbR$ that are independent of $\bv_h^n$ and $\bbF_h^n$. We have similarly to the proof of Lemma~\ref{lemma:stability_FE} that
\begin{align*}
    \mathbf{F}(\bx) \cdot \bx
    &= \frac{\rho}{2 \Delta t} \norm{\bv_h^n}_{L^2}^2 
    + \frac{\rho}{2 \Delta t} 
    \norm{\bv_h^n - \bv_h^{n-1}}_{L^2}^2
    + \frac{\mu}{2\Delta t}  \norm{\bbF_h^n}_{L^2}^2
    + \frac{\mu}{2 \Delta t}  \norm{\bbF_h^n- \bbF_h^{n-1}}_{L^2}^2
    \\ 
    &\quad
    + \nu \norm{\nabla\bv_h^n}_{L^2}^2
    + \frac{\mu^2}{2\lambda} \norm{\bbF_h^n (\bbF_h^n)^\top}_{L^2}^2
    + \mu \Delta t \norm{\nabla\bbF_h^n}_{L^2}^2
    \\
    &\quad 
    - \frac{\rho}{2 \Delta t}  \norm{\bv_h^{n-1}}_{L^2}^2 
    - \frac{\mu}{2\Delta t} \norm{\bbF_h^{n-1}}_{L^2}^2
    - \frac{\mu^2}{2\lambda} \norm{\bbF_h^n}_{L^2}^2
    \\
    &\geq \frac{\rho}{2 \Delta t} \norm{\bv_h^n}_{L^2}^2 
    + \frac\mu{2} \left(\frac{1}{\Delta t} - \frac{\mu}{\lambda}  \right)  \norm{\bbF_h^n}_{L^2}^2
    - c_2,
\end{align*}
for a positive constant $c_2 := \frac\rho{2 \Delta t}  \norm{\bv_h^{n-1}}_{L^2}^2 + \frac\mu{2\Delta t} \norm{\bbF_h^{n-1}}_{L^2}^2$ that depends only on $\bv_h^{n-1}$, $\bbF_h^{n-1}$ but not on $\bv_h^n$, $\bbF_h^n$.
Under the condition $\Delta t < \frac{\lambda}{\mu}$, and thanks to norm equivalence in finite-dimensional normed linear spaces, we have
\begin{align*}
    \mathbf{F}(\bx) \cdot \bx 
    \geq c_1 \abs{\bx}^2 - c_2,
\end{align*}
for some $c_1>0$ that is independent of $\bv_h^n$, $\bbF_h^n$. Hence, for $R>0$ large enough, we have $\mathbf{F}(\bx) \cdot \bx > 0$. Then, the existence of a zero of $\mathbf{F}$ follows from Lemma~\ref{lemma:zeros_vector_field}. This justifies the existence of a solution $(\bv_h^n, \bbF_h^n) \in \mathbf{U}_{k+1,\Div}^h \times \bbS_m^h$ to \eqref{eq:v_FE2}, \eqref{eq:F_FE}.

We note that \eqref{eq:v_FE} defines a linear functional $\mathbf{U}_{k+1}^h \to \bbR$ that vanishes on $\mathbf{U}_{k+1,\Div}^h$. Thus, the existence of a pressure $p_h^n \in \calS_k^h$, which is unique up to an additive constant, follows directly from \cite[Lem.~4.1]{girault_raviart_1986} since the inf-sup inequality \eqref{eq:LBB} is valid for the lowest order Taylor--Hood element $(\mathbf{U}_{k+1}^h, \calS_k^h)$. This justifies the existence of a solution $(\bv_h^n, p_h^n, \bbF_h^n) \in \mathbf{U}_{k+1}^h \times \calS_k^h \times \bbS_m^h$ to \eqref{eq:system_FE}, which satisfies \eqref{eq:stability_FE} and \eqref{eq:stability_FE_gronwall} thanks to Lemma~\ref{lemma:stability_FE}.
\end{proof}

\subsection{(Non-)Positivity of the determinant}
A major drawback of formulating the viscoelastic Giesekus model in terms of the deformation gradient $\bbF$, rather than the symmetric tensor $\bbB = \bbF\bbF^\top$ (see \eqref{eq:B}), is that it may not be possible to ensure the positivity of $\det(\bbF)$ at the fully discrete or time-discrete level. We now explain the underlying reason.

Although the scalar function $s\mapsto-\ln(s)=-\frac12\ln(s^2)$ is convex, its matrix-valued counterpart 
\begin{align*}
    \bbF\mapsto-\frac12\trace\ln(\bbF\bbF^\top) = -\frac12 \ln\det(\bbF\bbF^\top)
\end{align*}
is not convex in $\bbF \in\bbR^{d\times d}$, because $\bbF$ is not symmetric \cite[Theorem~5.1]{ball_1976}. However, the convexity of this mapping is essential for showing the integrablity of $\ln\det(\bbF\bbF^\top)$. 
In fact, using a diagonalization argument and the convexity of $s\mapsto-\ln(s)=-\frac12\ln(s^2)$ for all $s>0$, one can show for all symmetric and regular matrices $\mathbb{A},\mathbb{B} \in\bbR^{d\times d}$, that (cf.~\cite{barrett_boyaval_2009})
\begin{align}
    \label{eq:lndet_convex}
    (\bbA-\bbB) : (-\bbA^{-\top}) \geq -\frac12\ln\det(\bbA \bbA^\top) + \frac12\ln\det(\bbB \bbB^\top),
\end{align}
which, in general, does not hold for non-symmetric matrices.
A simple counter-example for non-symmetric matrices is provided by the choice
$\bbA = \left(\begin{smallmatrix}1 & 1 \\ 0 & 1\end{smallmatrix}\right)$
and 
$\bbB = \left(\begin{smallmatrix}2 & 0 \\ 1 & 2\end{smallmatrix}\right)$.
%

However, an inequality like \eqref{eq:lndet_convex} would be crucial to guarantee the positivity of $\det(\bbF)$ at the fully discrete or time-discrete level via energy estimates. Nevertheless, by first letting $h \to 0$ (for fixed $\Delta t$) and then $\Delta t \to 0$ in \eqref{eq:system_FE}, we can prove the subsequence convergence to a weak solution that indeed satisfies $\det(\bbF) > 0$, provided that the initial datum fulfills $\det(\bbF_0) > 0$ (see Theorem~\ref{theorem:convergence}).

If ensuring the positivity of $\det(\bbF)$ (or alternatively $\det(\bbB)$) is essential, one could instead discretize the model using linear finite elements for $\bbB = \bbF\bbF^\top$, rather than applying \eqref{eq:F_FE} directly to the deformation gradient $\bbF$. This would correspond to adapting the numerical approach developed for the Oldroyd-B model in \cite{barrett_boyaval_2009} to the present setting. The resulting discretization would still resemble \eqref{eq:F_FE}, but would require special numerical quadrature for technical reasons, as well as a convective term defined in the spirit of \eqref{eq:convective_Lambda}, though with a nonlinear definition instead (cf.~\cite{barrett_boyaval_2009}). 
The main benefit would be a discrete energy law similar to \eqref{eq:stability}, including the logarithmic part $\ln\det(\bbB)$ in contrast to \eqref{eq:stability_FE}. This would guarantee the positive definiteness of $\bbB$ at the fully discrete level.
However, in such a formulation, the convergence analysis presented in this work would most likely not be possible without adding an extra diffusive regularization term $\varepsilon \Delta \bbB$ to \eqref{eq:B} with fixed $\varepsilon > 0$.

\subsection{A fully linear and stable numerical approximation} \label{sec:numerical_scheme_linear}
So far, we have established existence and stability of discrete solutions to \eqref{eq:system_FE} in both two and three dimensions. In Section~\ref{sec:proof_2d}, we will first pass to the limit $h\to 0$ and then $\Delta t\to 0$ to prove (subsequence) convergence to a weak solution in two dimensions, as defined in Definition~\ref{def:weak_solution}.
Before that, we briefly discuss a stable numerical approximation that is fully linear (in contrast to \eqref{eq:system_FE}). This may be useful for computational experiments, although we are not able to prove its convergence.

Let $k,m\in\bbN$. Given the discrete initial data 
\[\bv_h^0 \in \mathbf{U}_{k+1,\Div}^h \quad \text{and} \quad
\bbF_h^0 \in \bbS_m^h,\] 
for $n\in\{1,\ldots,N_T\}$ we aim to find a solution 
\[(\bv_h^n, p_h^n, \bbF_h^n) \in \mathbf{U}_{k+1}^h \times \calS_k^h \times \bbS_m^h\] 
such that, for all test functions $(\bw_h, q_h, \bbG_h) \in \mathbf{U}_{k+1}^h \times \calS_k^h \times \bbS_m^h$,
\begin{subequations}
\label{eq:system_FE_k_linear}
\begin{align}
    \label{eq:v_FE_k_linear}
    \nonumber
    0 &= 
    \frac\rho{\Delta t}\skp{\bv_h^n - \bv_h^{n-1}}{\bw_h}_{L^2}
    + \frac\rho{2} \skp{(\bv_h^{n-1} \cdot \nabla) \bv_h^n}{\bw_h}_{L^2}
    - \frac\rho{2} \skp{\bv_h^n}{(\bv_h^{n-1} \cdot \nabla) \bw_h}_{L^2}
    \\
    &\quad
    + \nu \skp{\nabla \bv_h^n}{\nabla\bw_h}_{L^2}
    - \skp{p_h^n}{\Div\bw_h}_{L^2}
    + \mu \skp{\bbF_h^{n-1} (\bbF_h^{n})^\top }{\nabla\bw_h}_{L^2} ,
    \\
    \label{eq:div_FE_k_linear}
    0 &= \skp{\Div\bv_h^n}{q_h}_{L^2} ,
    \\
    \label{eq:F_FE_k_linear} \nonumber
    0 &= 
    \frac{1}{\Delta t} \skp{\bbF_h^n - \bbF_h^{n-1}}{\bbG_h}_{L^2}
    + \bc_h(\bv_h^{n-1}, \bbF_h^n, \bbG_h)
    + \frac{\mu}{2\lambda} \skp{ (\bbF_h^{n-1} (\bbF_h^{n})^\top \bbF_h^{n-1} - \bbF_h^n} {\bbG_h}_{L^2} 
    \\ 
    &\quad 
    - \skp{ (\nabla\bv_h^n) \bbF_h^{n-1} }{\bbG_h}_{L^2}
    + \Delta t \skp{\nabla\bbF_h^n}{\nabla\bbG_h}_{L^2},
\end{align}
\end{subequations}
where $\bc_h(\cdot,\cdot,\cdot)$ can be given by \eqref{eq:convective_skewsymm}, \eqref{eq:convective_higherorder} or \eqref{eq:convective_Lambda}.

As before, we can show the following existence and stability result with similar arguments. Note that this scheme is fully linear and the number of equations matches the number of unknowns, so the existence of a solution is equivalent to its uniqueness.

\begin{theorem} \label{theorem:existence_FE_k_linear}
Let $d\in\{2,3\}$ and let $k,m\in\bbN$ . Let the discrete initial data $\bv_h^0 \in \mathbf{U}_{k+1,\Div}^h$ and $\bbF_h^0 \in \bbS_{m}^h$ be given. In addition, suppose that $\Delta t < \frac{\lambda}{\mu}$ is satisfied. Then, for any $n\in\{1,\ldots,N_T\}$, there exists a unique solution
$(\bv_h^n, p_h^n, \bbF_h^n) \in \mathbf{U}_{k+1}^h \times \calS_k^h \times \bbS_m^h$
to \eqref{eq:system_FE_k_linear} with $\int_\Omega p_h^n = 0$, and we have that
\begin{align}
    \label{eq:stability_FE_k_linear_gronwall} \nonumber
    & \max_{n=1,\ldots,N_T} \left( 
    \norm{\bv_h^n}_{L^2}^2 
    + \norm{\bbF_h^n}_{L^2}^2 
    \right)
    + \sum_{j=1}^{N_T} \norm{\bv_h^j - \bv_h^{j-1}}_{L^2}^2
    + \sum_{j=1}^{N_T} \norm{\bbF_h^j- \bbF_h^{j-1}}_{L^2}^2
    \\ \nonumber
    &\quad
    + \Delta t \sum_{j=1}^{N_T} \norm{\nabla\bv_h^j}_{L^2}^2
    + \Delta t \sum_{j=1}^{N_T} \norm{\bbF_h^{j-1} (\bbF_h^j)^\top}_{L^2}^2
    + (\Delta t)^2 \sum_{j=1}^{N_T} \norm{\nabla\bbF_h^j}_{L^2}^2
    \\
    &\leq 
    \exp(CT) \left( \norm{\bv_h^0}_{L^2}^2 
    + \norm{\bbF_h^0}_{L^2}^2 \right),
\end{align}
where the constant $C>0$ is independent of $h, \Delta t>0$.
\end{theorem}

It remains unclear whether a convergence analysis analogous to Section~\ref{sec:proof_2d} can be carried out for System~\eqref{eq:system_FE_k_linear}.
In particular, identifying the weak limit of $|\bbF_h^{n-1} (\bbF_h^{n})^\top|^2$ with $|\bbF\bbF^\top|^2$ seems challenging, or even impossible, without a uniform bound on the discrete deformation tensor in $L^4(\Omega_T;\bbR^{d\times d})$. 
A possible way to overcome this difficulty would be to replace the stabilization term $\Delta t \skp{\nabla \bbF_h^n}{\nabla \bbG_h}_{L^2}$ in \eqref{eq:F_FE_k_linear} by $\varepsilon \skp{\nabla \bbF_h^n}{\nabla \bbG_h}_{L^2}$ with fixed $\varepsilon > 0$. However, this modification would naturally lead to limiting equations with stress diffusion.
Despite this drawback, the system~\eqref{eq:system_FE_k_linear} is fully linear and may still be of practical interest for computational purposes.


%
%
%
\section{Convergence to a weak solution in 2D} \label{sec:proof_2d}
As outlined earlier, the goal of this section is to first take the limit $h \to 0$ in \eqref{eq:system_FE}, and then let $\Delta t \to 0$. This process allows us to adapt the arguments from \cite{bulicek_2022_giesekus_2d} to the time-discrete setting.
The limit passage $h \to 0$ for fixed $\Delta t>0$ is presented in Section~\ref{sec:time_discrete}. The main result for the limit $\Delta t\to 0$ is stated in Theorem~\ref{theorem:convergence} and is proved in detail in Sections~\ref{sec:limit_vanishing_dt}--\ref{sec:positivity_F}.

\subsection{Limit passage in space} \label{sec:time_discrete}

We now introduce a time-discrete system in two and three dimensions, $d\in\{2,3\}$.

Let $\bv^0 \in H^1_{0,\Div}(\Omega;\bbR^d)$ and $\bbF^0 \in L^2(\Omega;\bbR^{d\times d})$ be given initial data acting as approximations of $\bv_0 \in L^2_{\Div}(\Omega;\bbR^d)$ and $\bbF_0 \in L^2(\Omega;\bbR^{d\times d})$ defined in \eqref{eq:system_init}.
Then, for $n\in\{1,\ldots,N_T\}$, we aim to find a time-discrete solution 
\[(\bv^n, \bbF^n) \in H^1_{0,\Div}(\Omega;\bbR^d) \times H^1(\Omega;\bbR^{d\times d})\] 
such that, for all test functions $(\bw, \bbG) \in H^1_{0,\Div}(\Omega;\bbR^d) \times H^1(\Omega;\bbR^{d\times d})$,
\begin{subequations}
\label{eq:system_time_discrete}
\begin{align}
    \nonumber
    \label{eq:v_time_discrete}
    0 &= 
    \frac{\rho}{\Delta t}\skp{\bv^n - \bv^{n-1}}{\bw}_{L^2}
    - \rho \skp{\bv^n \otimes \bv^{n-1}}{\nabla \bw}_{L^2}
    + \nu \skp{\nabla \bv^n}{\nabla\bw}_{L^2}
    \\*
    &\quad
    + \mu \skp{\bbF^n (\bbF^n)^\top }{\nabla\bw}_{L^2} ,
    \\
    \label{eq:F_time_discrete}
    \nonumber
    0 &= 
    \frac{1}{\Delta t} \skp{\bbF^n - \bbF^{n-1}}{\bbG}_{L^2}
    - \skp{\bbF^n \otimes \bv^{n-1}}{\nabla \bbG}_{L^2}
    + \frac{\mu}{2\lambda} \skp{ \bbF^n (\bbF^n)^\top \bbF^n - \bbF^n} {\bbG}_{L^2} 
    \\*
    &\quad 
    - \skp{ (\nabla\bv^n) \bbF^n }{\bbG}_{L^2}
    + \Delta t \skp{\nabla\bbF^n}{\nabla\bbG}_{L^2}.
\end{align}
\end{subequations}

The next theorem is proved by passing to the limit in \eqref{eq:system_FE} as $h \to 0$, while keeping the time step size $\Delta t > 0$ (and hence the time index $n \in \{1,\ldots,N_T\}$) fixed.

\begin{theorem} \label{theorem:existence_time_discrete}
Let $d\in\{2,3\}$.
Let $\bv^0 \in H^1_{0,\Div}(\Omega;\bbR^d)$ and $\bbF^0 \in L^2(\Omega;\bbR^{d\times d})$ be given initial data. In addition, suppose that $\Delta t < \frac{\mu}{\lambda}$.
Then, for $n\in\{1,\ldots,N_T\}$, there exists a time-discrete solution 
\[(\bv^n, \bbF^n) \in H^1_{0,\Div}(\Omega;\bbR^d) \times H^1(\Omega;\bbR^{d\times d})\] 
to \eqref{eq:system_time_discrete}. Moreover, we have that
\begin{align}
    \label{eq:stability_time_discrete_gronwall} \nonumber
    & \max_{n=1,\ldots,N_T} \left(
     \norm{\bv^n}_{L^2}^2 
    + \norm{\bbF^n}_{L^2}^2 
    \right)
    + \sum_{j=1}^{N_T} \norm{\bv^j - \bv^{j-1}}_{L^2}^2
    + \sum_{j=1}^{N_T} \norm{\bbF^j- \bbF^{j-1}}_{L^2}^2
    \\ \nonumber
    &\quad
    + \Delta t \sum_{j=1}^{N_T} \norm{\nabla\bv^j}_{L^2}^2
    + \Delta t \sum_{j=1}^{N_T} \norm{\bbF^j (\bbF^j)^\top}_{L^2}^2
    + (\Delta t)^2 \sum_{j=1}^{N_T} \norm{\nabla\bbF^j}_{L^2}^2
    \\
    &\leq 
    \exp(CT) \left( \norm{\bv^0}_{L^2}^2 
    + \norm{\bbF^0}_{L^2}^2 \right),
\end{align}
where the constant $C>0$ is independent of $\Delta t > 0$.
\end{theorem}

\begin{proof}
To establish the existence of a solution to \eqref{eq:system_time_discrete}, we proceed as follows.
We select the discrete initial data $\bv_h^0 \in \mathbf{U}^h_{k+1,\Div}$ and $\bbF_h^0 \in \bbS_m^h$ such that $\norm{\bv_h^0 - \bv^0}_{L^2}\to 0$ and $\norm{\bbF_h^0 - \bbF^0}_{L^2} \to 0$, as $h\to 0$.
This can be achieved by, e.g., standard projections of $\bv^0$, $\bbF^0$ to the respective finite element spaces.

For the limit passage in \eqref{eq:system_FE}, we consider arbitrary test functions $\bw \in C^\infty_c(\Omega;\bbR^d)$ with $\Div(\bw)=0$ and $\bbG\in C^\infty(\overline{\Omega};\bbR^d)$. We then approximate them by discrete test functions $\bw_h \in \mathbf{U}^h_{k+1,\Div}$ and $\bbG_h\in \bbS_m^h$ such that $(\bw_h,\bbG_h) \to (\bw,\bbG)$ strongly in $H^1_{0,\Div}(\Omega;\bbR^d) \times W^{1,\infty}(\Omega;\bbR^{d\times d})$ and uniformly on $\overline{\Omega}$. 

For each time step $n\in\{1,\ldots,N_T\}$ and fixed $\Delta t>0$, spatial compactness follows from the energy identity \eqref{eq:stability_FE}. In particular, we obtain for a (non-relabeled) subsequence:
\begin{alignat*}{3}
    \bv_h^n &\rightharpoonup \bv^n \quad &&\text{weakly } \quad 
    &&\text{in } H^1_{0,\Div}(\Omega;\bbR^d),
    \\*
    \bv_h^n &\to \bv^n \quad &&\text{strongly } \quad 
    &&\text{in } L^r_{\Div}(\Omega;\bbR^d),
    \\
    \bbF_h^n &\rightharpoonup \bbF^n \quad &&\text{weakly } \quad 
    &&\text{in } H^1(\Omega;\bbR^{d\times d}),
    \\*
    \bbF_h^n &\to \bbF^n \quad &&\text{strongly } \quad 
    &&\text{in } L^s(\Omega;\bbR^{d\times d}),
\end{alignat*}
for any $r,s\in[1,\frac{2d}{d-2})$.
%
Together with \eqref{eq:convective_convergence}, these convergence results are sufficient to pass to the limit in \eqref{eq:system_FE} and establish the existence of at least one solution to \eqref{eq:system_time_discrete}. Moreover, the specific form of the convective term in \eqref{eq:v_time_discrete} follows from integration by parts over $\Omega$. In addition, a density argument yields the precise formulation of \eqref{eq:system_time_discrete} with the correct test functions.

Finally, applying weak lower semicontinuity to \eqref{eq:stability_FE_gronwall} yields \eqref{eq:stability_time_discrete_gronwall}.
\end{proof}


%
%
%
We also present an estimate for the discrete time derivatives of $\bv^n$ and $\bbF^n$, which will be useful for the limit process $\Delta t\to 0$.
\begin{lemma}
Let $d\in\{2,3\}$. Assume that $\bv^0 \in H^1_{0,\Div}(\Omega;\bbR^d)$ and $\bbF^0 \in L^2(\Omega;\bbR^{d\times d})$ are given, and suppose that $\Delta t < \frac{\lambda}{\mu}$. Then, in addition to \eqref{eq:stability_time_discrete_gronwall}, we have that
\begin{align}
    \label{eq:time_derivatives_discrete} \nonumber
    &\Delta t \sum_{n=1}^{N_T} \nnorm{\frac{\bv^n - \bv^{n-1}}{\Delta t}}_{(H^1_{0,\Div})'}^{4/d}
    + \Delta t \sum_{n=1}^{N_T} \nnorm{\frac{\bbF^n - \bbF^{n-1}}{\Delta t}}_{(H^1)'}^{4/3} 
    \\*
    &\leq C(T) \left( 1 + \norm{\bv^0}_{L^2}^2 
    + \norm{\bbF^0}_{L^2}^2 
    + \Delta t \norm{\nabla\bv^0}_{L^2}^2 \right) .
\end{align}
Here, the constant $C>0$ depends on $T$ but not on $\Delta t$.
\end{lemma}

\begin{proof}
We only show the estimate for $d=3$. The result for $d=2$ can be shown in a similar manner.
For any $\bw\in H^1_{0,\Div}(\Omega;\bbR^d)$, we have
\begin{align*}
    \abs{\tfrac{1}{\Delta t}\skp{\bv^n - \bv^{n-1}}{\bw}_{L^2}}
    &\leq 
    \abs{\skp{\bv^n \otimes \bv^{n-1}}{\nabla \bw}_{L^2}}
    + \frac{\nu}{\rho} \abs{\skp{\nabla \bv^n}{\nabla\bw}_{L^2}}
    + \frac{\mu}{\rho} \abs{\skp{\bbF^n (\bbF^n)^\top }{\nabla\bw}_{L^2}}
    \\
    &\leq C \left( \norm{\bv^n}_{L^3} \norm{\bv^{n-1}}_{L^6} 
    + \norm{\bv^n}_{H^1} 
    + \norm{\bbF^n (\bbF^n)^\top}_{L^2} \right) \norm{\bw}_{H^1},
\end{align*}
which implies that
\begin{align*}
    \nnorm{\tfrac{1}{\Delta t}(\bv^n - \bv^{n-1})}_{(H^1_{0,\Div})'}
    &\leq C \left( \norm{\bv^n}_{L^3} \norm{\bv^{n-1}}_{L^6} 
    + \norm{\bv^n}_{H^1} 
    + \norm{\bbF^n (\bbF^n)^\top}_{L^2} \right)
    \\
    &\leq C \left( \norm{\bv^n}_{L^2}^{1/2} \norm{\bv^n}_{H^1}^{1/2} \norm{\bv^{n-1}}_{H^1} 
    + \norm{\bv^n}_{H^1} 
    + \norm{\bbF^n (\bbF^n)^\top}_{L^2} \right),
\end{align*}
where we have used the Gagliardo--Nirenberg inequality and the continuous Sobolev embedding $H^1(\Omega) \hookrightarrow L^6(\Omega)$, for $d\in\{2,3\}$, in the last step.
Taking the $4/3$ power on both sides, multiplying by $\Delta t$, summing over all $n\in\{1,\ldots,N_T\}$, and applying Hölder's and Young's inequalities leads to
\begin{align*}
    &\Delta t \sum_{n=1}^{N_T} \nnorm{\tfrac{1}{\Delta t}(\bv^n - \bv^{n-1})}_{(H^1_{0,\Div})'}^{4/3}
    \\*
    &\leq C \Delta t \sum_{n=1}^{N_T} \left( \norm{\bv^n}_{L^2}^{1/2} \norm{\bv^n}_{H^1}^{1/2} \norm{\bv^{n-1}}_{H^1} 
    + \norm{\bv^n}_{H^1} 
    + \norm{\bbF^n (\bbF^n)^\top}_{L^2} \right)^{4/3} 
    \\
    &\leq C \left( \max_{n=1,\ldots,N_T} \norm{\bv^n}_{L^2}^2 \right)^{1/3}
    \left( \Delta t \sum_{n=1}^{N_T} \norm{\bv^n}_{H^1}^2 \right)^{1/3}
    \left( \Delta t \sum_{n=1}^{N_T} \norm{\bv^{n-1}}_{H^1}^2 \right)^{2/3}
    \\*
    &\quad
    + C T^{1/3} \left( \Delta t \sum_{n=1}^{N_T} \norm{\bv^n}_{H^1}^2 \right)^{2/3}
    + C T^{1/3} \left( \Delta t \sum_{n=1}^{N_T} \norm{\bbF^n (\bbF^n)^\top}_{L^2}^2 \right)^{2/3}
    \\
    &\leq C(T) \left( 1 + \norm{\bv^0}_{L^2}^2 
    + \norm{\bbF^0}_{L^2}^2 
    + \Delta t \norm{\nabla\bv^0}_{L^2}^2 \right).
\end{align*}

Similarly, for any $\bbG \in H^1(\Omega;\bbR^{d\times d})$, we have  
\begin{align*}
    &\abs{\tfrac{1}{\Delta t} \skp{\bbF^n - \bbF^{n-1}}{\bbG}_{L^2} }
    \\
    &\leq \abs{\skp{\bbF^n \otimes \bv^{n-1}}{\nabla \bbG}_{L^2}}
    + \frac{\mu}{2\lambda} \abs{\skp{ \bbF^n (\bbF^n)^\top \bbF^n - \bbF^n} {\bbG}_{L^2}}
    \\*
    &\quad
    + \abs{\skp{ (\nabla\bv^n) \bbF^n }{\bbG}_{L^2}}
    + \Delta t \abs{\skp{\nabla\bbF^n}{\nabla\bbG}_{L^2}}
    \\
    &\leq C \left( \norm{\bbF^n}_{L^4} \norm{\bv^{n-1}}_{L^4} 
    + \norm{\bbF^n}_{L^4}^3 
    + \norm{\bbF^n}_{L^2}
    + \norm{\nabla\bv^n}_{L^2} \norm{\bbF^n}_{L^4} 
    + \Delta t \norm{\nabla\bbF^n}_{L^2} \right) \norm{\bbG}_{H^1},
\end{align*}
which implies the inequalities
\begin{align*}
    &\norm{\frac1{\Delta t}(\bbF^n - \bbF^{n-1})}_{(H^1)'}
    \\
    &\leq C \left( \norm{\bbF^n}_{L^4} \norm{\bv^{n-1}}_{L^4} 
    + \norm{\bbF^n}_{L^4}^3 
    + \norm{\bbF^n}_{L^2}
    + \norm{\nabla\bv^n}_{L^2} \norm{\bbF^n}_{L^4} 
    + \Delta t \norm{\nabla\bbF^n}_{L^2} \right)
    \\
    &\leq C \left( \norm{\bbF^n}_{L^4} \norm{\bv^{n-1}}_{H^1} 
    + \norm{\bbF^n}_{L^4}^3 
    + \norm{\bbF^n}_{L^2}
    + \norm{\bv^n}_{H^1} \norm{\bbF^n}_{L^4} 
    + \Delta t \norm{\nabla\bbF^n}_{L^2} \right) .
\end{align*}
Again, taking the $4/3$ power on both sides, multiplying both sides by $\Delta t$, summing over all $n\in\{1,\ldots,N_T\}$, and applying Hölder's and Young's inequalities give us
\begin{align*}
    &\Delta t \sum_{n=1}^{N_T} \norm{\frac1{\Delta t}(\bbF^n - \bbF^{n-1})}_{(H^1)'}^{4/3}
    \\*
    &\leq C \Delta t \sum_{n=1}^{N_T} \left( \norm{\bbF^n}_{L^4} \norm{\bv^{n-1}}_{H^1} 
    + \norm{\bbF^n}_{L^4}^3 
    + \norm{\bbF^n}_{L^2}
    + \norm{\bv^n}_{H^1} \norm{\bbF^n}_{L^4} 
    + \Delta t \norm{\nabla\bbF^n}_{L^2} \right)^{4/3}
    \\
    &\leq C \left( \Delta t \sum_{n=1}^{N_T} \norm{\bbF^n}_{L^4}^4  \right)^{1/3}
    \left( \Delta t \sum_{n=0}^{N_T}  \norm{\bv^n}_{H^1}^2 \right)^{2/3}
    + C \Delta t \sum_{n=1}^{N_T}  \norm{\bbF^n}_{L^4}^4
    \\*
    &\quad
    + C T^{1/3} \left( \Delta t \sum_{n=1}^{N_T} \norm{\bbF^n}_{L^2}^2  \right)^{2/3}
    + C T^{1/3} \left( (\Delta t)^2 \sum_{n=1}^{N_T} \norm{\nabla \bbF^n}_{L^2}^2  \right)^{2/3} 
    \\
    &\leq C(T) \left( 1 + \norm{\bv^0}_{L^2}^2 
    + \norm{\bbF^0}_{L^2}^2 
    + \Delta t \norm{\nabla\bv^0}_{L^2}^2 \right).
\end{align*}
This proves the lemma.
\end{proof}

\subsection{Interpolation in time}
We introduce the following notation for piecewise affine linear and piecewise constant extensions of time-discrete functions $a^n(\cdot)$, $n=0,...,N_T$:
\begin{alignat}{2}
    \label{eq:fun_Delta_t}
    a^{\Delta t}(\cdot, t) 
    &\coloneqq 
    \frac{t - t^{n-1}}{\Delta t} a^n(\cdot)
    + \frac{t^n - t}{\Delta t} a^{n-1}(\cdot)
    \quad\quad 
    && \mbox{for } t\in [t^{n-1},t^n], \, n\in \{1,...,N_T\},
    \\
    \label{eq:fun_Delta_t_pm}
    a^{\Delta t,+}(\cdot, t) 
    &\coloneqq a^n(\cdot),
    \quad\quad 
    a^{\Delta t,-}(\cdot, t) 
    \coloneqq a^{n-1}(\cdot)
    \quad\quad 
    && \mbox{for } t\in (t^{n-1},t^n], \, n\in \{1,...,N_T\}.
\end{alignat}
Note that we write $a^{\Delta t,\pm}$ for results that hold for both $a^{\Delta t,+}$ and $a^{\Delta t,-}$, and we write $a^{\Delta t(,\pm)}$ for results that are valid for $a^{\Delta t}$, $a^{\Delta t,+}$, and $a^{\Delta t,-}$, respectively. 

With this notation, we reformulate \eqref{eq:system_time_discrete} by interpolating between time levels.
In fact, for all $n\in\{1,\ldots,N_T\}$, we choose the test functions 
\[\bw^n = \frac{1}{\Delta t} \int_{t^{n-1}}^{t^n} \tilde{\bw}(\cdot,t)\dt\]
and 
\[\bbG^n = \frac{1}{\Delta t} \int_{t^{n-1}}^{t^n} \tilde{\bbG}(\cdot,t) \dt,\]
respectively, where $\tilde{\bw} \in C_c^\infty(\bbR;H^1_{0,\Div}(\Omega;\bbR^d))$ and $\tilde{\bbG} \in C_c^\infty(\bbR;H^1(\Omega;\bbR^{d\times d}))$. Then, we multiply each equation by $\Delta t$ and sum over $n\in\{1,\ldots,N_T\}$. 
%
%
%
Finally, we end up with the following system:
\begin{subequations}
\label{eq:system_time_continuous}
\begin{align}
    \nonumber
    \label{eq:v_time_continuous}
    0 &= 
    \rho \int_0^T \skp{\partial_t \bv^{\Delta t}}{\tilde{\bw}}_{L^2} \dt
    - \rho \int_0^T 
    \skp{\bv^{\Delta t,+} \otimes \bv^{\Delta t,-}}{\nabla\tilde{\bw}}_{L^2} \dt
    \\* 
    &\quad 
    + \nu \int_0^T \skp{\nabla \bv^{\Delta t,+}}{\nabla\tilde{\bw}}_{L^2} \dt
    + \mu \int_0^T \skp{\bbF^{\Delta t,+} (\bbF^{\Delta t,+})^\top }{\nabla\tilde{\bw}}_{L^2} \dt ,
    \\
    \label{eq:F_time_continuous}
    \nonumber
    0 &= 
    \int_0^T \skp{\partial_t \bbF^{\Delta t}}{\tilde{\bbG}}_{L^2} \dt
    - \int_0^T 
    \skp{\bbF^{\Delta t,+} \otimes \bv^{\Delta t,-}}{\nabla\tilde{\bbG}}_{L^2} \dt
    \\ \nonumber
    &\quad 
    + \frac{\mu}{2\lambda} \int_0^T \skp{ \bbF^{\Delta t,+} (\bbF^{\Delta t,+})^\top \bbF^{\Delta t,+} - \bbF^{\Delta t,+}}{\tilde{\bbG}}_{L^2}  \dt
    - \int_0^T \skp{ (\nabla\bv^{\Delta t,+}) \bbF^{\Delta t,+} }{\tilde{\bbG}}_{L^2} \dt
    \\*
    &\quad
    + \Delta t \int_0^T \skp{\nabla\bbF^{\Delta t,+}}{\nabla\tilde{\bbG}}_{L^2} \dt,
\end{align}
\end{subequations}
subject to the initial data $\bv^{\Delta t}(0) = \bv^0$ and $\bbF^{\Delta t}(0) = \bbF^0$.

Under the constraint $\Delta t < \frac{\mu}{\lambda}$, we directly obtain the following estimate from \eqref{eq:stability_time_discrete_gronwall}:
\begin{align}
    \label{eq:stability_time_continuous} \nonumber
    & \norm{\bv^{\Delta t,\pm}}_{L^\infty(0,T;L^2)}^2 
    + \norm{\bbF^{\Delta t,\pm}}_{L^\infty(0,T;L^2)}^2 
    \\ \nonumber
    &\quad
    + \tfrac{1}{\Delta t} \norm{\bv^{\Delta t,+} - \bv^{\Delta t,-}}_{L^2(0,T;L^2)}^2
    + \tfrac{1}{\Delta t} \norm{\bbF^{\Delta t,+}- \bbF^{\Delta t,-}}_{L^2(0,T;L^2)}^2
    \\ \nonumber
    &\quad
    + \norm{\nabla\bv^{\Delta t,+}}_{L^2(0,T;L^2)}^2
    + \norm{\bbF^{\Delta t,+} (\bbF^{\Delta t,+})^\top}_{L^2(0,T;L^2)}^2
    + \Delta t \norm{\nabla\bbF^{\Delta t+}}_{L^2(0,T;L^2)}^2
    \\
    &\leq 
    \exp(CT) \left( \norm{\bv^0}_{L^2}^2 
    + \norm{\bbF^0}_{L^2}^2 \right).
\end{align}
Moreover, the estimates in \eqref{eq:time_derivatives_discrete} imply
\begin{align}
    \label{eq:time_derivatives_continuous} \nonumber
    &\nnorm{\partial_t \bv^{\Delta t}}_{L^{4/d}(0,T;(H^1_{0,\Div})')}^{4/d} 
    + \nnorm{\partial_t \bbF^{\Delta t}}_{L^{4/3}(0,T;(H^1)')}^{4/3}  
    \\
    &\leq C(T) \left( 1 + \norm{\bv^0}_{L^2}^2 
    + \norm{\bbF^0}_{L^2}^2 
    + \Delta t \norm{\nabla\bv^0}_{L^2}^2 \right) .
\end{align}
Furthermore, it follows from the Cauchy--Schwarz inequality that 
\begin{align}
    \label{eq:equivalence_L4}
    \frac1d |\bbA|^4 = \frac1d \trace(\bbA\bbA^\top)^2 \leq \trace( (\bbA\bbA^\top)^2 ) = |\bbA\bbA^\top|^2 
    \leq |\bbA|^4 
    \qquad \forall\, \bbA\in\bbR^{d\times d}.
\end{align}
Combining \eqref{eq:equivalence_L4} with \eqref{eq:stability_time_continuous} yields
\begin{align} 
    \label{eq:F_in_L4L4}
    \frac1d \norm{\bbF^{\Delta t,+}}_{L^4(0,T;L^4)}^4
    &\leq 
    \norm{\bbF^{\Delta t,+} (\bbF^{\Delta t,+})^\top}_{L^2(0,T;L^2)}^2
    \leq 
    \exp(CT) \left( \norm{\bv^0}_{L^2}^2 
    + \norm{\bbF^0}_{L^2}^2 \right).
\end{align}

If additionally $\bbF^0 \in L^4(\Omega;\bbR^{d\times d})$, then using \eqref{eq:equivalence_L4} we have
\begin{align}   
    \label{eq:Fminus_in_L4L4}
    &\nonumber
    \norm{\nabla\bv^{\Delta t,-}}_{L^2(0,T;L^2)}^2
    + \frac1d \norm{\bbF^{\Delta t,-}}_{L^4(0,T;L^4)}^4
    \\ \nonumber
    &\leq
    \norm{\nabla\bv^{\Delta t,+}}_{L^2(0,T;L^2)}^2
    + \norm{\bbF^{\Delta t,+} (\bbF^{\Delta t,+})^\top}_{L^2(0,T;L^2)}^2
    + \Delta t \norm{\nabla\bv^0}_{L^2}^2 
    + \Delta t \norm{\bbF^0 (\bbF^0)^\top}_{L^2}^2
    \\
    &\leq d \exp(CT) \left( \norm{\bv^0}_{L^2}^2 
    + \norm{\bbF^0}_{L^2}^2 \right)
    + \Delta t \norm{\nabla\bv^0}_{L^2}^2 
    + \Delta t \norm{\bbF^0 (\bbF^0)^\top}_{L^2}^2.
\end{align}
Compared to \eqref{eq:stability_time_discrete_gronwall}, we have the additional terms $\Delta t \norm{\nabla\bv^0}_{L^2}^2$ and $\Delta t \norm{\bbF^0 (\bbF^0)^\top}_{L^2}^2$ for the initial data $\bv^0$, $\bbF^0$ from \eqref{eq:system_time_discrete} that we need to control.
%

%
%
%
\subsection{Convergence to a weak solution}

In what follows, we specify the requirements for the initial data of the time-discrete scheme \eqref{eq:system_time_discrete}. Let the discrete initial velocity and deformation gradient be denoted by, respectively, $\bv^0 \in H^1_{0,\Div}(\Omega;\bbR^d)$ and $\bbF^0 \in L^4(\Omega;\bbR^{d\times d})$. We assume that these functions approximate the continuous initial data $\bv_0 \in L^2_{\Div}(\Omega;\bbR^d)$ and $\bbF_0\in L^2(\Omega;\bbR^{d\times d})$ from \eqref{eq:system_init} so that they satisfy the following stability estimates, uniformly with respect to $\Delta t>0$:
\begin{align}\label{eq:v0_init_stability}
    \norm{\bv^0}_{L^2}^2 + \Delta t \norm{\nabla\bv^0}_{L^2}^2 
    &\leq C \norm{\bv_0}_{L^2}^2,
    \\
    \label{eq:F0_init_stability}
    \norm{\bbF^0}_{L^2}^2 + \Delta t \norm{\bbF^0 (\bbF^0)^\top}_{L^2}^2
    &\leq C \norm{\bbF_0}_{L^2}^2.
\end{align}
Furthermore, we assume strong convergence to the original data as the time step vanishes, i.e., for $\Delta t \to 0$,
\begin{alignat}{3}\label{eq:v0_init_convergence}
    \bv^0 &\to \bv_0\quad   
    &&\text{strongly in } L^2_{\Div}(\Omega;\bbR^d),
    \\
    \label{eq:F0_init_convergence}
    \bbF^0 &\to \bbF_0\quad
    &&\text{strongly in } L^2(\Omega;\bbR^{d\times d}).
\end{alignat}

In the next lemma, we show the existence of such approximations by providing a concrete construction.


\begin{lemma}\label{lemma:init_construction}
There exist functions $\bv^0 \in H^1_{0,\Div}(\Omega;\bbR^d)$ and $\bbF^0 \in L^4(\Omega;\bbR^{d\times d})$ satisfying the stability and convergence properties \eqref{eq:v0_init_stability}--\eqref{eq:F0_init_convergence}. Specifically, we define them as the unique solutions to the following variational problems:
\begin{align}
    \label{eq:v0_init_def}
    \skp{\bv^0}{\bw}_{L^2} + \Delta t \skp{\nabla\bv^0}{\nabla\bw}_{L^2}
    &= \skp{\bv_0}{\bw}_{L^2}\quad \forall\, \bw\in H^1_{0,\Div}(\Omega;\bbR^d),
    \\
    \label{eq:F0_init_def}
    \skp{\bbF^0}{\bbG}_{L^2} + \Delta t \skp{\bbF^0(\bbF^0)^\top \bbF^0}{\bbG}_{L^2}
    &= \skp{\bbF_0}{\bbG}_{L^2}\quad \forall\, \bbG\in L^4(\Omega;\bbR^{d\times d}).
\end{align}
\end{lemma}

\begin{proof}
Given $\bv_0\in L^2_{\Div}(\Omega;\bbR^d)$, the existence of a unique $\bv^0\in H^1_{0,\Div}(\Omega;\bbR^d)$ solving \eqref{eq:v0_init_def} follows from the Lax--Milgram theorem, for all $\Delta t>0$. By taking $\bw=\bv^0$ and using Hölder's and Young's inequalities, we directly obtain \eqref{eq:v0_init_stability} with $C=1$.
Next, for arbitrary $\bw \in H^1_{0,\Div}(\Omega;\bbR^d)$, we have with Hölder's inequality and \eqref{eq:v0_init_stability} with $C=1$, that
\begin{align*}
    \abs{\skp{\bv^0 - \bv_0}{\bw}_{L^2}} 
    \leq \Delta t \abs{\skp{\nabla\bv^0}{\nabla\bw}_{L^2}}
    \leq \sqrt{\Delta t} \norm{\bv_0}_{L^2} \norm{\nabla\bw}_{L^2}.
\end{align*}
By taking the supremum over all $\bw\in H^1_{0,\Div}(\Omega;\bbR^d)$ with $\norm{\nabla\bw}_{L^2}=1$, we obtain
\begin{align}
    \label{eq:v0_init_convergence_H1_0div'}
    \norm{\bv^0 - \bv_0}_{(H^1_{0,\Div})'}
    \leq \sqrt{\Delta t} \norm{\bv_0}_{L^2} 
    \to 0,
\end{align}
as $\Delta t\to 0$.
As $\norm{\bv^0}_{L^2} \leq \norm{\bv_0}_{L^2}$, which is due to \eqref{eq:v0_init_stability} with $C=1$, there exist a function $\bz \in L^2_{\Div}(\Omega;\bbR^d)$ and a non-relabeled subsequence such that $\bv^0 \rightharpoonup \bz$ weakly in $L^2_{\Div}(\Omega;\bbR^d)$. By the uniqueness of weak limits, together with \eqref{eq:v0_init_convergence_H1_0div'}, we have $\bz = \bv_0$.
Since 
\begin{align*}
    \norm{\bv^0}_{L^2}^2 + \Delta t \norm{\nabla\bv^0}_{L^2}^2
    = \skp{\bv_0}{\bv^0}_{L^2}
    = \skp{\bv_0}{\bv_0 - \bv^0}_{L^2} - \norm{\bv_0}_{L^2}^2
    + 2 \skp{\bv_0}{\bv^0}_{L^2},
\end{align*}
it follows that
\begin{align*}
    \norm{\bv^0 - \bv_0}_{L^2}^2 + \Delta t \norm{\nabla\bv^0}_{L^2}^2
    = \skp{\bv_0}{\bv_0 - \bv^0}_{L^2} .
\end{align*}
Hence, from the weak convergence of $\bv^0$ to $\bv_0$ in $L^2_{\Div}(\Omega;\bbR^d)$ we deduce that
\begin{align}
    \label{eq:v0_init_convergence_L2_proof}
    \norm{\bv^0 - \bv_0}_{L^2}^2
    \leq \abs{\skp{\bv_0}{\bv_0 - \bv^0}_{L^2}} \to 0,
\end{align}
as $\Delta t\to 0$. This proves \eqref{eq:v0_init_convergence}.

Next, we consider the functional $J\colon\, L^4(\Omega;\bbR^{d\times d})\to \bbR$ defined by
\begin{align*}
    J(\bbG) = \norm{\bbG}_{L^2}^2 + \frac{\Delta t}{4} \norm{\bbG\bbG^\top}_{L^2}^2 - \skp{\bbF_0}{\bbG}_{L^2} 
    \quad \forall\, \bbG\in L^4(\Omega;\bbR^{d\times d}) ,
\end{align*}
which is continuous, bounded, and strictly convex. Here we note that, thanks to \eqref{eq:equivalence_L4}, we have
\begin{align*}
    \frac1d \norm{\bbG}_{L^4}^4 
    \leq \norm{\bbG\bbG^\top}_{L^2}^2
    \leq \norm{\bbG}_{L^4}^4 \quad \forall\, \bbG\in L^4(\Omega;\bbR^{d\times d}).
\end{align*}
Thus, the existence of a unique minimizer $\bbF^0 \in L^4(\Omega;\bbR^{d\times d})$ follows by the direct method of the calculus of variations, for all $\Delta t>0$. Moreover, $\bbF^0$ uniquely solves \eqref{eq:F0_init_def}, which are the corresponding Euler--Lagrange equations. The inequality \eqref{eq:F0_init_stability} with $C=1$ follows directly from testing \eqref{eq:F0_init_def} with $\bbF^0$, and applying Hölder's and Young's inequalities. The strong convergence of $\bbF^0$ to $\bbF_0$ in $L^{4/3}(\Omega;\bbR^{d\times d}) \, \widetilde{=} \, (L^4(\Omega;\bbR^{d\times d}))'$, as $\Delta t\to 0$, can be shown analogously to \eqref{eq:v0_init_convergence_H1_0div'}. By \eqref{eq:F0_init_stability} and the uniqueness of weak limits, there exists a non-relabeled subsequence, such that $\bbF^0 \rightharpoonup \bbF_0$ weakly in $L^2(\Omega;\bbR^{d\times d})$. Finally, the strong convergence of $\bbF^0$ to $\bbF_0$ in $L^2(\Omega;\bbR^{d\times d})$ follows similarly to \eqref{eq:v0_init_convergence_L2_proof}. 
\end{proof}

The goal of the remainder of this section is to prove the following theorem.

\begin{theorem}\label{theorem:convergence}
Let $d\in \{2,3\}$. Assume that the time-discrete initial data $\bv^0 \in H^1_{0,\Div}(\Omega;\bbR^d)$ and $\bbF^0 \in L^4(\Omega;\bbR^{d\times d})$ satisfy \eqref{eq:v0_init_stability}--\eqref{eq:F0_init_convergence}, and suppose that $\Delta t < \frac{\lambda}{\mu}$. Then, there exists a (non-relabeled) subsequence of $\{\bv^{\Delta t(,\pm)}, \bbF^{\Delta t(,\pm)}\}_{\Delta t>0}$ solving \eqref{eq:system_time_continuous}, and there exist functions 
\begin{align*}
    &\bv \in C([0,T]; L^2_{\Div}(\Omega;\bbR^d)) \cap L^2(0,T; H^1_{0,\Div}(\Omega;\bbR^d))
    \cap W^{1,4/d}(0,T; (H^1_{0,\Div}(\Omega;\bbR^d))')\, ,
    \\
    &\bbF \in C([0,T]; L^2(\Omega;\bbR^{d\times d})) \cap L^4(\Omega_T;\bbR^{d\times d})
    \cap W^{1,4/3}(0,T; (H^1(\Omega;\bbR^{d\times d}))')\, ,
\end{align*}
with $\bv(0) = \bv_0$ and $\bbF(0)=\bbF_0$ a.e.~in $\Omega$, such that, as $\Delta t \to 0$,
\begin{alignat}{3}
    \label{eq:v_conv_LinfL2}
    \bv^{\Delta t(,\pm)}  &\rightharpoonup^* \bv \quad &&\text{weakly-$*$ } \quad 
    &&\text{in } L^\infty(0,T;L^2_{\Div}(\Omega;\bbR^d)),
    \\
    \label{eq:v_conv_L2H1}
    \bv^{\Delta t(,\pm)}  &\rightharpoonup \bv \quad &&\text{weakly } \quad 
    &&\text{in } L^2(0,T; H^1_{0,\Div}(\Omega;\bbR^d)),
    \\
    \label{eq:v_conv_H1H1'}
    \partial_t\bv^{\Delta t}  &\rightharpoonup \partial_t \bv \quad &&\text{weakly } \quad 
    &&\text{in } L^{4/d}(0,T; (H^1_{0,\Div}(\Omega;\bbR^d))'),
    \\
    \label{eq:v_conv_L2Lr}
	\bv^{\Delta t(,\pm)}  &\to \bv \quad &&\text{strongly } \quad 
	&&\text{in } L^r(\Omega_T; \bbR^{d}) \cap L^2(0,T; L^s(\Omega;\bbR^d)) ,
    \\
    \label{eq:F_conv_LinfL2}
    \bbF^{\Delta t(,\pm)}  &\rightharpoonup^* \bbF \quad &&\text{weakly-$*$ } \quad 
    &&\text{in } L^\infty(0,T;L^2(\Omega;\bbR^{d\times d})),
    \\
    \label{eq:F_conv_L4L4}
    \bbF^{\Delta t(,\pm)}  &\rightharpoonup \bbF \quad &&\text{weakly } \quad 
    &&\text{in } L^4(\Omega_T;\bbR^{d\times d}),
    \\
    \label{eq:F_conv_L2H1}
    \Delta t \, \nabla \bbF^{\Delta t,+} &\to \mathbb{O} \quad &&\text{strongly } \quad 
    &&\text{in } L^2(\Omega_T;\bbR^{d\times d\times d}),
    \\
    \label{eq:F_conv_H1H1'}
    \partial_t\bbF^{\Delta t}  &\rightharpoonup \partial_t\bbF \quad &&\text{weakly } \quad 
    &&\text{in } L^{4/3}(0,T; (H^1(\Omega;\bbR^{d\times d}))'),
\end{alignat}
for all $r\in[1, \frac{2(d+2)}{d})$ and $s\in[1,\tfrac{2d}{d-2})$. Moreover, we have that
\begin{align}
    \label{eq:stability_time_continuous_gronwall} \nonumber
    & \norm{\bv}_{L^\infty(0,T;L^2)}^2
    + \norm{\bbF}_{L^\infty(0,T;L^2)}^2
    + \norm{\nabla\bv}_{L^2(0,T;L^2)}^2
    + \norm{\bbF}_{L^4(0,T;L^4)}^4
    \\
    &\leq 
    \exp(CT) \left( \norm{\bv_0}_{L^2}^2 
    + \norm{\bbF_0}_{L^2}^2 \right),
\end{align}
and
\begin{align}
    \label{eq:time_derivatives_weak} 
    &\norm{\partial_t \bv}_{L^{4/d}(0,T;(H^1_{0,\Div})')}^{4/d} 
    + \norm{\partial_t \bbF}_{L^{4/3}(0,T;(H^1)')}^{4/3}  
    \leq C(T) \left( 1 + \norm{\bv_0}_{L^2}^2 
    + \norm{\bbF_0}_{L^2}^2  \right) .
\end{align}

Further, by restricting to dimension $d=2$, we have
\begin{alignat}{3}
    \label{eq:F_conv_strong}
    \bbF^{\Delta t(,\pm)}  &\to \bbF \quad &&\text{strongly } \quad 
    &&\text{in } L^2(\Omega_T;\bbR^{2\times 2}),
\end{alignat}
and the limit functions $\bv$ and $\bbF$ solve \eqref{eq:v_weak} and \eqref{eq:F_weak}.
Moreover, if the initial data additionally satisfy $\det(\bbF_0)>0$ a.e.~in $\Omega$ and $\ln\det(\bbF_0)\in L^1(\Omega)$, then the limit function $\bbF$ satisfies $\det(\bbF)>0$ a.e.~in $\Omega_T$, and the following estimate holds:
\begin{align}
    \label{eq:F_lndet}
    \norm{\ln\det(\bbF)}_{L^\infty(0,T;L^1)}
    \leq C \norm{\ln\det(\bbF_0)}_{L^1}.
\end{align}
\end{theorem}

We divide the proof of Theorem~\ref{theorem:convergence} into four steps, presented in the following subsections.
First, using a priori estimates, we apply weak-($*$) compactness results to extract limit functions as $\Delta t \to 0$ (Section~\ref{sec:limit_vanishing_dt}). This provides candidate solutions but does not immediately identify all nonlinear terms because of the lack of compactness for the deformation gradient $\bbF$.
Next, we establish the strong continuity of $\bbF$ in Section~\ref{sec:continuity_F}.
To address the compactness of $\bbF$ in two dimensions, we adapt the approach of \cite{bulicek_2022_giesekus_2d} to the time-discrete setting (Section~\ref{sec:compactness_F}), which is the most involved part of the proof.
Finally, under suitable assumptions on the initial data, we establish positivity of $\det(\bbF)$ at the level of weak solutions (Section~\ref{sec:positivity_F}).

\subsection{Limit passage in time} \label{sec:limit_vanishing_dt}

Let $d\in\{2,3\}$. The goal of this subsection is to send the time discretization parameter $\Delta t$ to zero and identify convergent subsequences of time-discrete solutions to \eqref{eq:system_time_discrete}.

The uniform bounds \eqref{eq:stability_time_continuous} and \eqref{eq:time_derivatives_continuous} ensure the existence of functions $\bv$ and $\bbF$ with
\begin{align*}
    &\bv \in L^\infty(0,T; L^2_{\Div}(\Omega;\bbR^d)) \cap L^2(0,T; H^1_{0,\Div}(\Omega;\bbR^d))
    \cap W^{1,4/d}(0,T; (H^1_{0,\Div}(\Omega;\bbR^d))')\, ,
    \\
    &\bbF \in L^\infty(0,T; L^2(\Omega;\bbR^{d\times d})) \cap L^4(\Omega_T;\bbR^{d\times d})
    \cap W^{1,4/3}(0,T; (H^1(\Omega;\bbR^{d\times d}))')\, ,
\end{align*}
such that, up to a non-relabeled subsequence, the convergence results \eqref{eq:v_conv_LinfL2}, \eqref{eq:v_conv_L2H1}, \eqref{eq:v_conv_H1H1'}, \eqref{eq:F_conv_LinfL2}, \eqref{eq:F_conv_L4L4}, \eqref{eq:F_conv_L2H1}, and \eqref{eq:F_conv_H1H1'} hold, as $\Delta t\to 0$.
Moreover, since 
\[\bv\in L^2(0,T;H^1_{0,\Div}(\Omega;\bbR^d))\quad \mbox{with}\quad \partial_t \bv \in L^{4/d}(0,T; (H^1_{0,\Div}(\Omega;\bbR^d))'),\] 
while 
\[\bbF \in L^\infty(0,T;L^2(\Omega;\bbR^{d\times d}))\quad \mbox{with}\quad\partial_t \bbF \in L^{4/3}(0,T; (H^1(\Omega;\bbR^{d\times d}))'),\] 
we have that
\begin{align}
    \label{eq:v_C(L2)}
    \bv &\in C([0,T]; L^2_{\Div}(\Omega;\bbR^d)),
    \\
    \label{eq:F_Cw(L2)}
    \bbF &\in C_{w}([0,T]; L^2(\Omega;\bbR^{d\times d})).
\end{align}
Using the Aubin--Lions theorem, we first obtain \eqref{eq:v_conv_L2Lr} for $\bv^{\Delta t}$.
The same statement holds for $\bv^{\Delta t,\pm}$ because of the Gagliardo--Nirenberg inequality and the estimate
\begin{align*}
    \norm{\bv^{\Delta t} - \bv^{\Delta t,\pm}}_{L^2(0,T;L^2)}^2 
    \leq C \Delta t \left( \norm{\bv_0}_{L^2}^2 + \norm{\bbF_0}_{L^2}^2 \right),
\end{align*}
which directly follows from \eqref{eq:stability_time_continuous}.

Combining the above results about weak(-$*$) and strong convergence yields
\begin{alignat}{3}
    \label{eq:v_otimes_v_conv}
    \bv^{\Delta t,+} \otimes \bv^{\Delta t,-}  &\rightharpoonup \bv \otimes \bv \quad &&\text{weakly } \quad 
    &&\text{in } L^{4/d}(0,T;L^2(\Omega;\bbR^{d\times d})),
    \\
    \label{eq:F_otimes_v_conv}
    \bbF^{\Delta t,+} \otimes \bv^{\Delta t,-}  &\rightharpoonup \bbF \otimes \bv \quad &&\text{weakly } \quad 
    &&\text{in } L^{4/3}(0,T;L^2(\Omega;\bbR^{d\times d\times d})).
\end{alignat}

In the following we use the notation $f(u^{\Delta t}) \rightharpoonup^{(*)} \overline{f(u)}$ weakly (or weakly-$*$, respectively), in some Banach space $X$ if $f(u^{\Delta t})$ is uniformly bounded in the norm of $X$ and $u^{\Delta t}\rightharpoonup u$ in some Banach space $Y$.
It follows from the above estimates that
\begin{alignat}{3}
    \label{eq:gradv_F_conv}
    (\nabla \bv^{\Delta t,+}) \bbF^{\Delta t,+}  &\rightharpoonup \overline{(\nabla\bv) \bbF} \quad &&\text{weakly } \quad 
    &&\text{in } L^{4/3}(0,T;L^2(\Omega;\bbR^{d\times d})),
    \\
    \label{eq:FF_conv}
    \bbF^{\Delta t,+} (\bbF^{\Delta t,+})^\top &\rightharpoonup 
    \overline{\bbF \bbF^\top} \quad &&\text{weakly } \quad 
    &&\text{in } L^2(\Omega_T;\bbR^{d\times d}),
    \\
    \label{eq:|F|^2_conv}
    \abs{\bbF^{\Delta t,+}}^2 &\rightharpoonup 
    \overline{\abs{\bbF}^2} \quad &&\text{weakly } \quad 
    &&\text{in } L^2(\Omega_T),
    \\
    \label{eq:FFtF_conv}
    \bbF^{\Delta t,+} (\bbF^{\Delta t,+})^\top \bbF^{\Delta t,+}
    &\rightharpoonup 
    \overline{\bbF \bbF^\top \bbF} \quad &&\text{weakly } \quad 
    &&\text{in } L^{4/3}(\Omega_T;\bbR^{d\times d}),
    \\
    \label{eq:|FFt|^2_conv}
    \abs{\bbF^{\Delta t,+} (\bbF^{\Delta t,+})^\top}^2 &\rightharpoonup^* 
    \overline{\abs{\bbF \bbF^\top}^2} \quad &&\text{weakly-$*$ } \quad 
    &&\text{in } \mathcal{M}(\overline{\Omega_T}),
    \\
    \label{eq:gradvFFt_conv}
    (\nabla \bv^{\Delta t,+}) \bbF^{\Delta t,+} (\bbF^{\Delta t,+})^\top
    &\rightharpoonup^*
    \overline{(\nabla\bv) \bbF \bbF^\top} \quad &&\text{weakly-$*$ } \quad 
    &&\text{in } \mathcal{M}(\overline{\Omega_T};\bbR^{d\times d}),
    \\
    \label{eq:|gradv|^2_conv}
    \abs{\nabla\bv^{\Delta t,+}}^2 &\rightharpoonup^*
    \overline{\abs{\nabla\bv}^2} \quad &&\text{weakly-$*$ } \quad 
    &&\text{in } \mathcal{M}(\overline{\Omega_T}).
\end{alignat}

Using these convergence results, we deduce that, for almost all $t\in(0,T)$ and all $\bw\in H^1_{0,\Div}(\Omega;\bbR^d)$ and $\bbG\in H^1(\Omega;\bbR^{d\times d})$, 
\begin{subequations}
\label{eq:system_weaklimit}
\begin{align}
    \label{eq:v_weaklimit}
    0 &= 
    \rho \dualp{\partial_t \bv}{\bw}_{(H^1_{0,\Div})',H^1_{0,\Div}} 
    - \rho \skp{\bv \otimes \bv}{\nabla\bw}_{L^2} 
    + \nu \skp{\nabla \bv}{\nabla\bw}_{L^2} 
    + \mu \skp{\overline{\bbF \bbF^\top} }{\nabla\bw}_{L^2},
    \\
    \label{eq:F_weaklimit}
    0 &= 
    \dualp{\partial_t \bbF}{\bbG}_{(H^1)',H^1} 
    - \skp{\bbF \otimes \bv}{\nabla\bbG}_{L^2}
    + \frac{\mu}{2\lambda} \skp{ \overline{\bbF \bbF^\top \bbF} - \bbF}{\bbG}_{L^2} 
    - \skp{ \overline{(\nabla\bv) \bbF} }{\bbG}_{L^2}.
\end{align}
\end{subequations}
Attainment of the initial data $\bv_0 \in L^2_{\Div}(\Omega;\bbR^d)$ and $\bbF_0 \in L^2(\Omega;\bbR^{d\times d})$, respectively, follows from (weak) continuity in time; that is, \eqref{eq:v_C(L2)} and \eqref{eq:F_Cw(L2)} hold. 

The estimates \eqref{eq:stability_time_continuous_gronwall} and \eqref{eq:time_derivatives_weak} follow from \eqref{eq:stability_time_continuous}, \eqref{eq:time_derivatives_continuous}, \eqref{eq:F_in_L4L4}, \eqref{eq:v0_init_stability}, \eqref{eq:F0_init_stability}, \eqref{eq:v_conv_LinfL2}--\eqref{eq:F_conv_H1H1'}, and the weak(-*) lower semicontinuity of the $L^p$ norms. 

\subsection{Strong continuity of the deformation gradient in time} \label{sec:continuity_F}

Let $d\in\{2,3\}$. Our next aim is to show strong continuity of the deformation gradient in time, i.e.,
\begin{align}
    \label{eq:F_C(L2)}
    \bbF \in C([0,T]; L^2(\Omega;\bbR^{d\times d})).
\end{align}

For arbitrary $t_0, t_1 \in (0,T)$, one can derive the identity
\begin{align}
    \label{eq:F_strong_continuity_proof}
    \norm{\bbF(t_1)}_{L^2}^2 - \norm{\bbF(t_0)}_{L^2}^2
    &=
    2 \int_{t_0}^{t_1} \skp{\overline{(\nabla\bv) \bbF}}{\bbF}_{L^2}
    - \int_{t_0}^{t_1} \skp{\overline{\bbF \bbF^\top \bbF}}{\bbF}_{L^2}
    + \int_{t_0}^{t_1} \norm{\bbF}_{L^2}^2.
\end{align}
A rigorous proof of \eqref{eq:F_strong_continuity_proof} is provided in \cite[Sect.~6.3]{bulicek_2022_giesekus_2d}. The main idea is to establish an analogous relation at an approximate level by mollifying the weak limit system \eqref{eq:F_weaklimit} in the spatial variable, and then letting the mollification parameter tend to zero. The reason behind this workaround is that $\bbF$ itself is not an admissible choice as a test function in \eqref{eq:F_weaklimit} because of its limited spatial regularity.

The integrals on the right-hand side of \eqref{eq:F_strong_continuity_proof} are well-defined thanks to \eqref{eq:gradv_F_conv}, \eqref{eq:FFtF_conv} and \eqref{eq:F_conv_L4L4}. This allows the following limit passages in \eqref{eq:F_strong_continuity_proof}:
\begin{gather*}
    \lim_{t_1 \to t_0} \norm{\bbF(t_1)}_{L^2}^2 = \norm{\bbF(t_0)}_{L^2}^2 
    \quad \forall\, t_0 \in(0,T),
    \\
    \lim_{t_1 \to 0^+} \norm{\bbF(t_1)}_{L^2}^2 = \norm{\bbF(0)}_{L^2}^2,
    \qquad
    \lim_{t_1 \to T^-} \norm{\bbF(t_1)}_{L^2}^2 = \norm{\bbF(T)}_{L^2}^2,
\end{gather*}
which, together with \eqref{eq:F_Cw(L2)}, imply \eqref{eq:F_C(L2)}.

\subsection{Compactness of the deformation gradient in two dimensions} \label{sec:compactness_F}
From now on we restrict our analysis to two space dimensions, i.e., we now fix $d=2$. For the sake of simplicity of the exposition, we also set the physical parameters $\rho$, $\nu$, $\mu$, and $\lambda$ to~1, as their specific values have no impact on any of our assertions.

The main goal of this subsection is to prove \eqref{eq:F_conv_strong}, i.e., the strong convergence of $\bbF^{\Delta t(,\pm)}$ to $\bbF$ in $L^2(\Omega_T;\bbR^{2\times 2})$. In this context, we note that
\begin{align}
    \label{eq:F_norm_compactness}
    \lim_{\Delta t \to 0} \norm{\bbF^{\Delta t(,\pm)} - \bbF}_{L^2(\Omega_T)}^2
    = \lim_{\Delta t \to 0} \int_{\Omega_T} \left( \abs{\bbF^{\Delta t(,\pm)}}^2 - 2 \bbF^{\Delta t(,\pm)} : \bbF + \abs{\bbF}^2 \right)
    = \int_{\Omega_T} \left( \overline{\abs{\bbF}^2} - \abs{\bbF}^2 \right),
\end{align}
where the last equality follows from \eqref{eq:F_conv_LinfL2} and \eqref{eq:|F|^2_conv}. Hence, thanks to \eqref{eq:F_norm_compactness}, it is sufficient to show that $\overline{\abs{\bbF}^2} = \abs{\bbF}^2$ almost everywhere in $\Omega_T$. 
As before, we follow the approach of \cite[Sect.~6.4]{bulicek_2022_giesekus_2d} and adapt the strategy therein to our time-discrete setting.

Our aim is to derive a differential inequality for $\overline{\abs{\bbF}^2} - \abs{\bbF}^2$. Once this inequality is established, a renormalization argument will imply that $\overline{\abs{\bbF}^2} = \abs{\bbF}^2$ almost everywhere in $\Omega_T$.
To obtain such an inequality, we begin by deriving suitable (in-)equalities for $\overline{\abs{\bbF}^2}$ and $\abs{\bbF}^2$, and then subtract them.
In this procedure, mixed terms appear that are measure-valued, e.g., the limit functions in \eqref{eq:|FFt|^2_conv}, \eqref{eq:gradvFFt_conv} and \eqref{eq:|gradv|^2_conv}. To estimate or control these terms, we need to perform additional testing in \eqref{eq:v_weaklimit} and \eqref{eq:v_time_continuous} with functions that are not divergence-free and hence not admissible. For this reason, we will reconstruct local pressures before continuing the argument.

\subsubsection{Differential (in-)equalities}
In the first step, we aim to derive differential (in-)equalities for $\overline{\abs{\bbF}^2}$ and $\abs{\bbF}^2$.

For future reference, we note the formula for discrete integration by parts with respect to the time variable, i.e., for any real-valued sequences $(a_n)_{n\in\bbN_0}$, $(b_n)_{n\in\bbN_0}$, one has, for any $m\in\bbN$,  
\begin{align}
    \label{eq:discrete_integration_by_parts_in_time}
    \sum_{n=1}^{m} (a_n - a_{n-1}) b_n 
    + \sum_{n=1}^{m} a_{n-1} (b_n - b_{n-1}) 
    = a_m b_m - a_0 b_0 .
\end{align}

We set $\tilde{\bbG} = \bbF^n \phi^n \in H^1(\Omega;\bbR^{2\times 2})$ in \eqref{eq:F_time_discrete} with $\phi^n = \frac{1}{\Delta t} \int_{t^{n-1}}^{t^n} \phi(\cdot,t)\dt$, and where the function $\phi \in C^\infty_c(\Omega \times (-\infty,T))$ with $\phi\geq0$. This gives us
\begin{align*}
    0 &= \frac{1}{2\Delta t} \skp{\abs{\bbF^n}^2 - \abs{\bbF^{n-1}}^2 + \abs{\bbF^n - \bbF^{n-1}}^2 }{\phi^n}_{L^2}
    - \frac12 \skp{\bv^{n-1} \cdot \nabla\phi^n}{\abs{\bbF^n}^2}_{L^2}
    \\
    &\quad
    + \frac12 \skp{ \abs{\bbF^n (\bbF^n)^\top}^2 -  \abs{\bbF^n}^2} {\phi^n}_{L^2}
    - \skp{ \nabla\bv^n : (\bbF^n (\bbF^n)^\top) }{\phi^n}_{L^2}
    \\
    &\quad
    + \Delta t \skp{\abs{\nabla\bbF^n}^2}{\phi^n}_{L^2}
    - \frac12 \Delta t \skp{\abs{\bbF^n}^2}{\Delta \phi^n}_{L^2},
\end{align*}
where we have used the chain rule and integration by parts over $\Omega$ to compute
\begin{align*}
    \skp{(\bv^{n-1} \cdot \nabla) \bbF^n}{\bbF^n \phi^n}_{L^2}
    = \frac12 \skp{\bv^{n-1} \cdot \nabla \abs{\bbF^n}^2}{\phi^n}_{L^2}
    = -\frac12 \skp{\bv^{n-1} \cdot \nabla\phi^n}{\abs{\bbF^n}^2}_{L^2} ,
\end{align*}
and
\begin{align*}
    \skp{\nabla \bbF^n}{\nabla(\bbF^n \phi^n)}_{L^2}
    = \skp{\abs{\nabla\bbF^n}^2}{\phi^n}_{L^2}
    - \frac12 \skp{\abs{\bbF^n}^2}{\Delta \phi^n}_{L^2} .
\end{align*}
After multiplying both sides by $2\Delta t$ and summing over $n\in\{1,\ldots,N_T\}$, we have 
\begin{align*}
    0 &= \sum_{n=1}^{N_T} \skp{\abs{\bbF^n}^2 - \abs{\bbF^{n-1}}^2 + \abs{\bbF^n - \bbF^{n-1}}^2 }{\phi^n}_{L^2}
    - \Delta t \sum_{n=1}^{N_T}  \skp{\bv^{n-1} \cdot \nabla\phi^n}{\abs{\bbF^n}^2}_{L^2}
    \\
    &\quad
    + \Delta t \sum_{n=1}^{N_T}  \skp{ \abs{\bbF^n (\bbF^n)^\top}^2 -  \abs{\bbF^n}^2 }{\phi^n}_{L^2}
    - 2 \Delta t \sum_{n=1}^{N_T}  \skp{ \nabla\bv^n : (\bbF^n (\bbF^n)^\top) }{\phi^n}_{L^2}
    \\
    &\quad
    + 2 (\Delta t)^2 \sum_{n=1}^{N_T}   \skp{\abs{\nabla\bbF^n}^2}{\phi^n}_{L^2}
    - (\Delta t)^2 \sum_{n=1}^{N_T} \skp{\abs{\bbF^n}^2}{\Delta \phi^n}_{L^2}.
\end{align*}
Now, by applying \eqref{eq:discrete_integration_by_parts_in_time} and using the notation introduced in \eqref{eq:fun_Delta_t} and \eqref{eq:fun_Delta_t_pm}, we obtain
\begin{align*}
    0 &= - \int_0^T \skp{\abs{\bbF^{\Delta t,-}}^2}{\partial_t \phi^{\Delta t}}_{L^2}
    + \int_0^T \skp{\frac{1}{\Delta t} \abs{\bbF^{\Delta t,+} - \bbF^{\Delta t,-}}^2}{\phi}_{L^2}
    \\
    &\quad
    + \skp{\abs{\bbF^{\Delta t}(T)}^2}{\phi^{\Delta t}(T)}_{L^2}
    - \skp{\abs{\bbF^{\Delta t}(0)}^2}{\phi^{\Delta t}(0)}_{L^2}
    - \int_0^T  \skp{\bv^{\Delta t,-} \cdot \nabla \phi}{\abs{\bbF^{\Delta t,+}}^2}_{L^2}
    \\
    &\quad
    + \int_0^T \skp{ \abs{\bbF^{\Delta t,+} (\bbF^{\Delta t,+})^\top}^2 -  \abs{\bbF^{\Delta t,+}}^2 }{\phi}_{L^2}
    - 2\int_0^T \skp{ \nabla\bv^{\Delta t,+} : (\bbF^{\Delta t,+} (\bbF^{\Delta t,+})^\top) }{\phi }_{L^2}
    \\
    &\quad
    + 2 \Delta t \int_0^T  \skp{\abs{\nabla\bbF^{\Delta t,+}}^2}{\phi}_{L^2}
    - \Delta t \int_0^T \skp{\abs{\bbF^{\Delta t,+}}^2}{\Delta \phi}_{L^2},
\end{align*}
where  $\phi^{\Delta t}$ is the continuous, piecewise affine interpolant in time of $\phi^n$ based on the notation \eqref{eq:fun_Delta_t}.
Using the assumed non-negativity of $\phi$, we then deduce the inequality
\begin{align*}
    &- \int_0^T \skp{\abs{\bbF^{\Delta t,-}}^2}{\partial_t \phi^{\Delta t}}_{L^2}
    - \skp{\abs{\bbF^{\Delta t}(0)}^2}{\phi^{\Delta t}(0)}_{L^2}
    - \int_0^T  \skp{\bv^{\Delta t,-} \cdot \nabla \phi}{\abs{\bbF^{\Delta t,+}}^2}_{L^2}
    \\
    &\quad
    + \int_0^T \skp{ \abs{\bbF^{\Delta t,+} (\bbF^{\Delta t,+})^\top}^2 -  \abs{\bbF^{\Delta t,+}}^2 }{\phi}_{L^2}
    - 2\int_0^T \skp{ \nabla\bv^{\Delta t,+} : (\bbF^{\Delta t,+} (\bbF^{\Delta t,+})^\top) }{\phi }_{L^2}
    \\
    &\leq 
    \Delta t \int_0^T \skp{\abs{\bbF^{\Delta t,+}}^2}{\Delta \phi}_{L^2}.
\end{align*}
By sending $\Delta t\to 0$ and noting \eqref{eq:|F|^2_conv}, \eqref{eq:|FFt|^2_conv}, \eqref{eq:gradvFFt_conv}, \eqref{eq:v_conv_L2Lr}, the strong convergence of $\bbF^{\Delta t}(0)$ to $\bbF_0$ in $L^2(\Omega;\bbR^{2\times 2})$, and the fact that \eqref{eq:|F|^2_conv} is also valid for $|\bbF^{\Delta t,-}|^2$ instead of $|\bbF^{\Delta t,+}|^2$ (which directly follows from \eqref{eq:stability_time_continuous}), 
%
%
%
%
we obtain the following inequality 
\begin{align}
    \label{eq:differential_inequality_overline|F|^2}
    \nonumber
    &- \int_{\Omega_T} \overline{\abs{\bbF}^2} \partial_t \phi 
    - \int_\Omega \abs{\bbF_0}^2 \phi(0)
    - \int_{\Omega_T} \left( \overline{\abs{\bbF}^2} \bv \right) \cdot \nabla\phi
    - 2 \dualp{ \overline{\nabla \bv : (\bbF\bbF^\top)} }{\phi}_{\mathcal{M}(\overline{\Omega_T}), C(\overline{\Omega_T})}
    \\*
    &\quad 
    + \dualp{ \overline{\abs{\bbF\bbF^\top}^2} }{\phi}_{\mathcal{M}(\overline{\Omega_T}), C(\overline{\Omega_T})}
    - \int_{\Omega_T} \overline{\abs{\bbF}^2} \phi 
    \leq 0,
\end{align}
which holds for all non-negative test functions $\phi\in C^\infty_c(\Omega\times(-\infty,T))$.

Using a similar (though technically more involved) approach, one can also obtain the following differential equality for $\abs{\bbF}^2$:
\begin{align}
    \label{eq:differential_inequality_|F|^2}
    \nonumber
    &- \int_{\Omega_T} \abs{\bbF}^2 \partial_t \phi 
    - \int_\Omega \abs{\bbF_0}^2 \phi(0)
    - \int_{\Omega_T} \left( \abs{\bbF}^2 \bv \right) \cdot \nabla\phi
    - 2 \int_{\Omega_T} \overline{\nabla \bv \bbF} : (\phi\bbF)
    \\*
    &\quad 
    + \int_{\Omega_T} \left( \overline{\bbF\bbF^\top\bbF}:\bbF - \abs{\bbF}^2 \right) \phi
    = 0,
\end{align}
which holds for all test functions $\phi\in C^\infty_c(\Omega\times(-\infty,T))$.
The proof of this can be found in \cite[Sect.~6.4]{bulicek_2022_giesekus_2d}. The main idea is to establish an analogue of \eqref{eq:differential_inequality_|F|^2} at an approximate level by mollifying the weak formulation \eqref{eq:F_weaklimit} in the spatial variable, and then passing to the limit as the mollification parameter goes to zero.
This detour is necessary since $\bbF \phi$ with $\phi \in C^\infty_c(\Omega \times (-\infty,T))$ cannot be used directly as a test function in \eqref{eq:F_weaklimit}, due to the limited spatial regularity of $\bbF$.

\subsubsection{Reconstruction of local pressures and their convergence in two dimensions}

As noted earlier, several terms involving yet-to-be-identified weak limits appear in \eqref{eq:differential_inequality_overline|F|^2} and \eqref{eq:differential_inequality_|F|^2}. To handle them, we will employ an additional testing procedure in \eqref{eq:v_weaklimit} and \eqref{eq:v_time_continuous}, which requires test functions that are not divergence-free and thus not admissible. To make this possible, we first reconstruct local pressures for \eqref{eq:v_weaklimit} and \eqref{eq:v_time_continuous}, which is summarized in the next lemma. Having done so, we can proceed to derive an inequality for $\overline{\abs{\bbF}^2} - \abs{\bbF}^2$ by subtracting \eqref{eq:differential_inequality_|F|^2} from \eqref{eq:differential_inequality_overline|F|^2}.

Our strategy is to adapt the argument in \cite[Lemma~6.1]{bulicek_2022_giesekus_2d} to the discrete setting \eqref{eq:system_time_continuous}.
For this reason, we also restrict ourselves to the two-dimensional case, since the required time-regularity, provided in the following lemma, is not available in three dimensions but will be essential for the subsequent computations.

\begin{lemma} \label{lemma:local_pressures}
Let $d=2$. The following assertions hold.
\begin{enumerate}
\item Let $\widetilde{\Omega} \subset \overline{\widetilde{\Omega}} \subset\Omega$ with $\partial\widetilde{\Omega}$ smooth. Then, for any $\Delta t< \frac{\lambda}{\mu}$, there exists a pressure $p^{\Delta t,+} = p_1^{\Delta t,+} + p_2^{\Delta t,+}$ with uniform bounds
\begin{align}
    \label{eq:p1_num_L2H2}
    \left\{ p_1^{\Delta t(,\pm)} \right\}_{\Delta t>0} &\subset L^2(0,T; H^2(\widetilde{\Omega})),
    \\
    \label{eq:p2_num_L2L2}
    \left\{ p_2^{\Delta t,+} \right\}_{\Delta t>0} &\subset L^2(0,T; L^2(\widetilde{\Omega})),
    \\
    \label{eq:dtp1_num_L2H10'}
    \left\{ \partial_t (\bv^{\Delta t} + \nabla p_1^{\Delta t})
    \right\}_{\Delta t>0} &\subset L^2(0,T; (H^1_0(\widetilde{\Omega}; \mathbb{R}^2))'),
\end{align}
and, for all $\bw\in H^1_0(\widetilde{\Omega};\mathbb{R}^2)$ and almost every $t\in(0,T)$, one has that
\begin{align}
    \label{eq:dtp1_num_eq}
    \dualp{\partial_t (\bv^{\Delta t} + \nabla p_1^{\Delta t})}{\bw}_{(H^1_0(\widetilde{\Omega}))',H^1_0(\widetilde{\Omega})}
    &= \int_{\widetilde{\Omega}} \bbS^{\Delta t,+} : \nabla \bw
    + \int_{\widetilde{\Omega}} p_2^{\Delta t,+} \Div(\bw),
\end{align}
where $\bbS^{\Delta t,+} := (\bv^{\Delta t,+} \otimes \bv^{\Delta t,-}) - \nabla \bv^{\Delta t,+} - \bbF^{\Delta t,+} (\bbF^{\Delta t,+})^\top$.
\item Next, there exists a pressure $p\,\colon\, \widetilde{\Omega}\times(0,T)\to \bbR$ of the form $p=p_1+p_2$ with
\begin{align}
    \label{eq:p1_L2H2}
    p_1 &\in L^2(0,T; H^2(\widetilde{\Omega})),
    \\
    \label{eq:p2_L2L2}
    p_2 &\in L^2(0,T; L^2(\widetilde{\Omega})),
    \\
    \label{eq:dtp1_L2H10'}
    \partial_t (\bv + \nabla p_1) &\in L^2(0,T; (H^1_0(\widetilde{\Omega}; \mathbb{R}^2))'),
\end{align}
and, for all $\bw\in H^1_0(\widetilde{\Omega}; \mathbb{R}^2)$ and almost every $t\in(0,T)$, one has
\begin{align}
    \label{eq:dtp1_eq}
    \dualp{\partial_t (\bv + \nabla p_1)}{\bw}_{(H^1_0(\widetilde{\Omega}))',H^1_0(\widetilde{\Omega})}
    &= \int_{\widetilde{\Omega}} \bbS : \nabla \bw
    + \int_{\widetilde{\Omega}} p_2 \Div(\bw),
\end{align}
where $\bbS := (\bv \otimes \bv) - \nabla \bv - \overline{\bbF \bbF^\top}$.
\item Moreover, in the limit $\Delta t\to 0$, we have
\begin{alignat}{3}
    \label{eq:p1_conv_L2H2loc}
    p_1^{\Delta t(,\pm)} &\to p_1 
    \quad &&\text{strongly } \quad 
    &&\text{in } L^2(0,T;H^2_{loc}(\widetilde{\Omega})),
    \\
    \label{eq:p2_conv_L2L2}
    p_2^{\Delta t,+} &\rightharpoonup p_2
    \quad &&\text{weakly } \quad 
    &&\text{in } L^2(0,T;L^2(\widetilde{\Omega})).
\end{alignat}
\item In addition, 
\begin{align}
    \label{eq:p1_C(H1)}
    p_1 \in C([0,T]; H^1(\widetilde{\Omega})),
\end{align}
 and $p_1^{\Delta t}(\cdot,0) \to p_1(\cdot,0)$ strongly in $H^1(\widetilde{\Omega})$, as $\Delta t \to 0$.
\end{enumerate}
\end{lemma}

\begin{proof}
\textbf{Step 1:} We let $\widetilde{\Omega}$ be an arbitrary smooth domain with $\widetilde{\Omega}\subset \overline{\widetilde{\Omega}} \subset\Omega$. We consider the following time-discrete Stokes problems:
\begin{subequations}
\begin{alignat}{3}
    - \Delta \bu_1^n + \nabla p_1^n &= \bv^n
    \qquad &&\text{in } \widetilde{\Omega},
    \\
    \Div(\bu_1^n) &= 0 
    \qquad &&\text{in } \widetilde{\Omega},
    \\
    \bu_1^n |_{\partial\widetilde{\Omega}} &= \mathbf{0}
    \qquad &&\text{on } \partial\widetilde{\Omega},
\end{alignat}
\end{subequations}
for any $n\in\{0,\ldots,N_T\}$, and
\begin{subequations}
\begin{alignat}{3}
    - \Delta \bu_2^n + \nabla p_2^n &= \Div(\bbS^n)
    \qquad &&\text{in } \widetilde{\Omega},
    \\
    \Div(\bu_2^n) &= 0 
    \qquad &&\text{in } \widetilde{\Omega},
    \\
    \bu_2^n |_{\partial\widetilde{\Omega}} &= \mathbf{0}
    \qquad &&\text{on } \partial\widetilde{\Omega},
\end{alignat}
\end{subequations}
for any $n\in\{1,\ldots,N_T\}$, with $\bbS^n = (\bv^n \otimes \bv^{n-1}) - \nabla\bv^n - \bbF^n (\bbF^n)^\top$.

Since $\bv^n \in H^1_{0,\Div}(\Omega;\bbR^2)$ for any $n\in\{0,\ldots,N_T\}$, we obtain from \cite[Lemma A.1]{bulicek_2022_giesekus_2d} the existence of a unique weak solution $\bu_1^n \in H^1_{0,\Div}(\widetilde{\Omega};\bbR^2) \cap H^3(\widetilde{\Omega};\bbR^2)$, $p_1^n \in H^2(\widetilde{\Omega})$ with $\int_{\widetilde{\Omega}} p_1^n  = 0$ satisfying
\begin{align}
    \label{eq:u1_weak}
    \skp{\nabla \bu_1^n}{\nabla \bw}_{L^2(\widetilde{\Omega})}
    - \skp{p_1^n}{\Div(\bw)}_{L^2(\widetilde{\Omega})}
    &= \skp{\bv^n}{\bw}_{L^2(\widetilde{\Omega})},
    \\
    \label{eq:divu1_weak}
    \skp{\Div(\bu_1^n)}{q}_{L^2(\widetilde{\Omega})} &= 0,
\end{align}
for all $\bw\in H^1_0(\widetilde{\Omega};\bbR^2)$ and $q\in L^2(\widetilde{\Omega})$, and the elliptic regularity estimate
\begin{align}
    \label{eq:u1p1_stability}
    \norm{\bu_1^n}_{W^{m+2,2}(\widetilde{\Omega})} 
    + \norm{p_1^n}_{W^{m+1,2}(\widetilde{\Omega})} 
    \leq C \norm{\bv^n}_{W^{m,2}(\widetilde{\Omega})}
    \quad \forall\, m\in\{-1,0,1\}.
\end{align}
By squaring both sides, summing over all $n\in\{0,\ldots,N_T\}$ and by multiplication with $\Delta t$, we directly obtain
\begin{align}
    \label{eq:u1p1_stability_time}
    \Delta t \sum_{n=0}^{N_T} \norm{\bu_1^n}_{H^3(\widetilde{\Omega})}^2
    + \Delta t \sum_{n=0}^{N_T} \norm{p_1^n}_{H^2(\widetilde{\Omega})}^2 
    &\leq 
    C \Delta t \sum_{n=0}^{N_T} \norm{\bv^n}_{H^1}^2 .
\end{align}
With similar arguments, and by using the linearity of \eqref{eq:u1_weak}--\eqref{eq:divu1_weak}, we also have
\begin{align}
    \label{eq:u1p1_stability_time_differences}
    \sum_{n=1}^{N_T} \norm{\bu_1^n - \bu_1^{n-1}}_{H^2(\widetilde{\Omega})}^2 
    + \sum_{n=1}^{N_T} \norm{p_1^n - p_1^{n-1}}_{H^1(\widetilde{\Omega})}^2 
    &\leq 
    C \sum_{n=1}^{N_T} \norm{\bv^n - \bv^{n-1}}_{L^2}^2.
\end{align}
We note that the right-hand sides in \eqref{eq:u1p1_stability_time} and \eqref{eq:u1p1_stability_time_differences} are uniformly bounded in $\Delta t>0$ thanks to the time-discrete energy inequality \eqref{eq:stability_time_discrete_gronwall}. Hence, \eqref{eq:p1_num_L2H2} follows from \eqref{eq:fun_Delta_t}, \eqref{eq:fun_Delta_t_pm}, \eqref{eq:stability_time_continuous} and \eqref{eq:u1p1_stability_time}, \eqref{eq:u1p1_stability_time_differences}.

%
Similarly, since $\bbS^n \in L^2(\Omega;\bbR^{2\times 2})$ for any $n\in\{1,\ldots,N_T\}$, we have $\Div(\bbS^n)\in (H^1_0(\Omega;\bbR^2))'$. Then, by \cite[Lemma A.1]{bulicek_2022_giesekus_2d} we deduce the existence of a unique weak solution $\bu_2^n \in H^1_{0,\Div}(\widetilde{\Omega};\bbR^2)$, $p_2^n \in L^2(\widetilde{\Omega})$ with $\int_{\widetilde{\Omega}} p_2^n  = 0$ satisfying
\begin{align}
    \label{eq:u2_weak}
    \skp{\nabla \bu_2^n}{\nabla \bw}_{L^2(\widetilde{\Omega})}
    - \skp{p_2^n}{\Div(\bw)}_{L^2(\widetilde{\Omega})}
    &= - \skp{\bbS^n}{\nabla \bw}_{L^2(\widetilde{\Omega})},
    \\
    \label{eq:divu2_weak}
    \skp{\Div(\bu_2^n)}{q}_{L^2(\widetilde{\Omega})} &= 0,
\end{align}
for all $\bw\in H^1_0(\widetilde{\Omega};\bbR^2)$ and $q\in L^2(\widetilde{\Omega})$, and the estimate
\begin{align}
    \label{eq:u2p2_stability}
    \norm{\bu_2^n}_{H^{1}(\widetilde{\Omega})} 
    + \norm{p_2^n}_{L^{2}(\widetilde{\Omega})} 
    \leq C\norm{\Div(\bbS^n)}_{(H^1_0(\Omega))'} .
\end{align}
Again, by squaring both sides, summing over all $n\in\{1,\ldots,N_T\}$ and by multiplication with $\Delta t$, we observe
\begin{align}
    \label{eq:u2p2_stability_time}
    \Delta t \sum_{n=1}^{N_T} \norm{\bu_2^n}_{H^1(\widetilde{\Omega})}^2
    + \Delta t \sum_{n=1}^{N_T} \norm{p_2^n}_{L^2(\widetilde{\Omega})}^2 
    &\leq 
    C \Delta t \sum_{n=1}^{N_T} \norm{\Div(\bbS^n)}_{(H^1_0(\Omega))'}^2 ,
\end{align}
where the right-hand side is uniformly bounded in $\Delta t>0$ due to \eqref{eq:stability_time_discrete_gronwall}, which follows similarly to the proof of \eqref{eq:time_derivatives_discrete}. Thus, \eqref{eq:p2_num_L2L2} follows from \eqref{eq:fun_Delta_t_pm}, \eqref{eq:stability_time_continuous} and \eqref{eq:u2p2_stability_time}.

%
Let $\bw\in H^1_{0,\Div}(\widetilde{\Omega};\bbR^2)$ be arbitrary. By taking the difference of \eqref{eq:u1_weak} with the right-hand sides $\bv^n$ and $\bv^{n-1}$, and by multiplication with $\frac{1}{\Delta t}$, we have
\begin{align}
    \frac{1}{\Delta t} \skp{\bv^n-\bv^{n-1}}{\bw}_{L^2(\widetilde{\Omega})}
    &= \frac{1}{\Delta t} \skp{\nabla \bu_1^n-\nabla \bu_1^{n-1}}{\nabla \bw}_{L^2(\widetilde{\Omega})}.
\end{align}
Adding this equality to \eqref{eq:u2_weak} and noting \eqref{eq:v_time_discrete} gives
\begin{align}
    \nonumber
    0 &=
    \frac{1}{\Delta t} \skp{\bv^n-\bv^{n-1}}{\bw}_{L^2(\widetilde{\Omega})}
    - \skp{\bbS^n}{\nabla\bw}_{L^2(\widetilde{\Omega})}
    \\
    &= \skp{\nabla \left( \tfrac{1}{\Delta t} (\bu_1^n-\bu_1^{n-1}) + \bu_2^n \right) }{\nabla \bw}_{L^2(\widetilde{\Omega})}.
\end{align}
Since $\bw\in H^1_{0,\Div}(\widetilde{\Omega};\bbR^2)$ is arbitrary, we have that 
\begin{align*}
    \bu^n \coloneqq \left( \tfrac{1}{\Delta t} (\bu_1^n-\bu_1^{n-1}) + \bu_2^n \right) \  \in H^1_{0,\Div}(\widetilde{\Omega};\bbR^2)
\end{align*}
solves the following Stokes problem:
\begin{align}
    \skp{\nabla \bu^n}{\nabla \bw}_{L^2(\widetilde{\Omega})} 
    = 0 
    \qquad \forall\, \bw\in H^1_{0,\Div}(\widetilde{\Omega};\bbR^2),
\end{align}
and the uniqueness of solutions (cf.~\cite[Lemma A.1]{bulicek_2022_giesekus_2d}) implies that
\begin{align*}
    \bu^n = \tfrac{1}{\Delta t} (\bu_1^n-\bu_1^{n-1}) + \bu_2^n = \mathbf{0}
    \quad \text{a.e.~in } \widetilde{\Omega}.
\end{align*}
This allows us to deduce that
\begin{align}
    \label{eq:v1p1_weak} \nonumber
    &\frac{1}{\Delta t} \skp{(\bv^n-\bv^{n-1}) + (\nabla p_1^n - \nabla p_1^{n-1})}{\bw}_{L^2(\widetilde{\Omega})}
    - \skp{\bbS^n}{\nabla\bw}_{L^2(\widetilde{\Omega})}
    + \skp{p_2^n}{\Div(\bw)}_{L^2(\widetilde{\Omega})}
    \\
    &= \skp{\nabla \left( \tfrac{1}{\Delta t} (\bu_1^n-\bu_1^{n-1}) + \bu_2^n \right) }{\nabla \bw}_{L^2(\widetilde{\Omega})}
    = 0,
\end{align}
for all $\bw\in H^1_0(\widetilde{\Omega};\bbR^2)$. 

Similarly to \eqref{eq:time_derivatives_discrete}, one can show the following estimate with the help of \eqref{eq:v1p1_weak}, \eqref{eq:u1p1_stability_time} and \eqref{eq:u2p2_stability_time}:
\begin{align}
    \label{eq:u1p1_stability_time_derivative}
    & \nonumber
    \Delta t \sum_{n=1}^{N_T} 
    \nnorm{\frac{1}{\Delta t} (\bv^n-\bv^{n-1}) + \frac{1}{\Delta t}(\nabla p_1^n - \nabla p_1^{n-1})}_{(H^1_0(\widetilde{\Omega}))'}^2
    \\
    &\leq 
    C(T) \left(1 + \norm{\bv^0}_{L^2}^2 + \norm{\bbF^0}_{L^2}^2 + \Delta t \norm{\nabla\bv^0}_{L^2}^2  \right).
\end{align}
Hence, \eqref{eq:dtp1_num_L2H10'} follows from \eqref{eq:fun_Delta_t}, \eqref{eq:u1p1_stability_time_derivative} and the uniform bounds for the time-discrete initial data \eqref{eq:v0_init_stability}, \eqref{eq:F0_init_stability}.

%
\textbf{Step 2:}
By recalling \eqref{eq:u1p1_stability_time} and \eqref{eq:u1p1_stability_time_differences} and adopting the notation \eqref{eq:fun_Delta_t} and \eqref{eq:fun_Delta_t_pm},
we infer the existence of functions $\bu_1 \in L^2(0,T;H^3(\widetilde{\Omega};\bbR^2) \cap H^1_{0,\Div}(\widetilde{\Omega};\bbR^2))$ and $p_1 \in L^2(0,T;H^2(\widetilde{\Omega}))$ with $\int_{\widetilde{\Omega}} p_1 = 0$, such that for a non-relabeled subsequence
\begin{alignat}{3}
    \bu_1^{\Delta t(,\pm)} &\rightharpoonup \bu_1 
    \quad &&\text{weakly } \quad 
    &&\text{in } L^2 \left(0,T;H^3(\widetilde{\Omega};\bbR^2) \cap H^1_{0,\Div}(\widetilde{\Omega};\bbR^2) \right),
    \\
    p_1^{\Delta t(,\pm)} &\rightharpoonup p_1 
    \quad &&\text{weakly } \quad 
    &&\text{in } L^2(0,T;H^2(\widetilde{\Omega})),
    \\
    \partial_t (\bv^{\Delta t} + \nabla p_1^{\Delta t}) &\rightharpoonup \partial_t (\bv + \nabla p_1) 
    \quad &&\text{weakly } \quad 
    &&\text{in } L^2(0,T;(H^1_0(\widetilde{\Omega};\bbR^2))'),
\end{alignat}
as $\Delta t\to 0$.
This allows us to show pass to the limit in \eqref{eq:u1_weak} and \eqref{eq:divu1_weak} to deduce that
\begin{align}
    \label{eq:u1_weak_time}
    \skp{\nabla \bu_1}{\nabla \bw}_{L^2(\widetilde{\Omega})}
    - \skp{p_1}{\Div(\bw)}_{L^2(\widetilde{\Omega})}
    &= \skp{\bv}{\bw}_{L^2(\widetilde{\Omega})},
    \\
    \label{eq:divu1_weak_time}
    \skp{\Div(\bu_1)}{q}_{L^2(\widetilde{\Omega})} &= 0,
\end{align}
for almost all $t\in(0,T)$ and all $\bw\in H^1_0(\widetilde{\Omega};\bbR^2)$ and $q\in L^2(\widetilde{\Omega})$.

%
With similar arguments as in Section~\ref{sec:limit_vanishing_dt}, by recalling \eqref{eq:u2p2_stability_time} and adopting the notation \eqref{eq:fun_Delta_t_pm},
we infer the existence of functions $\bu_2 \in L^2(0,T;H^1_{0,\Div}(\widetilde{\Omega};\bbR^2))$ and $p_2 \in L^2(0,T;L^2(\widetilde{\Omega}))$ with $\int_{\widetilde{\Omega}} p_2 = 0$, such that for a non-relabeled subsequence
\begin{alignat}{3}
    \bu_2^{\Delta t,\pm} &\rightharpoonup \bu_2 
    \quad &&\text{weakly } \quad 
    &&\text{in } L^2 \left(0,T;H^1_{0,\Div}(\widetilde{\Omega};\bbR^2) \right),
    \\
    p_2^{\Delta t,\pm} &\rightharpoonup p_2 
    \quad &&\text{weakly } \quad 
    &&\text{in } L^2(0,T;L^2(\widetilde{\Omega})),
    \\
    \bbS^{\Delta t,\pm} &\rightharpoonup \bbS 
    \quad &&\text{weakly } \quad 
    &&\text{in } L^2(0,T;(H^1_0(\widetilde{\Omega};\bbR^{2\times 2}))'),
\end{alignat}
as $\Delta t\to 0$.
On noting \eqref{eq:stability_time_continuous}, we have that
\begin{align}
    \label{eq:u2_weak_time}
    \skp{\nabla \bu_2}{\nabla \bw}_{L^2(\widetilde{\Omega})}
    - \skp{p_2}{\Div(\bw)}_{L^2(\widetilde{\Omega})}
    &= \skp{\bbS}{\nabla \bw}_{L^2(\widetilde{\Omega})},
    \\
    \label{eq:divu2_weak_time}
    \skp{\Div(\bu_2)}{q}_{L^2(\widetilde{\Omega})} &= 0,
\end{align}
for almost all $t\in(0,T)$ and all $\bw\in H^1_0(\widetilde{\Omega};\bbR^2)$ and $q\in L^2(\widetilde{\Omega})$, where we have defined $\bbS := (\bv \otimes \bv) - \nabla \bv - \overline{\bbF \bbF^\top}$.

Proceeding in the same way as before, we pass to the limit in \eqref{eq:v1p1_weak} and deduce that
\begin{align*}
    \dualp{\partial_t (\bv + \nabla p_1)}{\bw}_{(H^1_0(\widetilde{\Omega}))',H^1_0(\widetilde{\Omega})}
    &= \int_{\widetilde{\Omega}} \bbS : \nabla \bw
    + \int_{\widetilde{\Omega}} p_2 \Div(\bw),
\end{align*}
for almost all $t\in(0,T)$ and all $\bw\in H^1_0(\widetilde{\Omega};\bbR^2)$, which is \eqref{eq:dtp1_eq}.

%
\textbf{Step 3:} The statement \eqref{eq:p2_conv_L2L2} has already be shown.
Concerning the strong convergence result \eqref{eq:p1_conv_L2H2loc}, we proceed as follows. By noting \eqref{eq:u1_weak}, \eqref{eq:divu1_weak} and \eqref{eq:u1_weak_time}, \eqref{eq:divu1_weak_time}, we observe that $(\bu_1^{\Delta t(,\pm)}-\bu_1)$, $(p_1^{\Delta t(,\pm)} - p_1)$ solve the following Stokes problem:
\begin{align}
    \label{eq:u1p1_weak_difference}
    \skp{\nabla (\bu_1^{\Delta t(,\pm)}-\bu_1)}{\nabla \bw}_{L^2(\widetilde{\Omega})}
    - \skp{p_1^{\Delta t(,\pm)} - p_1}{\Div(\bw)}_{L^2(\widetilde{\Omega})}
    &= \skp{\bv^{\Delta t(,\pm)} - \bv}{\bw}_{L^2(\widetilde{\Omega})},
    \\
    \label{eq:divu1_weak_difference}
    \skp{\Div(\bu_1^{\Delta t(,\pm)}-\bu_1)}{q}_{L^2(\widetilde{\Omega})} &= 0,
\end{align}
for almost all $t\in(0,T)$ and for all $\bw\in H^1_0(\widetilde{\Omega})$ and $q\in L^2(\widetilde{\Omega})$.
Hence, from \cite[Lemma A.1]{bulicek_2022_giesekus_2d} we obtain the stability estimate for almost all $t\in(0,T)$ and all $m\in\{-1,0,1\}$:
\begin{align}
    \label{eq:u1_stability}
    \norm{\bu_1^{\Delta t(,\pm)}-\bu_1}_{W^{m+2,2}(\widetilde{\Omega})} 
    + \norm{p_1^{\Delta t(,\pm)} - p_1}_{W^{m+1,2}(\widetilde{\Omega})} 
    \leq C \norm{\bv^{\Delta t(,\pm)} - \bv}_{W^{m,2}(\widetilde{\Omega})} .
\end{align}
On noting the strong convergence of $\bv^{\Delta t(,\pm)}$ to $\bv$ in $L^2(0,T;L^2(\Omega;\bbR^2))$, we deduce from \eqref{eq:u1_stability} that 
\begin{alignat}{3}
    \label{eq:p1_conv_L2H1}
    p_1^{\Delta t(,\pm)} &\to p_1 
    \quad &&\text{strongly } \quad 
    &&\text{in } L^2(0,T;H^1(\widetilde{\Omega})),
\end{alignat}
as $\Delta t\to 0$.

%
The next goal is to show that the Hessian of $p_1^{\Delta t,(\pm)}$, i.e., $\nabla^2 p_1^{\Delta t(,\pm)}$ converges strongly to $\nabla^2 p_1$ in $L^2(0,T;L^2_{loc}(\widetilde{\Omega};\bbR^{2\times 2}))$, as $\Delta t \to 0$.
On noting the regularity \eqref{eq:u1_stability} and integrating by parts in \eqref{eq:u1p1_weak_difference} and using the fundamental theorem of calculus of variations, we have that
\begin{align}
    -\Delta (\bu_1^{\Delta t(,\pm)}-\bu_1)
    + \nabla(p_1^{\Delta t(,\pm)} - p_1)
    &= \bv^{\Delta t(,\pm)} - \bv
    \quad \text{a.e.~in } \widetilde{\Omega}\times(0,T).
\end{align}
By taking the dot product with $\nabla\phi$, where $\phi\in H^1_0(\widetilde{\Omega})$ is arbitrary, and integrating over $\widetilde{\Omega}$, we have that
\begin{align*}
    \skp{\nabla(p_1^{\Delta t(,\pm)} - p_1)}{\nabla\phi}_{L^2(\widetilde{\Omega})} = 0,
\end{align*}
where we have used that $\bu_1^{\Delta t(,\pm)}$, $\bu_1$, $\bv^{\Delta t(,\pm)}$ and $\bv$ are divergence-free a.e.~in $\widetilde{\Omega} \times (0,T)$. 
By applying interior $H^2$-regularity for elliptic problems (cf.~\cite[Theorem~8.8]{gilbarg_trudinger_2001}), we have
\begin{align*}
	\norm{\nabla^2 p_1^{\Delta t(,\pm)} - \nabla^2 p_1}_{L^2(U)} 
	\leq C \norm{ p_1^{\Delta t(,\pm)} - p_1}_{H^1(\widetilde{\Omega})} ,
\end{align*}
for all open subsets $U\subset \overline{U} \subset \widetilde{\Omega}$. This implies
\begin{align*}
    \norm{\nabla^2 p_1^{\Delta t(,\pm)} - \nabla^2 p_1}_{L^2(0,T; L^2(U))} 
    \leq C \norm{ p_1^{\Delta t(,\pm)} - p_1}_{L^2(0,T;H^1(\widetilde{\Omega}))} ,
\end{align*}
where the right-hand side converges to zero due \eqref{eq:p1_conv_L2H1}. From this and \eqref{eq:p1_conv_L2H1}, we deduce \eqref{eq:p1_conv_L2H2loc}.

\textbf{Step 4:} We now show the final statement of the lemma. Let $t_1, t_2 \in [0,T]$ be arbitrary. We observe that the differences of functions $\bu_1(t_1) - \bu_1(t_2)$, $p_1(t_1) - p_1(t_2)$ uniquely solve the Stokes problem with right-hand side $\bv(t_1) - \bv(t_2)$. Thus, using \cite[Lemma A.1]{bulicek_2022_giesekus_2d} allows us to deduce the stability estimate
\begin{align*}
    \norm{ \bu_1(t_1) - \bu_1(t_2) }_{H^2(\widetilde{\Omega})}
    + 
    \norm{ p_1(t_1) - p_1(t_2) }_{H^1(\widetilde{\Omega})}
    \leq C \norm{ \bv(t_1) - \bv(t_2) }_{L^2(\Omega)}.
\end{align*}
Since $\bv$ belongs to $C([0,T]; L^2_{\Div}(\Omega;\bbR^2))$, the right-hand side vanishes for $t_2 \to t_1$, from which we deduce that $p_1 \in C([0,T]; H^1(\widetilde{\Omega}))$. 

Lastly, from \eqref{eq:u1_stability} and \eqref{eq:fun_Delta_t}, it follows that
\begin{align*}
    \norm{ p_1^{\Delta t}(0) - p_1(0) }_{H^1(\widetilde{\Omega})}
    \leq C \norm{ \bv^{\Delta t}(0) - \bv(0) }_{L^2(\Omega)} ,
\end{align*}
where the right-hand side converges to zero, as $\Delta t \to 0$, which is due to \eqref{eq:v0_init_convergence}.
\end{proof}

\subsubsection{Another differential inequality}

Next, we shall derive the following differential inequality for $\overline{\abs{\bbF}^2} - \abs{\bbF}^2$:
\begin{align}
    \label{eq:differential_inequality_2d}
    &- \int_{\Omega_T} \left( \overline{\abs{\bbF}^2} - \abs{\bbF}^2 \right)  \partial_t \phi 
    - \int_{\Omega_T} \left( \overline{\abs{\bbF}^2} - \abs{\bbF}^2 \right) \bv \cdot \nabla\phi
    \leq \int_{\Omega_T} L \left( \overline{\abs{\bbF}^2} - \abs{\bbF}^2  \right) \phi,
\end{align}
for all $\phi\in C^\infty_c(\Omega\times(-\infty,T))$ with $\phi\geq0$. Here, $L\in L^2(\Omega_T)$ is a fixed function defined as $L \coloneqq 1 + 2\abs{\nabla\bv} + \frac12 \abs{\bbF}^2$; see below.

We start by subtracting the previous differential (in-)equalities \eqref{eq:differential_inequality_overline|F|^2} and \eqref{eq:differential_inequality_|F|^2} for $\overline{|\bbF|^2}$ and $|\bbF|^2$, respectively. This gives us 
\begin{align}
    \label{eq:differential_inequality_2d_with_terms}
    \nonumber
    &- \int_{\Omega_T} (\overline{\abs{\bbF}^2} - \abs{\bbF}^2) \partial_t \phi 
    - \int_{\Omega_T} \left( (\overline{\abs{\bbF}^2} - \abs{\bbF}^2) \bv \right) \cdot \nabla\phi
    \\ \nonumber
    &\leq \int_{\Omega_T} (\overline{\abs{\bbF}^2} - \abs{\bbF}^2) \phi
    + 2 \dualp{\overline{\nabla\bv:(\bbF\bbF^\top)} - \overline{(\nabla\bv)\bbF} : \bbF}{\phi}_{\mathcal{M}(\overline{\Omega_T}), C(\overline{\Omega_T})}
    \\
    &\quad
    - \dualp{ \overline{\abs{\bbF\bbF^\top}^2} - \overline{\bbF\bbF^\top\bbF}:\bbF }{\phi}_{\mathcal{M}(\overline{\Omega_T}), C(\overline{\Omega_T})},
\end{align}
for all $\phi\in C^\infty_c(\Omega\times(-\infty,T))$ with $\phi\geq0$.

By noting \eqref{eq:F_conv_L4L4}, \eqref{eq:FFtF_conv}, \eqref{eq:|FFt|^2_conv}, we have that
\begin{align}
    \label{eq:bulicek_(6.67)}
    \nonumber
    &\dualp{ \overline{\abs{\bbF\bbF^\top}^2} - \overline{\bbF\bbF^\top\bbF}:\bbF }{\phi}_{\mathcal{M}(\overline{\Omega_T}), C(\overline{\Omega_T})}
    \\
    \nonumber
    &= \lim_{\Delta t\to 0} \int_{\Omega_T} \left( \abs{\bbF^{\Delta t,+}(\bbF^{\Delta t,+})^\top}^2 - (\bbF^{\Delta t,+}(\bbF^{\Delta t,+})^\top \bbF^{\Delta t,+}) : \bbF \right) \phi
    \\
    \nonumber
    &= \lim_{\Delta t\to 0} \int_{\Omega_T} \left( \bbF^{\Delta t,+}(\bbF^{\Delta t,+})^\top \bbF^{\Delta t,+} \right) : (\bbF^{\Delta t,+} - \bbF) \phi
    - (\bbF\bbF^\top \bbF) : (\bbF^{\Delta t,+} - \bbF) \phi
    \\
    \nonumber
    &= \lim_{\Delta t\to 0} \int_{\Omega_T} \left( \bbF^{\Delta t,+}(\bbF^{\Delta t,+})^\top \bbF^{\Delta t,+} - \bbF\bbF^\top \bbF \right) : (\bbF^{\Delta t,+} - \bbF) \phi
    \\*
    &\geq 0,
\end{align}
where the inequality in the last step results from the non-negativity of $\phi$ and the monotonicity of the matrix function $\bbG\mapsto \bbG\bbG^\top \bbG$ for all $\bbG\in\bbR^{2\times 2}$; see \cite[Lemma~4.2]{bulicek_2022_giesekus_2d}.

The goal of the remainder is to prove the inequality
\begin{align}
    \label{eq:bulicek_(6.68)}
    \dualp{\overline{\nabla\bv:(\bbF\bbF^\top)} - \overline{(\nabla\bv)\bbF} : \bbF}{\phi}_{\mathcal{M}(\overline{\Omega_T}), C(\overline{\Omega_T})}
    \leq \int_{\Omega_T} \tilde{L} (\overline{\abs{\bbF}^2} - \abs{\bbF}^2) \phi,
\end{align}
for some $\tilde{L} \in L^2(\Omega_T)$ being a fixed function that will be specified as $\tilde{L} = \abs{\nabla\bv} + \frac14 \abs{\bbF}^2$; see below.
Then, we will obtain the desired inequality \eqref{eq:differential_inequality_2d} from \eqref{eq:differential_inequality_2d_with_terms}, \eqref{eq:bulicek_(6.67)} and \eqref{eq:bulicek_(6.68)}.

We first decompose the left-hand side of \eqref{eq:bulicek_(6.68)} into the sum $I_1 + I_2 + I_3$ with
\begin{align*}
    I_1 &\coloneqq \dualp{\overline{\nabla\bv:(\bbF\bbF^\top)} - \nabla\bv : \overline{\bbF \bbF^\top}}{\phi}_{\mathcal{M}(\overline{\Omega_T}), C(\overline{\Omega_T})},
    \\
    I_2 &\coloneqq \int_{\Omega_T} \left({\nabla\bv : \overline{\bbF \bbF^\top} - \nabla\bv : (\bbF \bbF^\top)}\right){\phi},
    \\
    I_3 &\coloneqq \int_{\Omega_T} \left({\nabla\bv : (\bbF \bbF^\top) - \overline{(\nabla\bv)\bbF} : \bbF}\right) {\phi} .
\end{align*}
For $I_2$, using \eqref{eq:F_conv_L4L4}, \eqref{eq:FF_conv}, \eqref{eq:|F|^2_conv} and the Cauchy--Schwarz inequality, we have that
\begin{align}
    \label{eq:I2_estimate}
    \nonumber
    I_2 
    &= \lim_{\Delta t\to 0} \int_{\Omega_T} \left( \nabla\bv : (\bbF^{\Delta t,+} (\bbF^{\Delta t,+})^\top) 
    - \nabla\bv : (\bbF\bbF^\top) \right) \phi
    \\
    \nonumber
    &= \lim_{\Delta t\to 0} \int_{\Omega_T} \nabla\bv : \left(  (\bbF^{\Delta t,+} - \bbF) (\bbF^{\Delta t,+} - \bbF)^\top \right) \phi
    \\
    \nonumber
    &\leq 
    \lim_{\Delta t\to 0} \int_{\Omega_T} \abs{\nabla\bv} \abs{\bbF^{\Delta t,+} - \bbF}^2  \phi
    \\
    &= \int_{\Omega_T} \abs{\nabla\bv} \left( \overline{\abs{\bbF}^2} - \abs{\bbF}^2 \right) \phi.
\end{align}
For $I_3$, we have by a similar argument,  based on \eqref{eq:F_conv_L4L4}, \eqref{eq:v_conv_L2H1}, \eqref{eq:|F|^2_conv}, \eqref{eq:|gradv|^2_conv}, together with the Cauchy--Schwarz inequality and Young's inequality, 
that
\begin{align}
    \label{eq:I3_estimate}
    \nonumber
    I_3 
    &= \lim_{\Delta t\to 0} \int_{\Omega_T} \left( (\nabla\bv^{\Delta t,+}) (\bbF-\bbF^{\Delta t,+}) \right) : (\phi\bbF)
    \\ \nonumber
    &= \lim_{\Delta t\to 0} \int_{\Omega_T} \left( (\nabla\bv^{\Delta t,+} - \nabla\bv) (\bbF-\bbF^{\Delta t,+}) \right) : (\phi\bbF)
    \\ \nonumber
    &\leq \lim_{\Delta t\to 0} \int_{\Omega_T} \left( \abs{\nabla\bv^{\Delta t,+} - \nabla\bv}^2 + \frac14 \abs{\bbF-\bbF^{\Delta t,+}}^2 \abs{\bbF}^2 \right) \phi
    \\
    &= \dualp{\overline{\abs{\nabla\bv}^2} - \abs{\nabla\bv}^2}{\phi}_{\mathcal{M}(\overline{\Omega_T}), C(\overline{\Omega_T})}
    + \frac14 \int_{\Omega_T} \abs{\bbF}^2 \left( \overline{\abs{\bbF}^2} - \abs{\bbF}^2 \right) \phi.
\end{align}

With the help of our (weak(-$*$)) compactness results from Section~\ref{sec:limit_vanishing_dt}, we aim to show for $I_1$ that
\begin{align}
    \label{eq:I1_estimate}
    I_1 \leq \dualp{ \abs{\nabla\bv}^2 - \overline{\abs{\nabla\bv}^2} }{\phi}_{\mathcal{M}(\overline{\Omega_T}), C(\overline{\Omega_T})} .
\end{align}
In contrast to \cite[(6.72)]{bulicek_2022_giesekus_2d}, we have an inequality here, which is sufficient for the subsequent analysis.
Importantly, this step can only be carried out in two space dimensions, where the required regularity provided by Lemma~\ref{lemma:local_pressures} is available. In three dimensions, the necessary regularity is missing, and this estimate \eqref{eq:I1_estimate} for $I_1$ cannot be shown.

To prove \eqref{eq:I1_estimate}, we proceed as follows.
For any $\phi \in C^\infty_c(\Omega\times(-\infty,T))$, we can choose a smooth subset $\widetilde{\Omega} \subset\overline{\widetilde{\Omega}} \subset\Omega$ such that $\phi \in C^\infty_c(\widetilde{\Omega}\times(-\infty,T))$. In other words, the compact support of $\phi(\cdot,t)$ is contained in the interior of a fixed smooth subset of $\Omega$.
Now fix $t \in (0,T)$. The general idea would be to test \eqref{eq:v_weaklimit} with $\bv \phi$ and \eqref{eq:v_time_continuous} with $\bv^{\Delta t,+}\phi$. However, these test functions are not divergence-free and hence not admissible. To address this, we construct local pressures $p = p_1 + p_2$ and $p^{\Delta t,+} = p_1^{\Delta t,+} + p_2^{\Delta t,+}$ using Lemma~\ref{lemma:local_pressures}.

With these in hand, we subtract \eqref{eq:dtp1_eq}, tested with $\bw = (\bv + \nabla p_1)\phi$, from \eqref{eq:dtp1_num_eq}, tested with $\bw = (\bv^{\Delta t,+} + \nabla p_1^{\Delta t,+})\phi$, and then integrate the resulting identity over $(0,T)$. Both test functions are admissible here.
Then, in the limit $\Delta t\to 0$, we obtain:
\begin{align*}
    I_1 &= \lim_{\Delta t\to 0}
    \int_{\Omega_T} \left( \nabla \bv^{\Delta t,+} : (\bbF^{\Delta t,+} (\bbF^{\Delta t,+})^\top) - \nabla\bv : (\bbF^{\Delta t,+} (\bbF^{\Delta t,+})^\top) \right) \phi
    \\
    &= \lim_{\Delta t\to 0} \sum_{j=1}^9 J_j^{\Delta t},
\end{align*}
where we introduced the following terms:
\begin{align*}
    J_1^{\Delta t} 
    &\coloneqq
    - \int_0^T \dualp{\partial_t (\bv^{\Delta t} + \nabla p_1^{\Delta t})}{\phi (\bv^{\Delta t,+} + \nabla p_1^{\Delta t,+})}_{(H^1_0(\widetilde{\Omega}))',H^1_0(\widetilde{\Omega})}
    \\*
    &\quad+ \int_0^T \dualp{\partial_t (\bv + \nabla p_1)}{\phi (\bv + \nabla p_1)}_{(H^1_0(\widetilde{\Omega}))',H^1_0(\widetilde{\Omega})} ,
    \\
    J_2^{\Delta t} 
    &\coloneqq
    \int_{\Omega_T} (\bv^{\Delta t,+} \otimes \bv^{\Delta t,-}) : (\phi \nabla \bv^{\Delta t,+}) 
    - \int_{\Omega_T} (\bv \otimes \bv) : (\phi \nabla \bv) ,
    \\
    J_3^{\Delta t} 
    &\coloneqq
    \int_{\Omega_T} (\bv^{\Delta t,+} \otimes \bv^{\Delta t,-}) : (\bv^{\Delta t,+} \otimes\nabla\phi) 
    - \int_{\Omega_T} (\bv \otimes \bv) : (\bv \otimes\nabla\phi) ,
    \\
    J_4^{\Delta t} 
    &\coloneqq
    \int_{\Omega_T} (\bv^{\Delta t,+} \otimes \bv^{\Delta t,-}) : (\nabla p_1^{\Delta t,+} \otimes\nabla\phi + \phi \nabla^2p_1^{\Delta t,+}) 
    \\
    &\quad
    - \int_{\Omega_T} (\bv \otimes \bv) : (\nabla p_1 \otimes\nabla\phi + \phi \nabla^2p_1) ,
    \\
    J_5^{\Delta t} 
    &\coloneqq 
    \int_{\Omega_T} \left( \abs{\nabla\bv}^2 - \abs{\nabla\bv^{\Delta t,+}}^2 \right) \phi ,
    \\
    J_6^{\Delta t} 
    &\coloneqq
    - \int_{\Omega_T} \nabla\bv^{\Delta t,+} : \left( (\bv^{\Delta t,+}+\nabla p_1^{\Delta t,+}) \otimes \nabla\phi + \phi \nabla^2 p_1^{\Delta t,+}\right) 
    \\
    &\quad 
    + \int_{\Omega_T} \nabla\bv : \left( (\bv+\nabla p_1) \otimes \nabla\phi + \phi \nabla^2 p_1\right),
    \\
    J_7^{\Delta t} 
    &\coloneqq
    \int_{\Omega_T} p_2^{\Delta t,+} \left( (\bv^{\Delta t,+}+\nabla p_1^{\Delta t,+}) \cdot\nabla\phi + \phi \Delta p_1^{\Delta t,+} \right) 
    \\
    &\quad 
    - \int_{\Omega_T} p_2 \left( (\bv+\nabla p_1) \cdot\nabla\phi + \phi \Delta p_1 \right),
    \\
    J_8^{\Delta t} 
    &\coloneqq
    \int_{\Omega_T} \bbF^{\Delta t,+} (\bbF^{\Delta t,+})^\top : \left( (\bv^{\Delta t,+}+\nabla p_1^{\Delta t,+})\otimes\nabla\phi \right) 
    \\
    &\quad 
    - \int_{\Omega_T} \bbF^{\Delta t,+} (\bbF^{\Delta t,+})^\top : \left( (\bv+\nabla p_1)\otimes\nabla\phi \right),
    \\
    J_9^{\Delta t} 
    &\coloneqq
    \int_{\Omega_T} \bbF^{\Delta t,+} (\bbF^{\Delta t,+})^\top : (\phi \nabla^2 p_1^{\Delta t,+})
    - \int_{\Omega_T} \bbF^{\Delta t,+} (\bbF^{\Delta t,+})^\top : (\phi \nabla^2 p_1).
\end{align*}
It follows from straightforward arguments based on \eqref{eq:v_conv_L2H1}, \eqref{eq:v_conv_L2Lr}, \eqref{eq:v_otimes_v_conv}, \eqref{eq:p1_conv_L2H2loc}, \eqref{eq:p2_conv_L2L2}, \eqref{eq:FF_conv} and \eqref{eq:|gradv|^2_conv} that
\begin{align*}
    0 = \lim_{\Delta t\to 0} J_3^{\Delta t} 
    = \lim_{\Delta t\to 0} J_4^{\Delta t} 
    = \lim_{\Delta t\to 0} J_6^{\Delta t} 
    = \lim_{\Delta t\to 0} J_7^{\Delta t} 
    = \lim_{\Delta t\to 0} J_8^{\Delta t} 
    = \lim_{\Delta t\to 0} J_9^{\Delta t} ,
\end{align*}
and
\begin{align*}
    \lim_{\Delta t\to 0} J_5^{\Delta t} 
    &= \dualp{ \abs{\nabla\bv}^2 - \overline{\abs{\nabla\bv}^2} }{\phi}_{\mathcal{M}(\overline{\Omega_T}), C(\overline{\Omega_T})}.
\end{align*}
Moreover, using integration by parts over $\Omega$ and noting $\Div(\bv^{\Delta t,\pm}) = \Div(\bv)=0$ a.e.~in $\Omega_T$, we obtain
\begin{align*}
    0 = \lim_{\Delta t\to 0} J_2^{\Delta t}.
\end{align*}
Concerning $J_1^{\Delta t}$, 
we note that it follows from \eqref{eq:p1_C(H1)}, \eqref{eq:v_C(L2)} and \eqref{eq:dtp1_L2H10'} that $(\bv + \nabla p_1) \in C([0,T]; L^2(\widetilde{\Omega};\bbR^2))$. This allows us to integrate by parts with respect to the time variable, so that
\begin{align}
    \label{eq:proof_J1dt_1}
    \nonumber
    &\int_0^T \dualp{\partial_t (\bv + \nabla p_1)}{\phi (\bv + \nabla p_1)}_{(H^1_0(\widetilde{\Omega}))',H^1_0(\widetilde{\Omega})}
    \\
    &= -  \frac12 \skp{\abs{\bv(0) + \nabla p_1(0)}^2 }{\phi(0)}_{L^2(\widetilde{\Omega})}
    -  \frac12 \int_0^T \skp{\abs{\bv + \nabla p_1}^2 }{\partial_t \phi}_{L^2(\widetilde{\Omega})}.
\end{align}
Similarly, we perform discrete integration by parts with respect to the time variable in the time-discrete formulation. In fact, by denoting $\bz^n := \bv^n + \nabla p_1^n$ and specifying $a_n := |\bz^n|^2$ and $b_n := \phi^n$ in \eqref{eq:discrete_integration_by_parts_in_time}, we have that
\begin{align}
    \label{eq:proof_J1dt_2}
    \nonumber
    &- \int_0^T \dualp{\partial_t (\bv^{\Delta t} + \nabla p_1^{\Delta t})}{\phi (\bv^{\Delta t,+} + \nabla p_1^{\Delta t,+})}_{(H^1_0(\widetilde{\Omega}))',H^1_0(\widetilde{\Omega})}
    = -\sum_{n=1}^{N_T} \skp{\bz^n - \bz^{n-1}}{\bz^n \phi^n}_{L^2(\widetilde{\Omega})}
    \\
    \nonumber
    &= -\frac12 \sum_{n=1}^{N_T} \skp{ \abs{\bz^n}^2 - \abs{\bz^{n-1}}^2}{\phi^n}_{L^2(\widetilde{\Omega})}
    - \frac12 \sum_{n=1}^{N_T} \skp{ \abs{\bz^n - \bz^{n-1}}^2}{\phi^n}_{L^2(\widetilde{\Omega})}
    \\
    \nonumber
    &=  \frac12 \sum_{n=1}^{N_T} \skp{ \abs{\bz^{n-1}}^2}{\phi^n - \phi^{n-1}}_{L^2(\widetilde{\Omega})}
    - \frac12 \skp{\abs{\bz^{N_T}}^2}{\phi^{N_T}}_{L^2(\widetilde{\Omega})}
    \\*
    \nonumber
    &\quad
    + \frac12 \skp{\abs{\bz^{0}}^2}{\phi^{0}}_{L^2(\widetilde{\Omega})}
    - \frac12 \sum_{n=1}^{N_T} \skp{ \abs{\bz^n - \bz^{n-1}}^2}{\phi^n}_{L^2(\widetilde{\Omega})}
    \\
    \nonumber
    &=  \frac12 \int_0^T \skp{\abs{\bz^{\Delta t,-}}^2}{\partial_t \phi^{\Delta t}}_{L^2(\widetilde{\Omega})}
    - \frac12 \skp{\abs{\bz^{\Delta t}(T)}^2}{\phi^{\Delta t}(T)}_{L^2(\widetilde{\Omega})}
    \\
    &\quad
    + \frac12 \skp{\abs{\bz^{\Delta t}(0)}^2}{\phi^{\Delta t}(0)}_{L^2(\widetilde{\Omega})}
    - \frac12 \int_0^T \skp{ \frac{1}{\Delta t} \abs{\bz^{\Delta t,+} - \bz^{\Delta t,-}}^2}{\phi}_{L^2(\widetilde{\Omega})} ,
\end{align}
where $\phi^n \coloneqq \frac{1}{\Delta t} \int_{t^{n-1}}^{t^n} \phi(\cdot,t) \dt$, and $\phi^{\Delta t}$ is the continuous, piecewise affine interpolant in time of $\phi^n$ based on the notation \eqref{eq:fun_Delta_t}.
Here we note that, as $\Delta t \to 0$,
\begin{alignat*}{3}
    \phi^{\Delta t,+}(\cdot,0) &= \frac{1}{\Delta t} \int_{-\Delta t}^0 \phi(\cdot,t)\dt \to \phi(\cdot,0), \quad &&\text{uniformly on $\overline{\Omega}$, }
    \\
    \bz^{\Delta t}(\cdot,0) &\to \bv(\cdot,0) + \nabla p_1(\cdot,0) \quad &&\text{strongly in $L^2(\widetilde{\Omega};\bbR^2)$,}
    \\
    \partial_t \phi^{\Delta t}(\bx,t) &= \frac{1}{\Delta t}(\phi(\bx,t)-\phi(\bx,t- \Delta t )) \to \partial_t\phi(\bx,t)\quad &&\text{uniformly in $(\bx,t)\in\overline{\Omega_T}$} .
\end{alignat*}
Adding \eqref{eq:proof_J1dt_1} and \eqref{eq:proof_J1dt_2}, and noting the non-negativity of $\phi$, 
we deduce that
\begin{align*}
    \lim_{\Delta t \to 0} J_1^{\Delta t} \leq 0.
\end{align*} 
This finally shows \eqref{eq:I1_estimate}.

Then, using the inequalities \eqref{eq:I2_estimate}, \eqref{eq:I3_estimate} and \eqref{eq:I1_estimate} for $I_1$, $I_2$ and $I_3$, respectively, we obtain
\begin{align}
    \nonumber
    &\dualp{\overline{\nabla\bv:(\bbF\bbF^\top)} - \overline{(\nabla\bv)\bbF} : \bbF}{\phi}_{\mathcal{M}(\overline{\Omega_T}), C(\overline{\Omega_T})}
    \\* \nonumber
    &=I_1 + I_2 + I_3 
    \\ \nonumber
    &\leq 
    \dualp{ \abs{\nabla\bv}^2 - \overline{\abs{\nabla\bv}^2} }{\phi}_{\mathcal{M}(\overline{\Omega_T}), C(\overline{\Omega_T})}
    + \int_{\Omega_T} \abs{\nabla\bv} \left( \overline{\abs{\bbF}^2} - \abs{\bbF}^2 \right) \phi
    \\ \nonumber
    &\quad 
    + \dualp{\overline{\abs{\nabla\bv}^2} - \abs{\nabla\bv}^2}{\phi}_{\mathcal{M}(\overline{\Omega_T}), C(\overline{\Omega_T})}
    + \frac14 \int_{\Omega_T} \abs{\bbF}^2 \left( \overline{\abs{\bbF}^2} - \abs{\bbF}^2 \right) \phi
    \\
    &= \int_{\Omega_T} \tilde{L}  \left( \overline{\abs{\bbF}^2} - \abs{\bbF}^2 \right) \phi,
\end{align}
with $\tilde{L} \coloneqq \abs{\nabla\bv} + \frac14 \abs{\bbF}^2 \in L^2(\Omega_T)$, which is \eqref{eq:bulicek_(6.68)}.

Combining this with \eqref{eq:differential_inequality_2d_with_terms} and \eqref{eq:bulicek_(6.67)}, we have finally shown the desired result \eqref{eq:differential_inequality_2d}.

\subsubsection{Renormalization and compactness}

Since $\bbF^{\Delta t(,\pm)} \rightharpoonup \bbF$ weakly in $L^2(\Omega_T; \mathbb{R}^{2\times 2})$, the weak lower semicontinuity of the $L^2$ norm guarantees that $\overline{\abs{\bbF}^2} - \abs{\bbF}^2 \geq 0$ almost everywhere in $\Omega_T$. 
Applying the renormalization result \cite[Proposition 1.5]{bulicek_2025_giesekus_3d} to \eqref{eq:differential_inequality_2d} then yields the equality $\overline{\abs{\bbF}^2} = \abs{\bbF}^2$ a.e.~in $\Omega_T$.
From this equality, we have the strong convergence $\bbF^{\Delta t(,\pm)} \to \bbF$ strongly in $L^2(\Omega_T; \mathbb{R}^{2\times 2})$ (i.e., \eqref{eq:F_conv_strong}), which in turn implies convergence almost everywhere in $\Omega_T$ after extracting a further (non-relabeled) subsequence.
Consequently, this ensures that the limits of products coincide with the products of limits a.e.~in $\Omega_T$:
\begin{align*}
    \overline{\bbF\bbF^\top} = \bbF\bbF^\top,
    \quad
    \overline{\bbF\bbF^\top \bbF} = \bbF\bbF^\top \bbF,
    \quad
    \overline{(\nabla\bv)\bbF} = (\nabla\bv) \bbF.
\end{align*}
Hence, the limit functions $\bv$ and $\bbF$ solve \eqref{eq:v_weak}, \eqref{eq:F_weak}.

\subsection{Positivity of the determinant} \label{sec:positivity_F}

To complete the proof of Theorem~\ref{theorem:convergence} and hence to verify Definition~\ref{def:weak_solution} fully, we let $d=2$ and suppose that $\det(\bbF_0) > 0$ almost everywhere in $\Omega$ and $\ln \det(\bbF_0) \in L^1(\Omega)$. Then, it remains to show the positivity of $\det(\bbF)$ almost everywhere in $\Omega_T$. 
For the proof, we make use of computations for the limiting system \eqref{eq:system_weaklimit} presented in \cite[Section 6.5]{bulicek_2022_giesekus_2d}, which rely on mollification, cut-off techniques, uniform estimates, and suitable limit passages. In particular, these arguments yield the estimate \eqref{eq:F_lndet}, from which one then infers that $\det(\bbF) > 0$ almost everywhere in~$\Omega_T$.
This completes the proof of Theorem~\ref{theorem:convergence} and we also deduce that $\bv$ and $\bbF$ form a weak solution to \eqref{eq:system} in the sense of Definition~\ref{def:weak_solution}.


%
%
%
%
%
\section{Numerical results}
\label{sec:numerical_results}

In this section, we present numerical results for the scheme \eqref{eq:system_FE} in two dimensions. We first perform convergence tests using manufactured solutions, following which we shall investigate the performance of the scheme by simulating flow through a planar domain exhibiting $4{:}1$ contraction.

\subsection{Discretization and solver setup}
For all computations, we utilize the lowest-order Taylor--Hood pair of elements ($\mathcal{P}_2/\mathcal{P}_1$) for the velocity and pressure approximation, and continuous piecewise linear elements ($\mathcal{P}_1$) for the deformation gradient (corresponding to $k = m = 1$ in \eqref{eq:system_FE}). 
To ensure energy stability, the convective term in the evolution equation \eqref{eq:F_FE} is implemented in the skew-symmetric form \eqref{eq:convective_skewsymm}.

The nonlinear system \eqref{eq:system_FE} is solved at each time step using Newton's method. The iteration is terminated when the $\ell^\infty$-norm of the Newton increment drops below a given tolerance, which is chosen to be $10^{-12}$. We observe that the Newton solver typically terminates within 2 to 4 iterations. The resulting linear sub-problems are solved exactly using the direct \texttt{UMFPACK LU} solver provided by the \texttt{PETSc} framework \cite{petsc-user-ref_2021}.
It is worth mentioning that the Newton iterations remained robust even for high Weissenberg numbers (e.g., $\textit{Wi}=8$).

\subsection{Convergence tests in two dimensions}
We perform convergence tests on the unit square domain $\Omega=(0,1)^2 \subset\bbR^2$ with a final time $T=0.1$. The physical parameters are set to one: $\rho=\nu=\mu=\lambda=1$.
To construct a manufactured solution, we introduce appropriate source terms into \eqref{eq:v} and \eqref{eq:F} such that the exact solution $(\tilde\bv, \tilde p, \tilde\bbF)$ matches:
\begin{align*}
    \tilde\bv(\mathbf{x},t) &= \mathrm{e}^{-t}
    \left(\begin{smallmatrix}
        x_1^2 (x_1-1)^2 x_2 (x_2-1) (2x_2-1)\\
        -x_1 (x_1-1) (2x_1-1) x_2^2 (x_2-1)^2
    \end{smallmatrix}\right),
    \quad
    \tilde p(\mathbf{x},t) = \mathrm{e}^{-t} (2x_1-1)(2x_2-1),
    \\
    \tilde\bbF(\mathbf{x},t) &= 
    \left(\begin{smallmatrix}
    1 & 0\\ 0 & 1
    \end{smallmatrix}\right)
    + \tfrac{1}{6} \mathrm{e}^{-t} \cos(4\pi x_1)\cos(4\pi x_2)
    \left(\begin{smallmatrix}
       1 & 0\\
        0 & -1
    \end{smallmatrix}\right),
\end{align*}
for all $(\bx,t)\in \Omega\times(0,T)$.

The discretization parameters are chosen as $\Delta t_\ell = \tfrac{1}{5} T \cdot 2^{-\ell}$ with $\ell\in\{1,\ldots,7\}$, and mesh sizes $h_j = 2^{-j}$ with $j\in\{2,\ldots,7\}$. The computational domain is discretized using uniform Friedrichs--Keller triangulations with $2^j \times 2^j$ vertices. The discrete initial data for \eqref{eq:system_FE} is obtained via $L^2$ projection of the exact solution onto the finite element spaces, using quadrature rules exact for polynomials up to degree eight.

While the alternative energy-stable form \eqref{eq:convective_Lambda} theoretically limits the accuracy of the spatial approximation of $\bbF$ to $\mathcal{O}(h)$ in the $L^2$-norm (see \eqref{eq:error_Lambda}), numerical experiments showed nearly identical error curves when compared to the skew-symmetric form \eqref{eq:convective_skewsymm}. A similar observation holds for the variant \eqref{eq:convective_higherorder}, which is not guaranteed to be energy-stable for $k=m=1$. Consequently, we limit our presentation to results obtained using the skew-symmetric form.

\begin{figure}[ht!]
\centering
\includegraphics[width=0.30\textwidth]{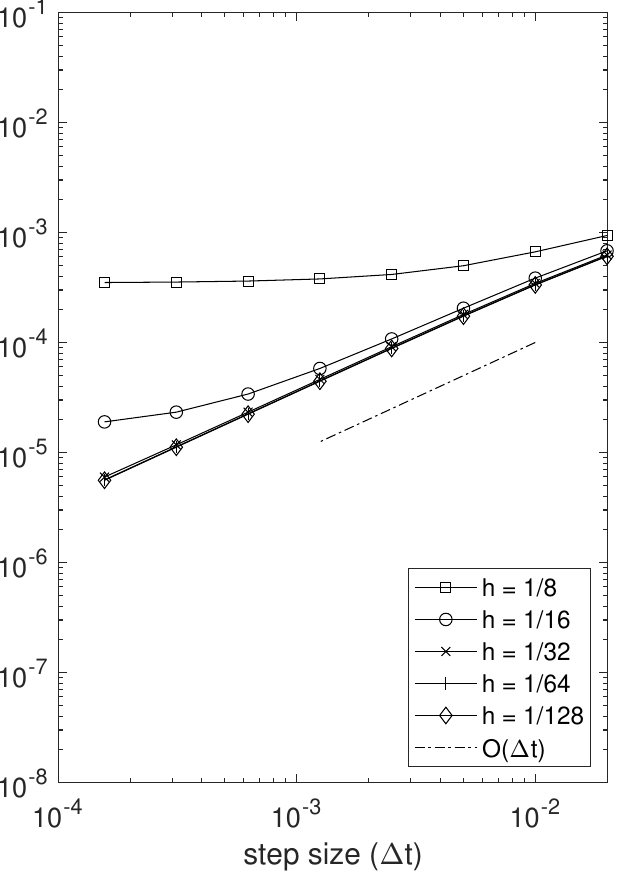}
\includegraphics[width=0.30\textwidth]{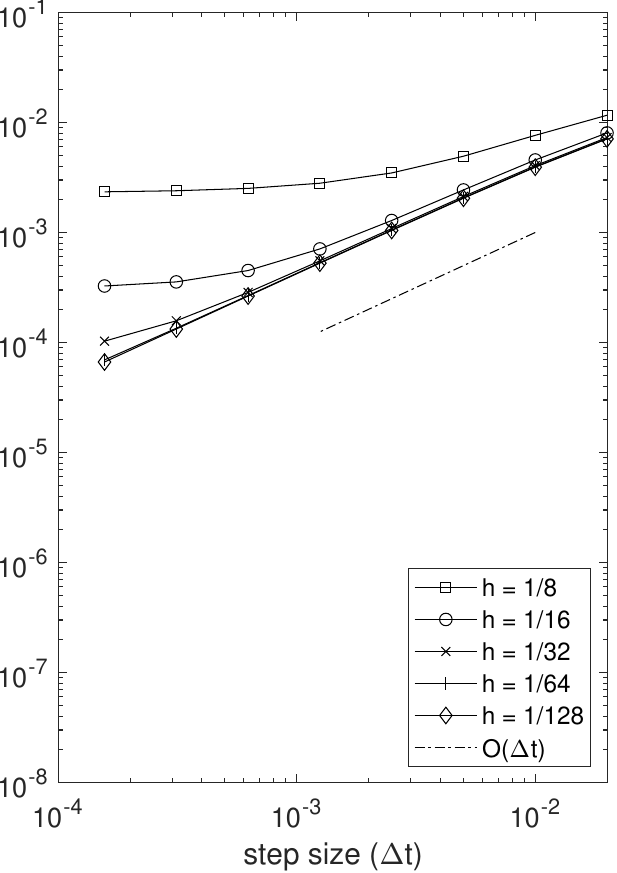}
\includegraphics[width=0.30\textwidth]{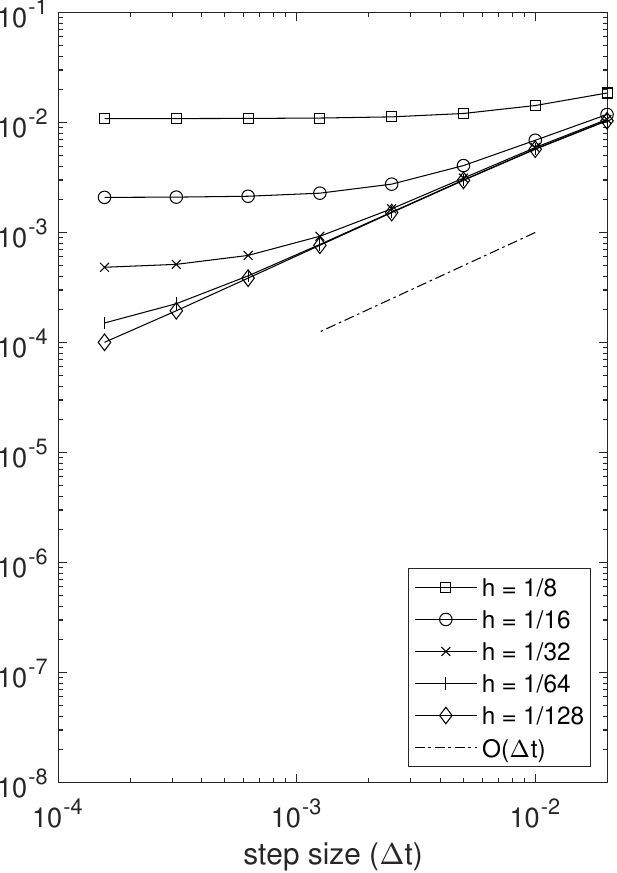}
\caption{Studying the temporal convergence of the $L^2(\Omega_T)$ errors for the velocity $\bv$ (left), the pressure $p$ (center), and the deformation gradient $\bbF$ (right).}
\label{fig:errors_L2L2_time}      
\end{figure}

\begin{figure}[ht!]
\centering
\includegraphics[width=0.30\textwidth]{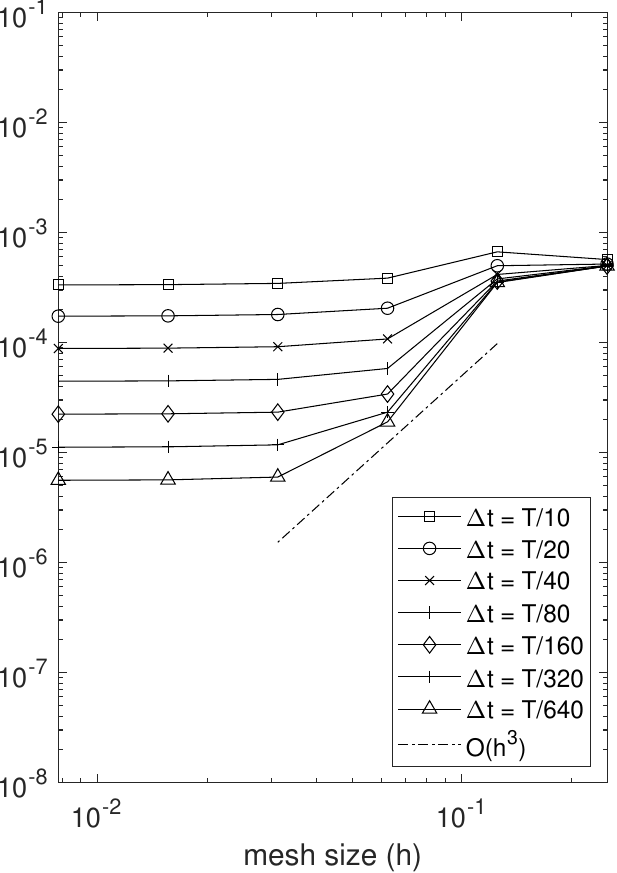}
\includegraphics[width=0.30\textwidth]{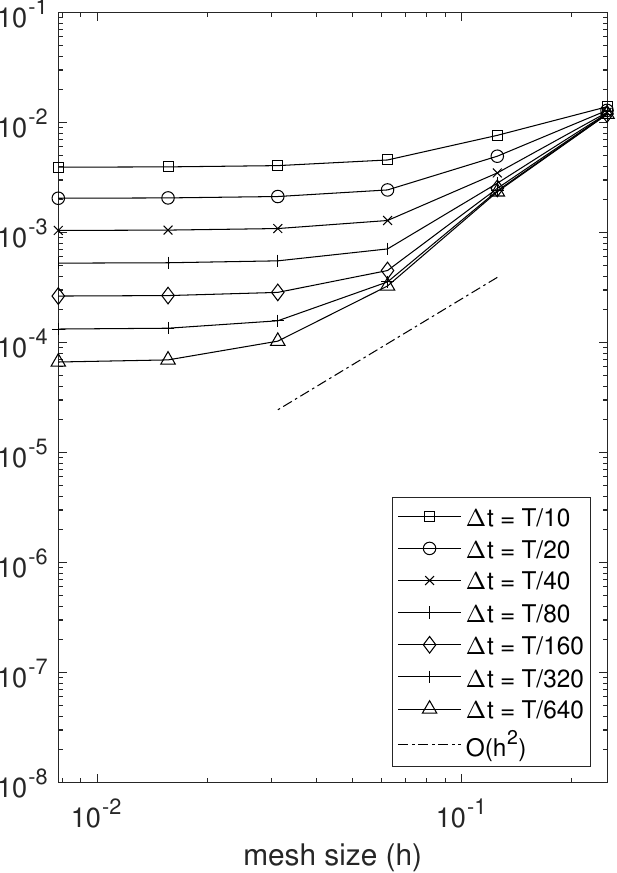}
\includegraphics[width=0.30\textwidth]{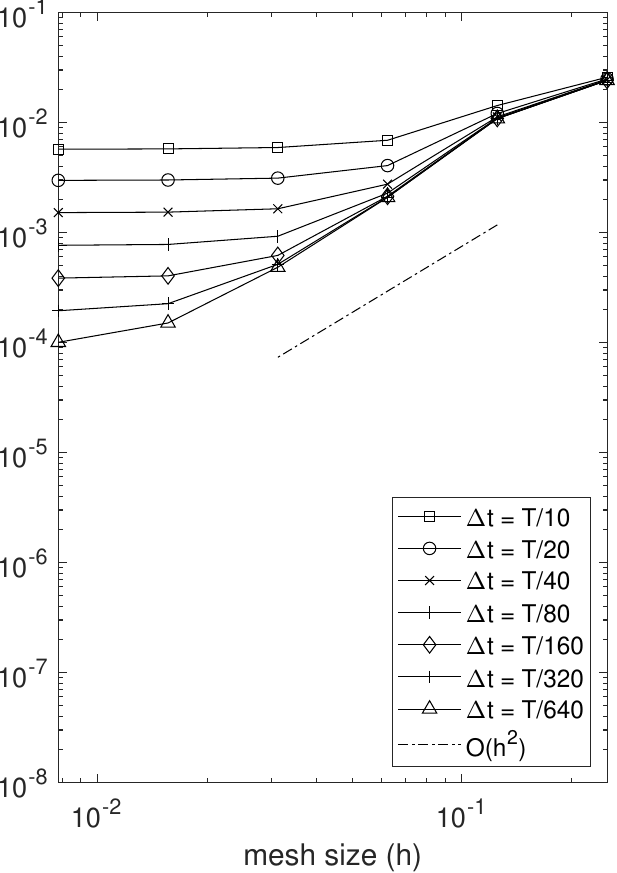}
\caption{Studying the spatial convergence of the $L^2(\Omega_T)$ errors for the velocity $\bv$ (left), the pressure $p$ (center), and the deformation gradient $\bbF$ (right).}
\label{fig:errors_L2L2_space}      
\end{figure}

The theoretically expected convergence rates in the $L^2(\Omega_T)$ space-time norm are $\mathcal{O}(h^3+\Delta t)$ for velocity, and $\mathcal{O}(h^2+\Delta t)$ for both pressure and the deformation gradient. 
Figure~\ref{fig:errors_L2L2_time} displays the error decay with respect to the time step size $\Delta t$. We clearly observe an asymptotic regime with $\mathcal{O}(\Delta t)$ convergence where temporal errors dominate, followed by a plateau where spatial discretization errors prevail.
Complementarily, Figure~\ref{fig:errors_L2L2_space} shows the error as a function of the mesh size $h$. The results confirm spatial convergence rates of $\mathcal{O}(h^3)$ for the velocity and $\mathcal{O}(h^2)$ for the pressure and deformation gradient. These rates agree with the expected approximation properties of the chosen finite element spaces and the first-order Euler-type time discretization.

\subsection{Planar 4:1 contraction}
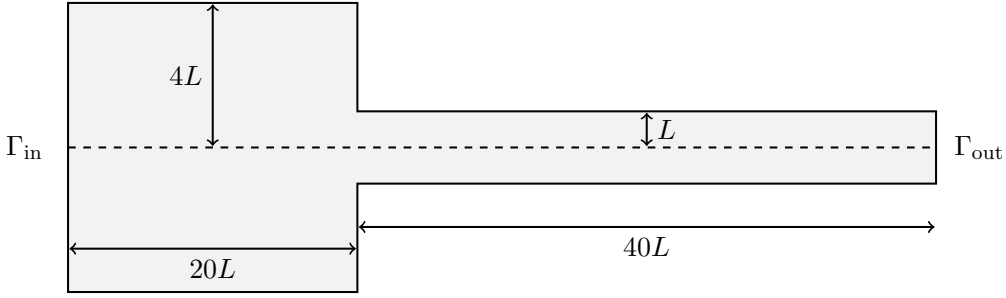
\begin{figure}
\centering
\resizebox{!}{0.25\textwidth}{%
\begin{tikzpicture}[scale=1.0, thick]
\def\links{-4}
\def\rechts{8}
\def\hoch{2}
\coordinate (A) at (\links,\hoch);
\coordinate (B) at (0,\hoch);
\coordinate (C) at (0,\hoch/4);
\coordinate (D) at (\rechts,\hoch/4);
\coordinate (E) at (\rechts,-\hoch/4);
\coordinate (F) at (0,-\hoch/4);
\coordinate (G) at (0,-\hoch);
\coordinate (H) at (\links,-\hoch);
\draw[fill=gray!10] (A)--(B)--(C)--(D)--(E)--(F)--(G)--(H)--cycle;
\draw[dashed] (\links,0)--(\rechts,0);
\draw[<->] (\links/2,0.02) -- (\links/2,\hoch-0.02)
node[midway, left] {$4L$};
\draw[<->] (\rechts/2,0.02) -- (\rechts/2,\hoch/4-0.02)
node[midway, right] {$L$};
\draw[<->] (\links+0.02,-\hoch+0.6) -- (-0.02,-\hoch+0.6)
node[midway, below] {$20L$};
\draw[<->] (0.02,-\hoch/4-0.6) -- (\rechts-0.02,-\hoch/4-0.6)
node[midway, below] {$40L$};
\node at (\links-0.6,0) {$\Gamma_\text{in}$};
\node at (\rechts+0.6,0) {$\Gamma_\text{out}$};
\end{tikzpicture}
}%
\caption{Schematic sketch of the $4{:}1$ planar contraction geometry.}
\label{fig:geometry_4:1_contraction}
\end{figure}

We now consider the benchmark problem of viscoelastic flow through a planar $4{:}1$ contraction, a standard test case for numerical schemes for viscoelastic models \cite{ alves_2021_viscoelastic_review, niethammer_marschall_bothe_2019_viscoelastic, pimenta_alves_2017_viscoelastic}. 
The geometry consists of an upstream channel of half-width $4L$ and length $20L$, and a downstream channel of half-width $L$ and length $40L$. The contraction occurs at the plane $x_1=0$. For a schematic sketch, we refer to Figure~\ref{fig:geometry_4:1_contraction}.

Boundary conditions are imposed as follows.
At the inlet $\Gamma_\text{in} = \{(-20L, x_2)^\top \in\bbR^2 \mid x_2 \in (-4L,4L)\}$, a parabolic flow profile
\begin{align*}
    \bv(\bx,t) = \left(V_{\text{in}} \left(1- \abs{\frac{x_2}{4L}}^2\right), 0 \right)^\top
    \quad \forall\, \bx \in \Gamma_\text{in}, \ t>0,
\end{align*}
is prescribed for some given $V_{\text{in}}>0$, together with $\bbF=\bbI$, resulting in zero elastic stress at the inlet. At the outlet $\Gamma_\text{out} = \{(40L, x_2)^\top \in\bbR^2 \mid x_2 \in (-L,L)\}$, the flow is assumed to be fully developed, and the zero normal flux boundary condition
\begin{align*}
    (\nu\nabla\bv - p\bbI + \mu(\bbF\bbF^\top - \bbI))\bn = \mathbf{0}
\end{align*}
is applied there, where $\bn=(1,0)^\top$ is the outer unit normal to $\Gamma_\text{out}$. 
On the remaining part $\partial\Omega \setminus (\Gamma_\text{in} \cup \Gamma_\text{out})$ of the boundary, the no-slip boundary condition $\bv=\mathbf{0}$ is applied.

Due to our use of the skew-symmetric form for the convective terms, the boundary integrals do not vanish on $\Gamma_\text{out}$ where Dirichlet conditions are absent. Therefore, to maintain consistency with the variational formulation, we include the following terms in \eqref{eq:v_FE} and \eqref{eq:F_FE}:
\begin{align*}
    \frac{\rho}{2}\int_{\Gamma_\text{out}} (\bn\cdot \bv_h^{n-1}) (\bv_h^n \cdot \bw_h)
    \quad 
    \text{and}
    \quad
    \frac12 \int_{\Gamma_\text{out}} (\bn\cdot \bv_h^{n-1}) (\bbF_h^n : \bbG_h) .
\end{align*}
The simulation starts from a fluid at rest ($\bv_0=\mathbf{0}$) with no elastic deformation ($\bbF_0=\bbI$). No symmetry is enforced along the center-line $x_2=0$ in order to allow possible flow asymmetries to develop.

\begin{figure}
\centering
\includegraphics[width=0.7\textwidth]{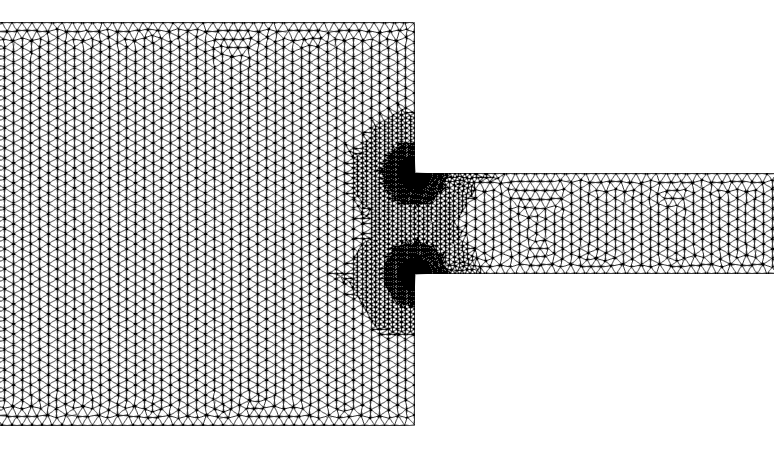}
\caption{Zoomed view of the computational mesh near the contraction plane $x_1=0$. The triangulation shows the unstructured base mesh ($h=1/20$) with five levels of local refinement applied near the re-entrant corners $(0, \pm L)^\top$ to resolve the stress singularity.}
\label{fig:mesh}
\end{figure}

The domain $\Omega$ is triangulated with an unstructured mesh created with \texttt{Gmsh} \cite{gmsh_library}, with local mesh size $h=\frac{1}{20}$. For better resolution at the re-entrant corners $(0,\pm L)^\top$, we subsequently perform local refinement by bisection of all cells within a ball of radius $r=0.6\cdot 2^{-i}$, $i\in\{0,1,2,3,4\}$, using the routines provided by \texttt{FEniCs} \cite{fenics_book_2012}. The refined part of the grid close to the re-entrant corners is shown in Figure~\ref{fig:mesh}.

For the experiment, we consider the model parameters based on \cite{pimenta_alves_2017_viscoelastic}. For the nondimensional parameters in \eqref{eq:nondim_parameters}, we choose the characteristic length $x_c=L$ and the characteristic velocity $v_c=4 V_{\text{in}}$, which is approximately the velocity magnitude in the center of the outlet tunnel. Hence, the Reynolds number, the Weissenberg number, and the elastic-to-viscous viscosity ratio are, respectively,
\begin{align*}
    \textit{Re} = \frac{4 \rho V_{\text{in}} L}{\nu + \lambda},
    \quad
    \textit{Wi} = \frac{4 \lambda V_{\text{in}}}{\mu L},
    \quad
    \alpha = \frac{\lambda}{\nu + \lambda} \in (0,1).
\end{align*}
In the following, we fix $\alpha=8/9$ and $\textit{Re}=0.01$, while the Weissenberg number \textit{Wi} is varied. We achieve this by setting
\begin{align*}
    \nu = \frac19 \lambda, \quad
    \rho = \frac{1}{18} \lambda, \quad 
    \mu = 1, \quad
    V_{\text{in}} = 0.1, \quad
    L = 0.5, 
\end{align*}
and adjusting the relaxation parameter $\lambda \geq 0$.
In the case $\lambda=0$, the fluid is purely Newtonian, i.e., the deformation gradient is fixed to the identity matrix ($\bbF=\bbI$) and only the Navier--Stokes equations \eqref{eq:v_FE} and \eqref{eq:div_FE} are solved. Moreover, we then set $\rho=0.05 \nu$ and $\nu=1$ to obtain $\textit{Re}=0.01$.

The relevant time scale for this physical setup is the elastic relaxation time $t_\lambda = \lambda/\mu$ \cite{pimenta_alves_2017_viscoelastic}. To ensure that we capture the full dynamic behaviour, we set the simulation duration to $T=10$, which satisfies $T \geq t_\lambda$ for all tested parameters. We vary the relaxation parameter $\lambda \in \{0, 0.125, 1.25, 3.75, 6.25, 10\}$ to obtain Weissenberg numbers $\textit{Wi} \in \{0, 0.1, 1, 3, 5, 8\}$.

\begin{figure}
\centering
\subfloat[Standard scheme ($\Delta t=0.01$).]{\label{fig:energy_total_a}
    \includegraphics[width=0.48\textwidth]{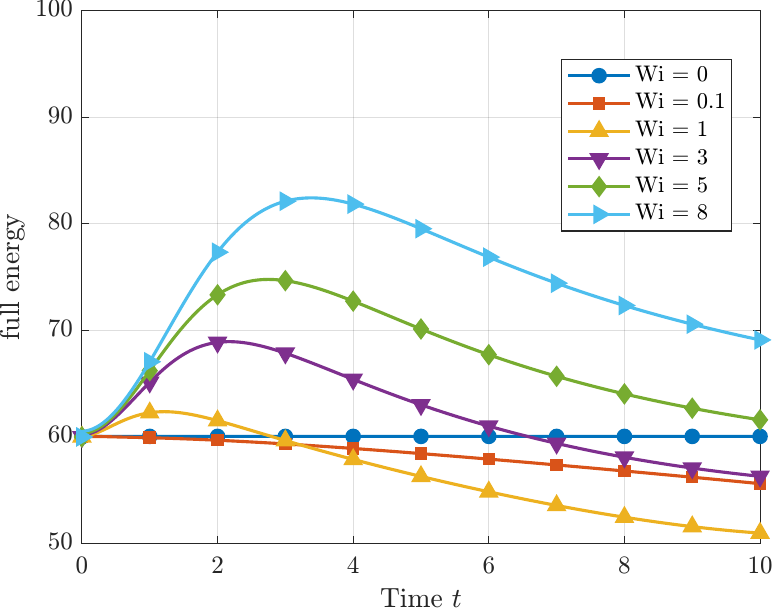}}
\hfill
\subfloat[Smaller time step size ($\Delta t=0.0025$).]{\label{fig:energy_total_b}
    \includegraphics[width=0.48\textwidth]{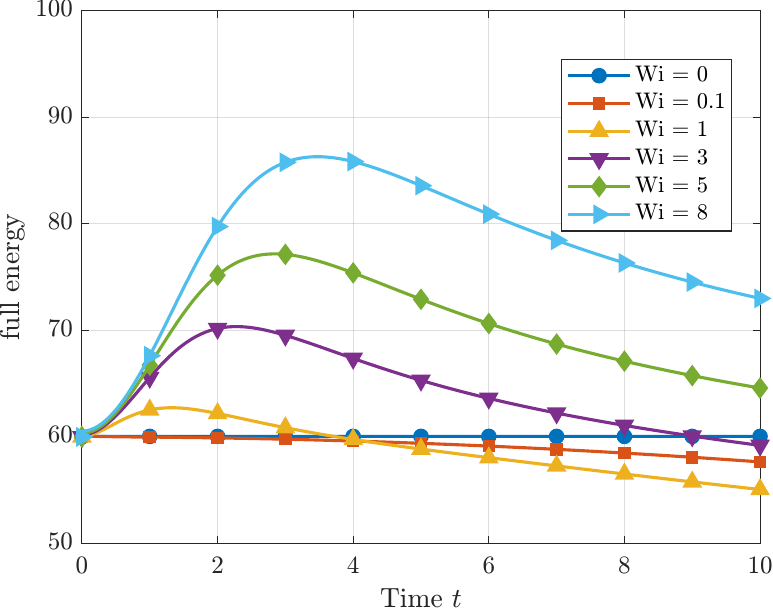}}

\vspace{0.2cm}

\subfloat[Stress diffusion $\sim (\Delta t)^2$.]{\label{fig:energy_total_c}
    \includegraphics[width=0.48\textwidth]{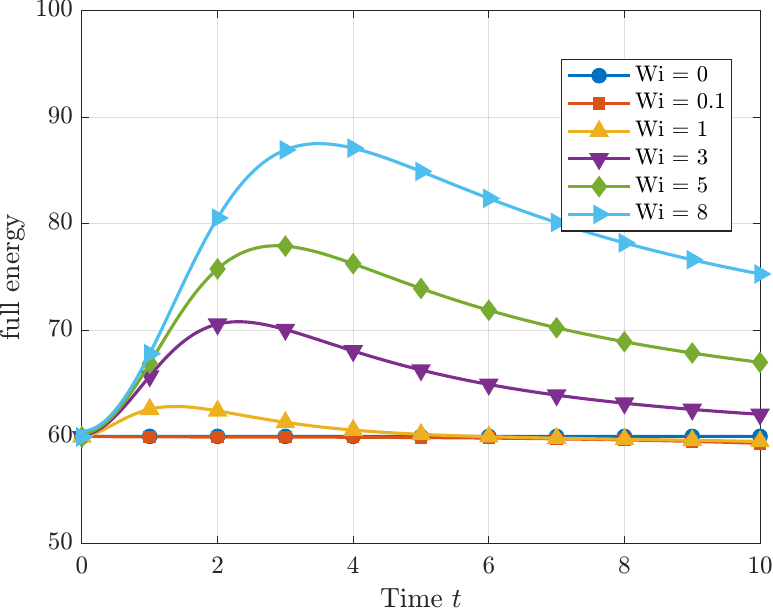}}
\hfill
\subfloat[No stress diffusion.]{\label{fig:energy_total_d}
    \includegraphics[width=0.48\textwidth]{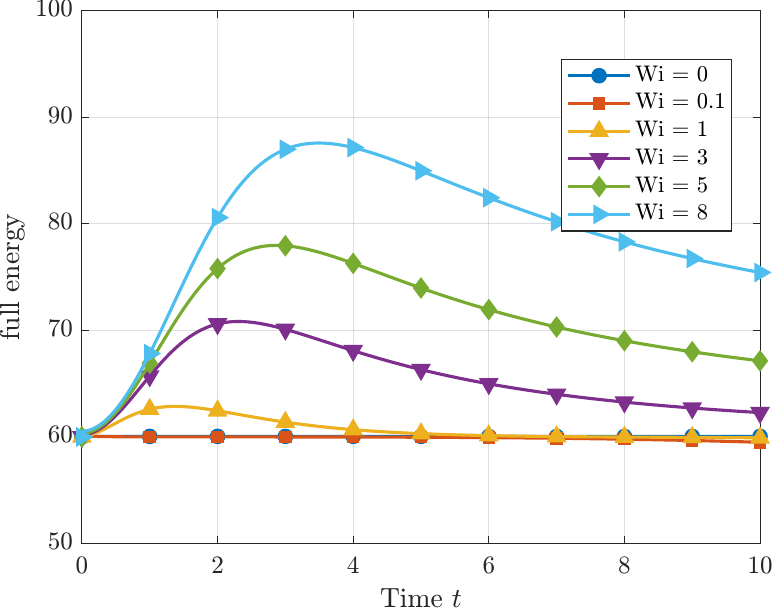}
}

\caption{Evolution of the energy $\int_\Omega \left( \frac\rho2 |\bv|^2 + \frac\mu2 |\bbF|^2 \right)$ for different Weissenberg numbers $\textit{Wi} \in \{0, 0.1, 1, 3, 5, 8\}$.
\textbf{(a)}~The standard scheme \eqref{eq:system_FE}.
\textbf{(b)}~The standard scheme with a refined time step size $\Delta t = 0.0025$.
\textbf{(c)}~The scheme with stress diffusion scaled by $(\Delta t)^2$ instead of $\Delta t$.
\textbf{(d)}~The scheme without stress diffusion.
In case (a), the stabilization term results in strong dissipation and decaying total energy in the long run. Reducing the time step size in (b) does not fully eliminate this dissipation. However, in case (c), the additional dissipation appears negligible, and the curves coincide with the non-stabilized case (d). All results agree with the stability result \eqref{eq:stability_FE}.
}
\label{fig:energy}
\end{figure}

\begin{figure}
\centering
\subfloat[Standard scheme ($\Delta t=0.01$).]{\label{fig:energy_log_a}
    \includegraphics[width=0.48\textwidth]{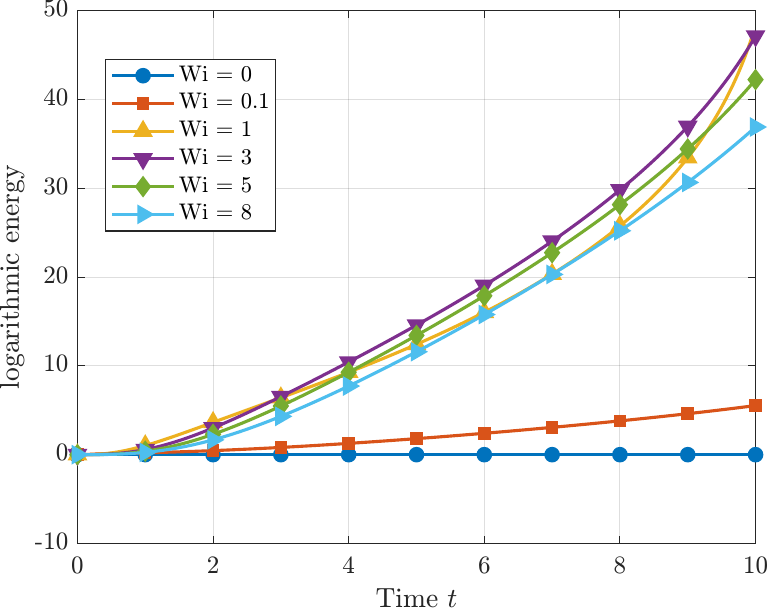}}
\hfill
\subfloat[Smaller time step size ($\Delta t=0.0025$).]{\label{fig:energy_log_b}
    \includegraphics[width=0.48\textwidth]{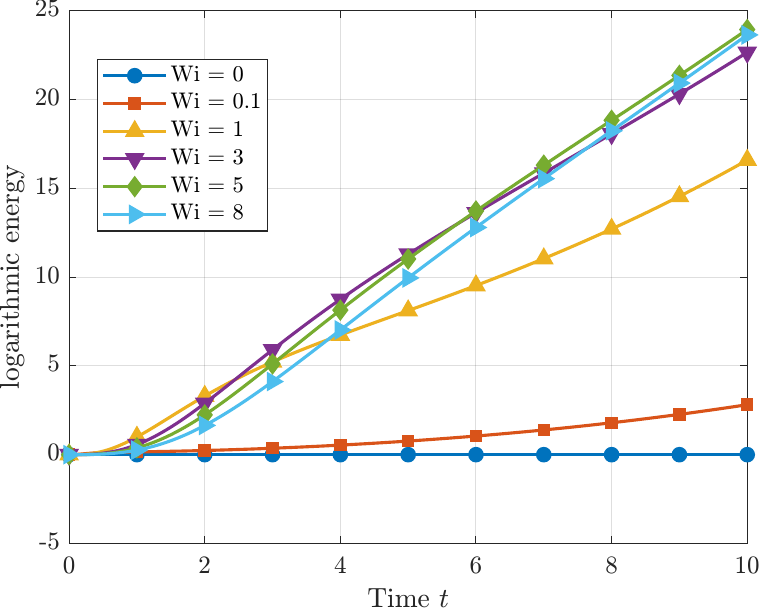}}

\vspace{0.2cm}

\subfloat[Stress diffusion $\sim (\Delta t)^2$.]{\label{fig:energy_log_c}
    \includegraphics[width=0.48\textwidth]{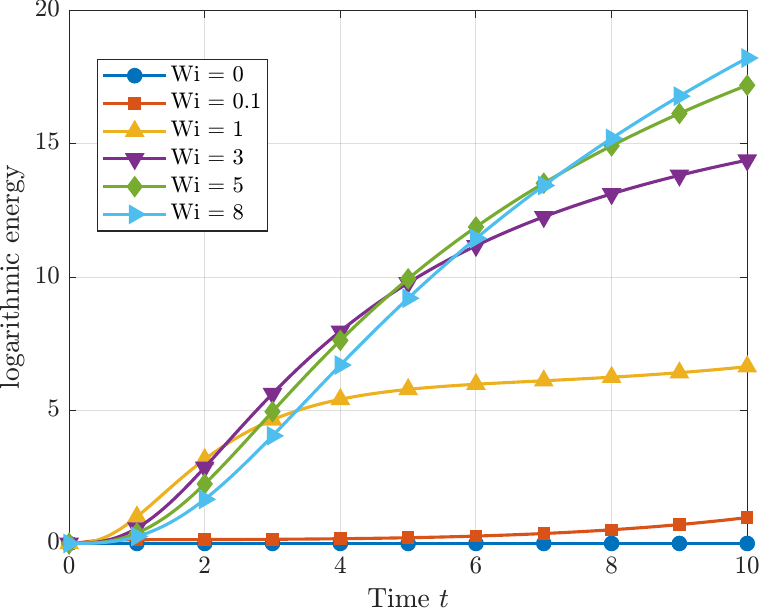}}
\hfill
\subfloat[No stress diffusion.]{\label{fig:energy_log_d}
    \includegraphics[width=0.48\textwidth]{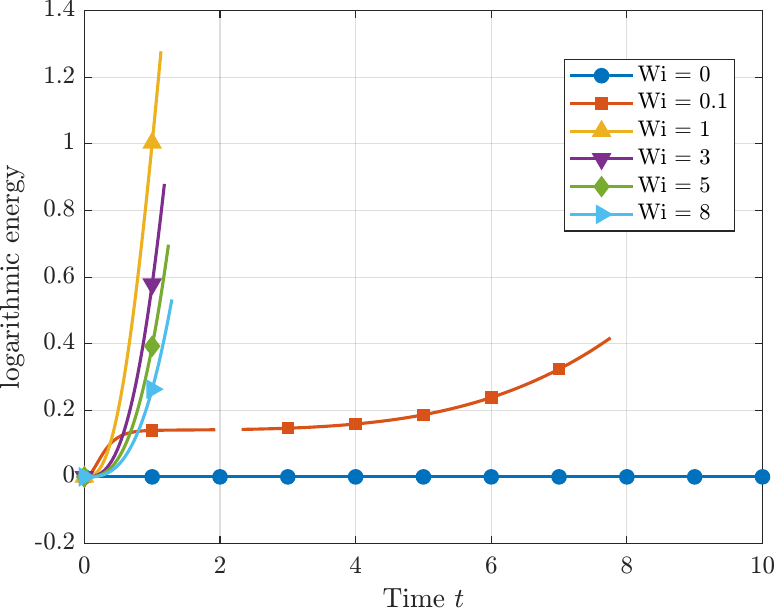}}

\caption{Time evolution of the logarithmic energy term $-\int_\Omega \frac{\mu}{2} \ln\det(\bbF\bbF^\top)$ for Weissenberg numbers $\textit{Wi} \in \{0, 0.1, 1, 3, 5, 8\}$.
\textbf{(a)}~The standard scheme \eqref{eq:system_FE}.
\textbf{(b)}~The standard scheme with a refined time step size $\Delta t = 0.0025$.
\textbf{(c)}~The scheme with stress diffusion scaled by $(\Delta t)^2$ instead of $\Delta t$.
\textbf{(d)}~The scheme without stress diffusion.
In all elastic cases ($\textit{Wi} > 0$), the logarithmic energy grows over time. Note that in (d), the curves terminate early because the determinant becomes non-positive, causing the logarithmic energy to blow up. The different stabilization in (c) and the smaller time steps in (b) reduce this growth and prevent numerical failure.
}
\label{fig:energy_log}
\end{figure}

We investigate four numerical configurations to analyze the energy stability based on \eqref{eq:stability_FE}. First, we employ the standard scheme \eqref{eq:system_FE} with a time step size $\Delta t = 10^{-2}$ (case a). Second, we test the same scheme with a refined time step size $\Delta t = 2.5 \cdot 10^{-3}$ (case b). Third, we modify the stabilization term in \eqref{eq:F_FE} so that the stress diffusion scales with $(\Delta t)^2$ instead of $\Delta t$, using the original time step size $\Delta t = 10^{-2}$ (case c). Finally, we compare these against the scheme without any stress diffusion (case d), again with $\Delta t = 10^{-2}$.

Figure~\ref{fig:energy} shows the evolution of the total energy $\int_\Omega \left( \frac\rho2 |\bv|^2 + \frac\mu2 |\bbF|^2 \right)$. In the standard case (a), the linear stabilization introduces strong artificial dissipation, causing the total energy to decay rapidly. Reducing the time step size in (b) does not fully eliminate this artificial dissipation. However, when the stress diffusion is scaled by $(\Delta t)^2$ in case (c), the artificial dissipation becomes negligible. The energy curves in (c) virtually coincide with the non-stabilized case (d). In all cases, the total energy remains finite, which agrees with the stability result \eqref{eq:stability_FE}.

\begin{figure}
\centering
\subfloat[Standard scheme ($\textit{Wi} = 0.1$).]{\label{fig:stress_a}
    \includegraphics[width=0.43\textwidth]{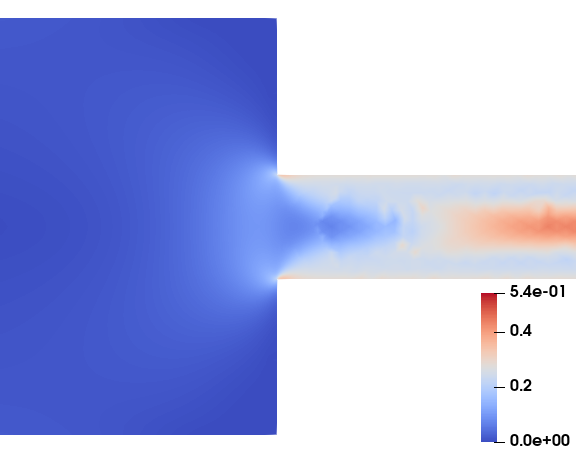}}
\hfill
\subfloat[Standard scheme ($\textit{Wi} = 1$).]{\label{fig:stress_b}
    \includegraphics[width=0.43\textwidth]{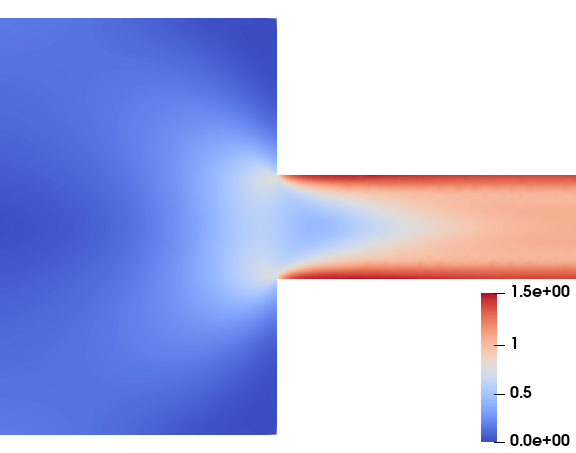}}


\subfloat[Standard scheme ($\textit{Wi} = 8$).]{\label{fig:stress_c}
    \includegraphics[width=0.43\textwidth]{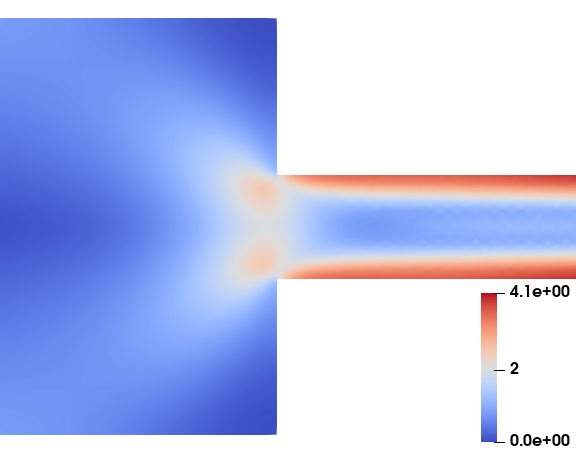}}
\hfill
\subfloat[Smaller time step size ($\textit{Wi} = 8$).]{\label{fig:stress_d}
    \includegraphics[width=0.43\textwidth]{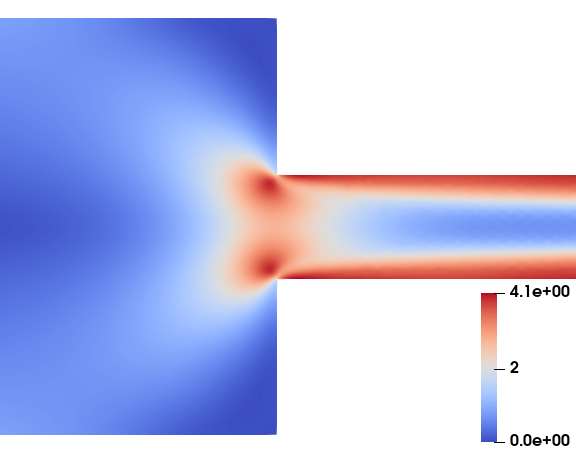}}


\subfloat[Stress diffusion $\sim (\Delta t)^2$ ($Wi = 8$).]{\label{fig:stress_e}
    \includegraphics[width=0.43\textwidth]{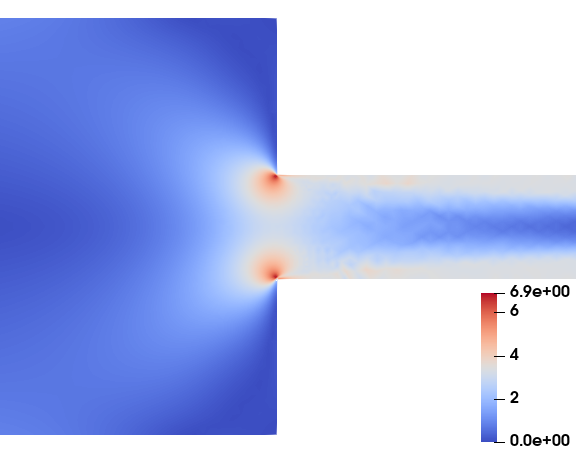}}
\hfill
\subfloat[No stress diffusion ($Wi = 8$).]{\label{fig:stress_f}
    \includegraphics[width=0.43\textwidth]{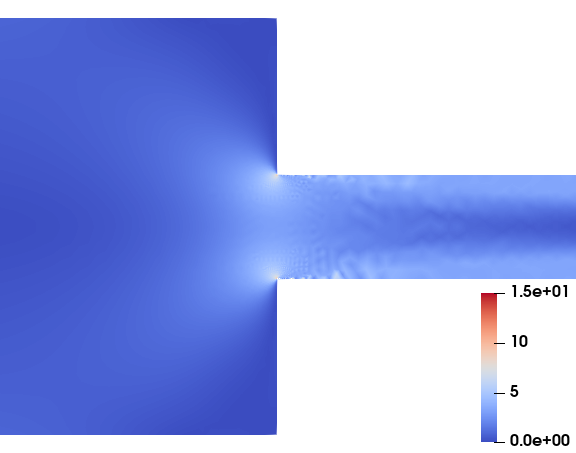}}
    
\caption{
Zoomed view of the contraction plane ($x_1=0$) showing the magnitude of the elastic stress $\mu|\bbF\bbF^\top - \bbI|$.
\textbf{(a)--(c)} Evolution with increasing Weissenberg number using the standard scheme. As $\textit{Wi}\in\{0.1, 1, 8\}$ increases, the stress magnitude near the re-entrant corners grows significantly.
\textbf{(c)--(f)} Comparison of numerical treatments for the highly elastic case ($\textit{Wi}=8$). Compared to the standard scheme (c), both the refined time step size (d) and the stress diffusion scaled by $(\Delta t)^2$ (e) resolve the stress more sharply, leading to more concentrated magnitudes around the corners. In contrast, the non-stabilized case (f) exhibits spurious peaks, indicating numerical instability.
}
\label{fig:stress}
\end{figure}

Figure~\ref{fig:energy_log} illustrates the preservation of the positivity of the determinant, monitored via the logarithmic energy term $-\int_\Omega \frac{\mu}{2} \ln\det(\bbF\bbF^\top)$. This integral is evaluated using eighth-order quadrature rules. For all elastic cases ($\textit{Wi} > 0$), this energy naturally grows over time due to the stress singularities at the sharp re-entrant corners (see Figure~\ref{fig:stress}). In the non-stabilized case (d), the deformation gradient violates the physical constraint $\det(\bbF) > 0$, causing the logarithmic energy to blow up. Both the refined time step in (b) and the quadratic stress diffusion scaling in (c) successfully prevent this failure. While the logarithmic energy curves continue to grow in these stabilized cases, the growth rate is significantly reduced compared to the non-stabilized case (d).

Importantly, these numerical findings do not contradict our theoretical stability and convergence results. Our analysis relies solely on the standard energy functional $\int_\Omega \left( \frac\rho2 |\bv|^2 + \frac\mu2 |\bbF|^2 \right)$, which controls the $L^2$-norms but does not include the logarithmic energy term. Consequently, strict preservation of $\det(\bbF) > 0$ is not guaranteed at the discrete level. As discussed in Section~\ref{sec:positivity_F}, the logarithmic energy is not convex with respect to the general (non-symmetric) tensor $\bbF \in \bbR^{d\times d}$. This lack of convexity prevents the derivation of a discrete energy law that would enforce positivity of the discrete counterpart of $\det(\bbF)$.

Figure~\ref{fig:stress} shows a zoomed view of the contraction plane ($x_1=0$), visualizing the magnitude of the elastic stress $\mu|\bbF\bbF^\top - \bbI|$. First, we examine the physical evolution using the standard scheme in panels (a)--(c). As the Weissenberg number increases ($\textit{Wi}\in\{0.1, 1, 8\}$), the magnitude of the stress accumulates significantly near the re-entrant corners, which is consistent with the expected stress singularity in these regions. Next, we compare the different numerical treatments for the high elasticity case ($\textit{Wi}=8$) in panels (c)--(f). Compared to the standard scheme (c), both the refined time step size (d) and the stress diffusion scaled by $(\Delta t)^2$ (e) resolve the stress more sharply. In these cases, the stress is concentrated directly at the corners instead of being smeared out. In contrast, the non-stabilized simulation (f) shows spurious stress peaks, which indicates the onset of numerical instability.
This happens because, without stress diffusion, the viscoelastic equation lacks any damping mechanism at high $\textit{Wi}$. Near the two re-entrant corners, the velocity gradient $\nabla \bv$ becomes extremely large and stays high in the surrounding region. This strong gradient drives the elastic stress to grow without bound, making the numerical scheme unstable downstream in the narrow part of the channel.

\begin{figure}
\centering
\subfloat[Standard scheme ($\textit{Wi} = 0.1$).]{\label{fig:velocity_a}
    \includegraphics[width=0.43\textwidth]{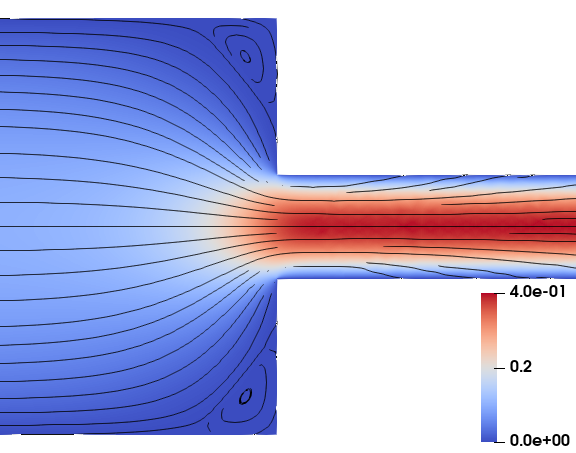}}
\hfill
\subfloat[Standard scheme ($\textit{Wi} = 1$).]{\label{fig:velocity_b}
    \includegraphics[width=0.43\textwidth]{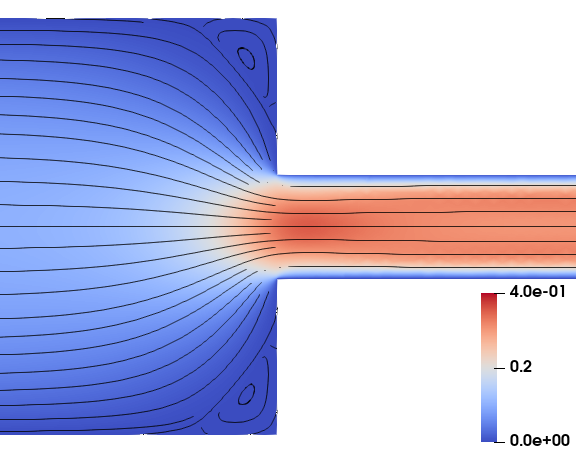}}


\subfloat[Standard scheme ($\textit{Wi} = 8$).]{\label{fig:velocity_c}
    \includegraphics[width=0.43\textwidth]{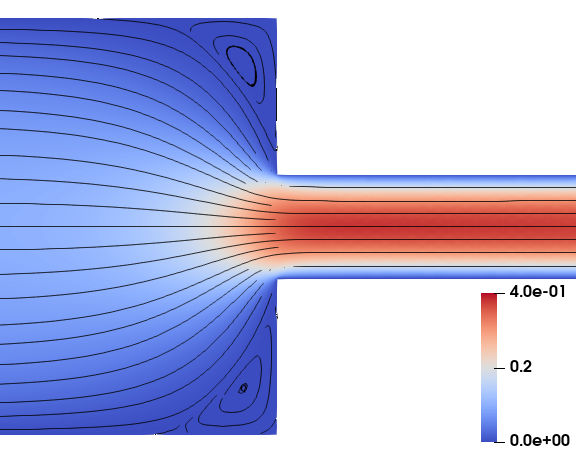}}
\hfill
\subfloat[Smaller time step size ($\textit{Wi} = 8$).]{\label{fig:velocity_d}
    \includegraphics[width=0.43\textwidth]{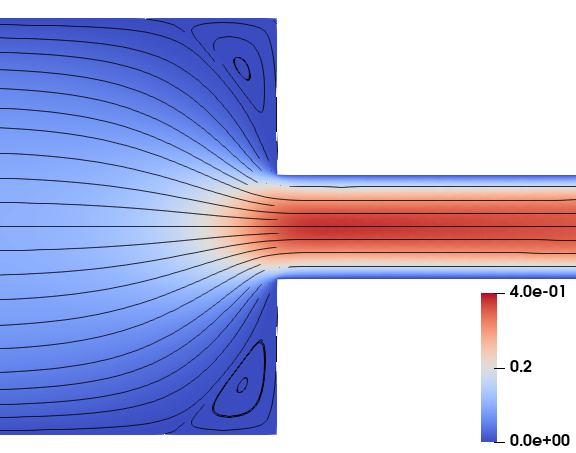}}


\subfloat[Stress diffusion $\sim (\Delta t)^2$ ($Wi = 8$).]{\label{fig:velocity_e}
    \includegraphics[width=0.43\textwidth]{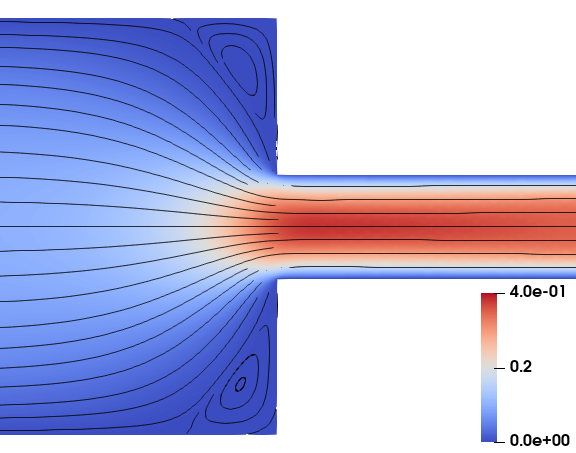}}
\hfill
\subfloat[No stress diffusion ($Wi = 8$).]{\label{fig:velocity_f}
    \includegraphics[width=0.43\textwidth]{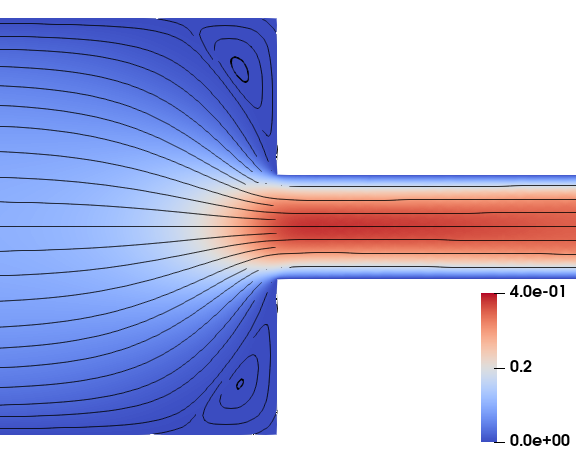}}
    
\caption{
Streamlines of the velocity field near the contraction plane ($x_1=0$). 
\textbf{(a)--(c)} Evolution with increasing $\textit{Wi}\in\{0.1, 1, 8\}$ using the standard scheme. The primary physical change is the growth of the corner vortices as $\textit{Wi}$ increases.
\textbf{(c)--(f)} Comparison of numerical schemes at high elasticity ($Wi=8$). Unlike the sensitive stress fields shown in Figure~\ref{fig:stress}, the velocity profiles remain robust and qualitatively similar across all numerical configurations.
}
\label{fig:velocity}
\end{figure}

Figure~\ref{fig:velocity} displays the corresponding streamlines of the velocity field. The evolution from low to high Weissenberg numbers (panels a--c) highlights the primary physical change, i.e., the growth of the corner vortices as elasticity increases. When comparing the numerical schemes at $\textit{Wi}=8$ (panels c--f), we observe a distinct difference compared to the stress results. Unlike the stress field, which is highly sensitive to the stabilization method, the velocity profiles remain robust. All numerical configurations produce qualitatively similar streamlines, regardless of whether stabilization by stress diffusion is applied.
